\documentclass{article}

\title{Braid groups and quiver mutation}
\author{Joseph Grant and Robert J. Marsh} 
\date{}

\usepackage{latexsym}
\usepackage{verbatim}
\usepackage{amsmath}
\usepackage{amssymb}
\usepackage{mathrsfs}
\usepackage{amsthm}
\usepackage{sfmath}
\usepackage{hyperref}

\usepackage{mathdots}

\input xy
\xyoption{all}
\usepackage{tikz}
\usetikzlibrary{calc,decorations.markings,%
arrows,backgrounds}
\tikzset{->-/.style={decoration={
  markings,
  mark=at position .5 with {\arrow{>}}},postaction={decorate}}}

\usepackage{color} 

\newcommand{\lettera}{a}
\newcommand{\letterb}{b}
\newcommand{\letterc}{c}
\newcommand{\letterd}{d}
\newcommand{\lettere}{e}

\newcommand{\uu}[1]{\widetilde{s}_{#1}}

\newcommand{\be}{\begin{equation}}
\newcommand{\ee}{\end{equation}}

\newcommand{\fs}{\text{.}}

\newcommand{\da}{\text{-}}

\newcommand{\Hom}{\operatorname{Hom}\nolimits}

\newcommand{\dert}{\otimes^{\textbf{L}}} 

\newcommand{\Aut}{\operatorname{Aut}\nolimits}

\newcommand{\cone}{\operatorname{cone}\nolimits}

\newcommand{\Z}{\mathbb{Z}}
\newcommand{\F}{\mathbb{F}}

\newcommand{\id}{\text{id}}

\newcommand{\D}{\operatorname{D}\nolimits}

\newcommand{\T}{\mathcal{T}}
\newcommand{\gen}[1]{\langle#1\rangle}

\newcommand{\into}{\hookrightarrow}
\newcommand{\abs}[1]{\left|#1\right|}

\newcommand{\arr}[1]{\stackrel{#1}{\to}}

\newcommand{\K}{\operatorname{K}\nolimits}

\newcommand{\im}{\operatorname{im}\nolimits}

\newcommand{\op}{{\operatorname{op}\nolimits}}

\newcommand*{\rom}[1]{\expandafter\@slowromancap\romannumeral #1@}
\newcommand{\inv}{^{-1}}

\newcommand{\p}{\operatorname{\textbf{p}}}

\newcommand{\Dfd}{\operatorname{D_{fd}}\nolimits}
\newcommand{\Mod}{\operatorname{-Mod}}
\newcommand{\rMod}{\operatorname{Mod-}}
\newcommand{\biMod}{\operatorname{-Mod-}}
\newcommand{\ev}{\operatorname{ev}}

\newcommand{\cofib}{\operatorname{-cofib}\nolimits}

\newcommand{\st}{\mid}

\newtheorem{thm}{Theorem}[section]
\newtheorem{lem}[thm]{Lemma}
\newtheorem{prop}[thm]{Proposition}

{\theoremstyle{definition} \newtheorem{defn}[thm]{Definition}
\newtheorem{rmk}[thm]{Remark}
}
\newenvironment{pf}{Proof: \em}{ \hfill $\Box$}

\newtheorem*{thma}{Theorem A}
\newtheorem*{thmb}{Theorem B}
\newtheorem*{thmc}{Theorem C}

\newcommand{\vertbowtie}{\rotatebox[origin=c]{90}{$\bowtie$}}

\begin{document}

\maketitle

\begin{abstract}
We describe presentations of braid groups of type $ADE$ and show how these
presentations are compatible with mutation of quivers.
In types $A$ and $D$ these presentations can be understood geometrically using triangulated surfaces.
We then give a categorical interpretation of the presentations, with the new
generators acting as spherical twists at simple modules on derived categories
of Ginzburg dg-algebras of quivers with potential.

\end{abstract}

\tableofcontents

\setlength{\parindent}{0pt}
\setlength{\parskip}{1em plus 0.5ex minus 0.2ex}

\section{Introduction}

Braid groups are fundamental objects in mathematics. Although they are of a
topological and geometric nature, they have an algebraic interpretation:
a simple presentation by generators and relations which is just based on
adjacency of integers~\cite{artin}.
This can be encoded in a line graph, and from there one can generalize to
define a group from any finite graph, known as the Artin braid group.

The most well-known groups defined from graphs are the Coxeter groups
(we restrict to the simply laced cases, for simplicity).
These are closely related to Artin braid groups: each Coxeter group
is a quotient of a corresponding Artin braid group in a natural way.
In particular, the symmetric group on $n$ letters is a quotient of the
classical braid group on $n$ strands. Coxeter groups naturally split into
two distinct classes: those of finite type, corresponding to the Dynkin
diagrams of type ADE, and those of infinite type.
Although all Artin braid groups are infinite, the Artin braid groups of
Dynkin type have a different character to those not of Dynkin type,
and are known as Artin groups of `finite type'.

This dichotomy also arises in another area of mathematics
which has generated a lot of interest in the recent years: cluster algebras.
In this theory, there is a notion of finite-type cluster algebras, which again correspond to the Dynkin diagrams~\cite{fz2}.
Cluster algebras are specified by a directed graph, known as a quiver, together with other information.  A key ingredient in the definition is the notion of mutation, which changes the arrows in a quiver in a non-obvious manner which generalizes reflection at a source or sink. Barot and Marsh~\cite{bm} have
given new presentations of Coxeter groups of finite type based on quivers
obtained from Dynkin diagrams under finite sequences of mutations.
Our first result generalizes this to braid groups:

\begin{thma} (Theorem~\ref{t:groupinvariance})
Let $Q$ be a quiver, with vertices $1,2,\ldots ,n$, obtained from a Dynkin quiver by a finite sequence of mutations. Let $B_Q$ be the group with generators $s_1,s_2,\ldots ,s_n$,
subject to the relations:
\begin{itemize}
\item[(a)] $s_is_j=s_js_i$ if there is no arrow between $i$ and $j$ (in either
direction);
\item[(b)] $s_is_js_i=s_js_is_j$ if there is an arrow between $i$ and $j$ (in either direction);
\item[(c)] $s_{i_1}s_{i_2}\cdots s_{i_n}s_{i_1}\cdots s_{i_{n-2}}=s_{i_2}s_{i_3}\cdots s_{i_n}s_{i_1}\cdots s_{i_{n-1}}=\cdots $,
whenever $i_1\rightarrow i_2\rightarrow \cdots \rightarrow i_n\rightarrow i_1$
is a chordless cycle in $Q$.
\end{itemize}
Then $B_Q$ is isomorphic to the Artin braid group of the same Dynkin type
as $Q$.
\end{thma}

We prove our result via isomorphisms between abstractly defined groups, which can be thought of as mutations of groups, even though the resulting groups are isomorphic. The Artin group presentations we obtain induce presentations of the corresponding Coxeter groups which are distinct from those in~\cite{bm}; we also give a compatibility result which shows the relationship between the two
presentations. Why have we chosen to use presentations which don't agree with the earlier work?  This is explained in the following two sections of the paper, as we now detail.

Certain cluster algebras can be understood using pictures.
A (tagged) triangulation of a Riemann surface with marked points on its boundary defines a quiver~\cite{fz1,fst08}; see also~\cite{ccs}.
Then mutation of the quiver has a natural interpretation in terms of swapping one diagonal of a given quadrilateral for the other.  So these cluster algebras have a topological interpretation. In particular, such descriptions are available for the infinite families of Dynkin type.  A natural question is: can we understand the generators above, and the isomorphisms corresponding to mutations, in terms of the geometry of the surface? The answer is yes:

\begin{thmb} (Theorem~\ref{t:braidisomorphism})
Let $\Delta$ be a Dynkin diagram of type $A_n$ or type $D_n$. In the former
case, let $(X,M)$ be a disk with $n+3$ marked points on its boundary. In the
latter case, let $(X,M)$ be a disk with $n$ marked points on its boundary and one
marked point in its interior, taken to be a cone point of order $2$ (so $X$ is
an orbifold in this case).

Let $\T$ be a tagged triangulation of $(X,M)$. Let $G_{\T}$ be the graph
in $(X,M)$ dual to $\T$. For each vertex $i$ of the quiver $Q_{\T}$ associated
to $\T$ as in~\cite{fz1,fst08},
let $\sigma_i$ be the braid of $(X,M)$ associated to the edge of
$G_{\T}$ crossing the tagged arc in $\T$ corresponding to $i$ (see Definition~\ref{d:pathbraid}). Then there
is an isomorphism between the subgroup $H_{\T}$ of the braid group generated by the $\sigma_i$ and the group $B_Q$ defined above, taking $\sigma_i$
to $s_i$.

Furthermore, in type $A_n$, $H_{\T}$ coincides with the braid group of
$(X,M)$, while in type $D_n$, $H_{\T}$ is of index two in the braid group of
$(X,M)$.
\end{thmb}

As well as the original combinatorial and commutative algebraic approach to cluster algebras, and the geometric approach described above, there is a third approach which has proved very powerful: the representation theoretic approach~\cite{bmrrt,ccs}. This approach uses finite dimensional (noncommutative) algebras and ideas from categorification to better understand cluster algebras, and has received intense study.  Braid groups also appear in representation theory and categorification~\cite{rz,st}: in many important situations there are actions of braid groups on derived categories via spherical twists.  One example of this is given by certain derived categories of differential graded algebras~\cite{g-cy,ky} which are known to cover the categories appearing in the representation theoretic approach to cluster algebras~\cite{amiot}.
One might hope that these categorical braid group actions are related to our presentations of braid groups, and we show that this is indeed the case.

First, we make a connection between the categorical and the geometric
situations.  The relevant differential graded algebras are defined by use of a 
quiver together with a formal sum of cycles in that quiver known as a
potential~\cite{g-cy}. Mutation of quivers of potential has been defined~\cite{dwz1} and, in the situations where our cluster algebra comes from a Riemann surface, the mutation of potentials also has a geometric interpretation~\cite{LF1}.  Relying heavily on results of Labardini-Fragoso~\cite{LF1,LF4}, we observe that the potential defined on mutation-Dynkin quivers according to the geometric procedure is
equivalent to the `obvious' potential that one might guess (Proposition~\ref{prop:pot-is-sum-cyc}). So, while the potential is important, it is in fact entirely determined by the quiver in types A and D. Note that this result could also be proved relatively easily via a direct calculation.

Next we show that we do indeed obtain an action of the groups $B_Q$
(defined using mutation-Dynkin quivers) on derived categories of Ginzburg
differential graded algebras in which the generators act via spherical twists.
After setting up all the technical machinery correctly, the main difficulty in proving this is to check that the mutation procedure for the groups $B_Q$, which relates the group associated to a quiver to the group associated to a mutated quiver, actually lifts to the categorical setting as a natural isomorphism of functors.  We do this, using important results of Keller and Yang~\cite{ky}. From here, we can use the earlier theory developed here
to show that the generators $s_i$ of finite type Artin braid groups from
Theorem A can be viewed as derived autoequivalences:

\begin{thmc} (Theorem~\ref{thm:mutation-action})
Let $(Q,W)$ be a mutation-Dynkin quiver with potential of type $ADE$, and
let $\Gamma_{Q,W}$ be the corresponding Ginzburg differential graded algebra.
Let $\Dfd(\Gamma_{Q,W})$ denote the full subcategory of the
derived category $\D(\Gamma_{Q,W})$ on objects with finite-dimensional total homology. Then there is a group homomorphism
\begin{align*}
B_Q &\to \Aut\Dfd(\Gamma_{Q,W})\\
s_i &\mapsto F_i
\end{align*}
sending the group generator associated to the vertex $i$ of $Q$ to the spherical twist $F_i$ at the simple $\Gamma_{Q,W}$-module $S_i$.
\end{thmc}

Since we started work on this project, we have become
aware of independent work by other authors.
A. King and Y. Qiu have a related project, and were aware of 
the new relations between spherical twists and a topological
interpretation of the spherical twist group;
see \cite{qiu}, particularly Section 10.1. In particular, an 
independent proof of a version of Theorem~\ref{thm:groupiso} 
in types $A$ and $D$ was announced in~\cite{qiu}. A key 
difference in our approach is the use of an orbifold with 
cone point of degree two in type $D$.
In~\cite[\S2.2]{nag-triang},
K. Nagao refers to an action of the mapping class group
of a marked surface on the derived category of a Ginzburg
dg-algebra associated to a triangulation.

Since we released the first draft of this article, the preprint~\cite{HHLP} has 
appeared, where the authors give a presentation (different from the one given 
here) of the Artin braid group for each diagram of finite type (in the cluster-theoretic sense). This includes the non-simply-laced cases (not considered here) 
but does not include a topological or categorical interpretation.

\section{Presentations of braid groups}

\subsection{Braid groups}
Let $\Delta$ be a graph of $ADE$ Dynkin type, i.e., $\Delta$ is a graph of type $A_n$ for $n\geq1$, $D_n$ for $n\geq4$, $E_6$, $E_7$, or $E_8$.  
$$
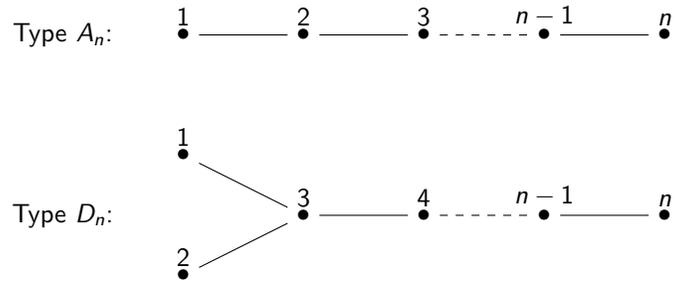
\begin{tikzpicture}[scale=0.8,baseline=(bb.base),
  quivarrow/.style={black, -latex}] 
\path (0,0) node (bb) {}; 

\node at (-2,0) {Type $A_n$:};

\node (1) at (0,0) {$\bullet$};
\node [above] at (1) {$1$}; 

\node (2) at (2,0) {$\bullet$};
\node [above] at (2) {$2$}; 

\node (3) at (4,0) {$\bullet$};
\node [above] at (3) {$3$}; 

\node (n1) at (6,0) {$\bullet$};
\node [above] at (n1) {$n-1$}; 

\node (n) at (8,0) {$\bullet$};
\node [above] at (n) {$n$}; 

\draw (1) -- (2) -- (3);
\draw [dashed] (3) to (n1);
\draw (n1) -- (n);

\node at (-2,-3) {Type $D_n$:};

\node (d1) at (0,-2) {$\bullet$};
\node [above] at (d1) {$1$}; 

\node (d2) at (0,-4) {$\bullet$};
\node [above] at (d2) {$2$}; 

\node (d3) at (2,-3) {$\bullet$};
\node [above] at (d3) {$3$}; 

\node (d4) at (4,-3) {$\bullet$};
\node [above] at (d4) {$4$}; 

\node (dn1) at (6,-3) {$\bullet$};
\node [above] at (dn1) {$n-1$}; 

\node (dn) at (8,-3) {$\bullet$};
\node [above] at (dn) {$n$}; 

\draw (d1) -- (d3) -- (d4);
\draw (d2) -- (d3);
\draw [dashed] (d4) to (dn1);
\draw (dn1) -- (dn);

\end{tikzpicture}
$$
In particular, $\Delta$ has no double edges or cycles.
Let $I$ be the set of vertices of $\Delta$.  We can associate a group $B_\Delta$ to $\Delta$, which we call \emph{the braid group of $\Delta$}.  It has a distinguished set of generators $S_\Delta=\{s_i\}_{i\in I}$, and the relations depend on whether or not two vertices are connected by an edge.  They are as follows:
\begin{enumerate}
\item\label{reln:comm} $s_is_j=s_js_i$ if $\xymatrix{ \stackrel i \bullet  & \stackrel j \bullet }$;
\item\label{reln:reid3} $s_is_js_i=s_js_is_j$ if $\xymatrix{ \stackrel i \bullet \ar@{-}[r] & \stackrel j \bullet }$;
\end{enumerate}
If $\Delta$ is of type $A_n$ then we recover the ``usual'' braid group, sometimes denoted $B_{n+1}$.  Its generators can be visualized as follows:
$$
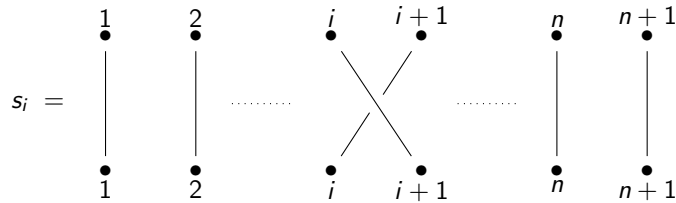
\begin{tikzpicture}[scale=0.6,baseline=(bb.base),
  quivarrow/.style={black, -latex}] 
\path (0,0) node (bb) {}; 

\node at (-1.5,-1.5) {$s_i \; = $};

\node (1) at (0,0) {$\bullet$};
\node [above] at (1) {$1$}; 

\node (2) at (2,0) {$\bullet$};
\node [above] at (2) {$2$}; 

\node (i) at (5,0) {$\bullet$};
\node [above] at (i) {$i$}; 

\node (i1) at (7,0) {$\bullet$};
\node [above] at (i1) {$i+1$}; 

\node (n) at (10,0) {$\bullet$};
\node [above] at (n) {$n$}; 

\node (n1) at (12,0) {$\bullet$};
\node [above] at (n1) {$n+1$};

\node (b1) at (0,-3) {$\bullet$};
\node [below] at (b1) {$1$}; 

\node (b2) at (2,-3) {$\bullet$};
\node [below] at (b2) {$2$}; 

\node (bi) at (5,-3) {$\bullet$};
\node [below] at (bi) {$i$}; 

\node (bi1) at (7,-3) {$\bullet$};
\node [below] at (bi1) {$i+1$}; 

\node (bn) at (10,-3) {$\bullet$};
\node [below] at (bn) {$n$}; 

\node (bn1) at (12,-3) {$\bullet$};
\node [below] at (bn1) {$n+1$}; 

\draw (1) -- (b1);
\draw (2) -- (b2);
\draw (n) -- (bn);
\draw (n1) -- (bn1);

\draw [dotted] (2.8,-1.5) -- (4.2,-1.5);
\draw [dotted] (7.8,-1.5) -- (9.2,-1.5);

\draw (i) -- (bi1);
\draw (i1) -- (6.15,-1.3);
\draw (5.85,-1.7) -- (bi);

\end{tikzpicture}
$$
and the relations of type \ref{reln:comm} record the fact that crossings of far-apart adjacent pairs of strings commute, while relations of type \ref{reln:reid3} record a Reidemeister 3 move.

If we also impose the relation that $s_i^2=1$ for all $i\in I$ then we recover the Coxeter group of type $\Delta$.  More information on Coxeter groups and braid groups can be found in \cite{hum} and \cite{kt-braid}.

\subsection{Mutation of quivers}

A quiver is just a directed graph.  Throughout this article we will only work with quivers with finitely many vertices and finitely many arrows that have no loops or oriented $2$-cycles.  For a given quiver $Q$, we again denote its set of vertices by $I$.

There is a procedure to obtain one quiver from another, called \emph{quiver mutation}, due to Fomin and Zelevinsky \cite[\S4]{fz1}.  Fix $Q$ and let $k\in I$.  Then we obtain the mutated quiver $\mu_k(Q)$ as follows:
\begin{enumerate}
\item for each pair of arrows $i\to k\to j$ through $k$, add a formal composite $i\to j$;
\item reverse the orientation of all arrows incident with $k$;
\item remove a maximal set of $2$-cycles (we may have created $2$-cycles in the previous two steps).
\end{enumerate}
It is a basic but important observation that quiver mutation does not change the set of vertices.  One can also check that mutation is an involution.

We call a cycle in an unoriented graph (or in the underlying unoriented graph of a quiver) \emph{chordless} if the full subgraph on the vertices of the cycle contains no edges which are not part of the cycle.  We call a quiver \emph{Dynkin} if its underlying unoriented graph is a Dynkin graph of type $ADE$, and \emph{mutation-Dynkin} if it can be obtained by mutating a Dynkin quiver finitely many times.  By a theorem of Fomin and Zelevinsky \cite[Thm.\ 1.4]{fz2}, there are only finitely many quivers that can be obtained by mutating a given Dynkin quiver.

The following fact will be useful to us.
\begin{prop}[Fomin-Zelevinsky]\label{prop:Dynkinquiv}
In any mutation-Dynkin quiver, there are no double arrows and all chordless cycles are oriented.

\begin{pf}
By \cite[Theorem 1.8]{fz2}, the entries in the corresponding exchange matrix $B$ satisfy
$\abs{B_{xy}B_{yx}}\leq3$ for all $x,y$ (known as being $2$-finite).
Hence there cannot be any double arrows in the quiver.

Now let $Q$ be a mutation-Dynkin quiver and $C$ a chordless cycle in $Q$.
Then, since $Q$ is 2-finite, 
so is $C$. 
By \cite[Proposition 9.7]{fz2}, $C$ must be an oriented cycle.
\end{pf}
\end{prop}

\subsection{Groups from quivers}\label{ss:define-B_Q}

Let $Q$ be a mutation-Dynkin quiver.

\begin{defn}
From the quiver $Q$ with vertex set $I$, we define the group $B_Q$ as follows: it has a distinguished generating set $S_Q=\{s_i\}_{i\in I}$ and the relations given by:
\begin{enumerate}
\item\label{reln:q-comm} $s_is_j=s_js_i$ if $\xymatrix{ \stackrel i \bullet  & \stackrel j \bullet }$;
\item\label{reln:q-reid3} $s_is_js_i=s_js_is_j$ if $\xymatrix{ \stackrel i \bullet \ar[r] & \stackrel j \bullet }$ or $\xymatrix{ \stackrel i \bullet  & \stackrel j \bullet \ar[l]}$;
\item\label{reln:q-cycle} if we have an oriented chordless $n$-cycle
$$
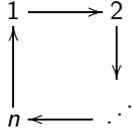
\xymatrix{
1\ar[r] &2\ar[d]\\
n\ar[u] & 
{}\iddots\ar[l]
}$$
then
\begin{align*}
s_1s_2\ldots s_ns_1s_2\ldots s_{n-2} &= s_2s_3\ldots s_ns_1\ldots s_{n-1}\\
 &=\cdots\\
 &=s_ns_1s_2\ldots  s_ns_1s_2\ldots s_{n-3}
\end{align*}
\end{enumerate}
\end{defn}

\begin{rmk}\label{r:artin}
If $Q$ is a Dynkin quiver, then $B_Q$ is (isomorphic to) the Artin braid group of the corresponding Dynkin type.
\end{rmk}

This presentation is symmetric but not minimal:
\begin{lem}\label{lem:one-implies-all}
For each single chordless $n$-cycle, in the presence of the relations of type \ref{reln:q-comm} and \ref{reln:q-reid3}, any one of the relations of type \ref{reln:q-cycle} implies all the others.

\begin{pf}
It is enough to show that if the relation
\begin{equation}
\label{e:firstrel}
s_1s_2\cdots s_ns_1s_2\cdots s_{n-2}
=
s_2s_3\cdots s_ns_1s_2\cdots s_{n-1}
\end{equation}

holds then

\begin{equation*}
s_1s_2\cdots s_ns_1s_2\cdots s_{n-2}
=
s_3s_4\cdots s_ns_1s_2\cdots s_n.
\end{equation*}
So, we assume that~\eqref{e:firstrel}
holds. Then we have:
\begin{equation*}
s_2^{-1}s_1s_2\cdots s_ns_1s_2\cdots s_{n-2}s_n
=
s_3\cdots s_ns_1s_2\cdots s_{n-1}s_n.
\end{equation*}
The left hand side can be rewritten,
using relations of type (i) and (ii),
as:
\begin{equation*}
\begin{split}
s_2^{-1}s_1s_2\cdots s_ns_1s_2\cdots s_{n-2}s_n &=
s_1s_2s_1^{-1}s_3\cdots s_ns_1s_2\cdots s_{n-2}s_n \\
&=s_1s_2s_3\cdots s_{n-1}s_1^{-1}s_ns_1s_2\cdots s_{n-2}s_n \\
&=s_1s_2\cdots s_{n-1}s_ns_1s_n^{-1}s_2
\cdots s_{n-2}s_n \\
&=s_1s_2\cdots s_ns_1s_2\cdots s_{n-2},
\end{split}
\end{equation*}
and the result follows.
\end{pf}
\end{lem}

Though the relations look different, by taking an appropriate quotient we can obtain the groups defined by Barot and Marsh \cite{bm} directly:

\begin{lem}
If we also impose the relations $s_i^2=1$ for all $i\in I$, then the group $B_Q$ becomes isomorphic to the group $\Gamma_{U(Q)}$ defined in \cite[Section 3]{bm}, where $U(Q)$ is the underlying graph of $Q$.

\begin{pf}
As our definition gives the usual definition of the braid group for a Dynkin quiver, this will follow from results in \cite{bm} and the results below on how our groups change with quiver mutation.  But it is straightforward to give a direct proof, so we do so.

We need to show that, in the presence of relations \ref{reln:q-comm}, \ref{reln:q-reid3}, and $s_i^2=1$ for all $i\in I$, our extra relation \ref{reln:q-cycle} holds if and only if the relation
$$(s_1s_2\ldots s_{n-1}s_ns_{n-1}\ldots s_2)^2=1$$
and its rotations hold for each $n$-cycle $1\to 2\to \cdots\to n\to 1$.  By symmetry, it is enough to check that the relation above is equivalent to $s_1s_2\ldots s_ns_1s_2\ldots s_{n-2} = s_2s_3\ldots s_ns_1\ldots s_{n-1}$.

Using that $s_i=s_i^{-1}$, we see that our our relation is equivalent to
$$s_1s_2\ldots s_ns_1s_2\ldots s_{n-2}s_{n-1}s_{n-2}\ldots s_1s_n\ldots s_3 s_2=1.$$
Multiplying out the Barot-Marsh relation, we see that it is equivalent to
$$s_1s_2\ldots s_{n-1}s_ns_{n-1}\ldots s_2s_1s_2\ldots s_{n-1}s_ns_{n-1}\ldots s_2=1\fs$$
Cancelling out $n$ terms on the left and $n-1$ terms on the right of these two expressions, it just remains to show that
$$s_1s_2\ldots s_{n-2}s_{n-1}s_{n-2}\ldots s_2s_1=s_{n-1}s_{n-2}\ldots s_2s_1s_2\ldots s_{n-2}s_{n-1}\fs$$
As there is an arrow $i\to i+1$ for each $i$ and the cycle is chordless,
the symmetric group on $n$ letters maps onto the subgroup generated by $s_1,\ldots,s_{n-1}$ with the transposition which swaps $i$ and
$i+1$ being sent to $s_i$. It is easy to see that the corresponding relation holds in the symmetric group, with both sides of the equation representing the transposition which swaps $1$ and $n$.
\end{pf}
\end{lem}
We will justify our choice of relations in Remark~\ref{r:singletwist}.

\subsection{Mutation of groups}
Let $B_Q$ be the group associated to the mutation-Dynkin quiver $Q$, as above, and let $k$ be a vertex of $Q$.  Denote $\mu_k(Q)$ by $Q'$. 
Our aim in
this section is to show that $B_Q$ is
isomorphic to $B_{Q'}$. We will do this
by using a group homomorphism $\varphi_k:B_Q\rightarrow B_{Q'}$ defined using a formula which lifts the formula used in~\cite[\S5]{bm}.

The following lemma follows from results
in~\cite{fz2} (see~\cite[\S2]{bm}).

\begin{figure}
$$
\begin{tikzpicture}[scale=1,baseline=(bb.base),
  quivarrow/.style={black, -latex}] 
\path (0,0) node (bb) {}; 


\node at (-0.9,1.4) {(a)};

\node (i1) at (0,0) {$\circ$};
\node [below] at (i1) {\small {$i$}}; 

\node (j1) at (1,0) {$\circ$};
\node [below] at (j1) {\small {$j$}};

\node (k1) at (0.5,1) {$\bullet$};
\node[above] at (k1) {\small {$k$}};

\draw [quivarrow, shorten <=-1pt, shorten >=-1pt] (i1) -- (k1);
\draw [quivarrow, shorten <=-1pt, shorten >=-1pt] (j1) -- (k1);

\draw [<->] (1.8,0.5) to node[above] {$\mu_k$} (2.3,0.5);

\begin{scope}[shift={(3.1,0)}]

\node (i2) at (0,0) {$\circ$};
\node [below] at (i2) {\small {$i$}}; 

\node (j2) at (1,0) {$\circ$};
\node [below] at (j2) {\small {$j$}};

\node (k2) at (0.5,1) {$\bullet$};
\node[above] at (k2) {\small {$k$}};

\draw [quivarrow, shorten <=-1pt, shorten >=-1pt] (k2) -- (i2);
\draw [quivarrow, shorten <=-1pt, shorten >=-1pt] (k2) -- (j2);

\end{scope}


\begin{scope}[shift={(7,0)}]

\node at (-0.9,1.4) {(b)};

\node (i1) at (0,0) {$\circ$};
\node [below] at (i1) {\small {$i$}}; 

\node (j1) at (1,0) {$\circ$};
\node [below] at (j1) {\small {$j$}};

\node (k1) at (0.5,1) {$\bullet$};
\node[above] at (k1) {\small {$k$}};

\draw [quivarrow, shorten <=-1pt, shorten >=-1pt] (i1) -- (k1);
\draw [quivarrow, shorten <=-1pt, shorten >=-1pt] (k1) -- (j1);

\draw [<->] (1.8,0.5) to node[above] {$\mu_k$} (2.3,0.5);

\begin{scope}[shift={(3.1,0)}]

\node (i2) at (0,0) {$\circ$};
\node [below] at (i2) {\small {$i$}}; 

\node (j2) at (1,0) {$\circ$};
\node [below] at (j2) {\small {$j$}};

\node (k2) at (0.5,1) {$\bullet$};
\node[above] at (k2) {\small {$k$}};

\draw [quivarrow, shorten <=-1pt, shorten >=-1pt] (k2) -- (i2);
\draw [quivarrow, shorten <=-1pt, shorten >=-1pt] (j2) -- (k2);
\draw [quivarrow, shorten <=-1pt, shorten >=-1pt] (i2) -- (j2);

\end{scope}
\end{scope}

\end{tikzpicture}
$$
\caption{Mutation of a quiver of mutation
Dynkin type.}
\label{f:localmutation}
\end{figure}
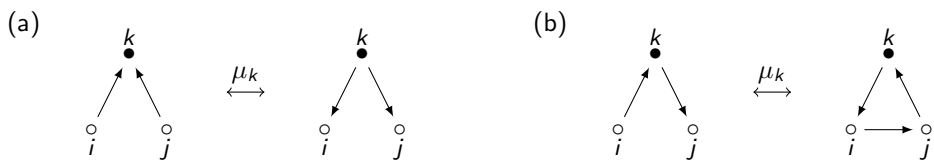

\begin{lem} \label{l:localmutation}
Let $Q$ be a quiver of mutation-Dynkin
type, and fix a vertex $k$ of $Q$. Suppose
that $k$ has two neighbouring vertices.
Then the possibilities for the induced
subquiver of $Q$ containing vertex $k$
and its neighbours are shown in Figure~\ref{f:localmutation}. The effect of
mutation is shown in each case.
\end{lem}

The following lemma follows from~\cite[Lemma 2.5]{bm}.

\begin{lem} \label{l:cycles}
Let $Q$ be a quiver of mutation-Dynkin type,
and fix a vertex $k$ of $Q$.
Let $C$ be an oriented cycle in $Q$. Then $C$
is one of the following. In each case we
indicate what happens locally under mutation at
$k$.
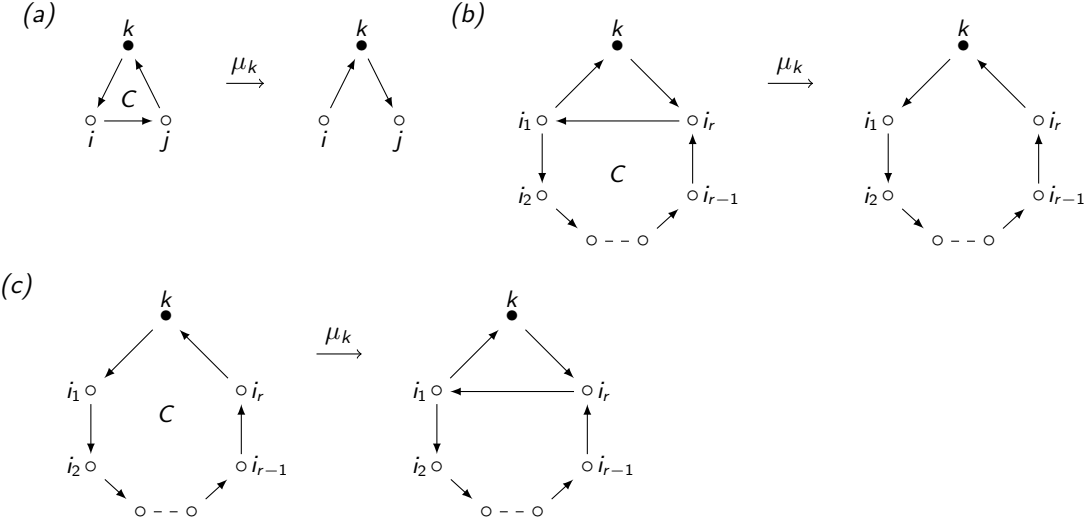
$$
\begin{tikzpicture}[scale=1,baseline=(bb.base),
  quivarrow/.style={black, -latex}] 
\path (0,0) node (bb) {}; 

\node at (-0.7,1.4) {(a)};

\node (i1) at (0,0) {$\circ$};
\node [below] at (i1) {\small {$i$}}; 

\node (j1) at (1,0) {$\circ$};
\node [below] at (j1) {\small {$j$}};

\node (k1) at (0.5,1) {$\bullet$};
\node[above] at (k1) {\small {$k$}};

\node (C) at (0.5,0.3) {\small $C$};

\draw [quivarrow, shorten <=-1pt, shorten >=-1pt] (k1) -- (i1);
\draw [quivarrow, shorten <=-1pt, shorten >=-1pt] (i1) -- (j1);
\draw [quivarrow, shorten <=-1pt, shorten >=-1pt] (j1) -- (k1);

\draw [->] (1.8,0.5) to node[above] {$\mu_k$} (2.3,0.5);

\begin{scope}[shift={(3.1,0)}]

\node (i2) at (0,0) {$\circ$};
\node [below] at (i2) {\small {$i$}}; 

\node (j2) at (1,0) {$\circ$};
\node [below] at (j2) {\small {$j$}};

\node (k2) at (0.5,1) {$\bullet$};
\node[above] at (k2) {\small {$k$}};

\draw [quivarrow, shorten <=-1pt, shorten >=-1pt] (i2) -- (k2);
\draw [quivarrow, shorten <=-1pt, shorten >=-1pt] (k2) -- (j2);
\end{scope}

\begin{scope}[shift={(6,0)},xscale=2]

\node at (-0.5,1.4) {(b)};

\node (i1A) at (0,0) {$\circ$};
\node [left] at (i1A) {\small {$i_1$}}; 

\node (irA) at (1,0) {$\circ$};
\node [right] at (irA) {\small {$i_r$}};

\node (kA) at (0.5,1) {$\bullet$};
\node[above] at (kA) {\small {$k$}};

\node (i2A) at (0,-1) {$\circ$};
\node [left] at (i2A) {\small {$i_2$}};

\node (i3A) at (0.33,-1.6) {$\circ$};

\node (ir1A) at (1,-1) {$\circ$};
\node [right] at (ir1A) {\small {$i_{r-1}$}};

\node (ir2A) at (0.67,-1.6) {$\circ$};

\node (C) at (0.5,-0.7) {\small $C$};

\draw[quivarrow, shorten <=-1pt, shorten >=-1pt] (i1A) -- (kA);
\draw[quivarrow, shorten <=-1pt, shorten >=-1pt] (irA) -- (i1A);
\draw[quivarrow, shorten <=-1pt, shorten >=-1pt] (kA) -- (irA);
\draw [quivarrow, shorten <=-1pt, shorten >=-1pt] (i1A) -- (i2A);
\draw [quivarrow, shorten <=-1pt, shorten >=-1pt] (i2A) -- (i3A);
\draw [quivarrow, shorten <=-1pt, shorten >=-1pt] (ir2A) -- (ir1A);
\draw [quivarrow, shorten <=-1pt, shorten >=-1pt] (ir1A) -- (irA);
\draw [dashed, shorten <=-1pt, shorten >=-1pt] (i3A) -- (ir2A);

\draw [->] (1.5,0.5) to node[above] {$\mu_k$} (1.8,0.5);

\begin{scope}[shift={(2.3,0)}]

\node (i1B) at (0,0) {$\circ$};
\node [left] at (i1B) {\small {$i_1$}}; 

\node (irB) at (1,0) {$\circ$};
\node [right] at (irB) {\small {$i_r$}};

\node (kB) at (0.5,1) {$\bullet$};
\node[above] at (kB) {\small {$k$}};

\node (i2B) at (0,-1) {$\circ$};
\node [left] at (i2B) {\small {$i_2$}};

\node (i3B) at (0.33,-1.6) {$\circ$};

\node (ir1B) at (1,-1) {$\circ$};
\node [right] at (ir1B) {\small {$i_{r-1}$}};

\node (ir2B) at (0.67,-1.6) {$\circ$};


\draw [quivarrow, shorten <=-1pt, shorten >=-1pt] (kB) -- (i1B);
\draw [quivarrow, shorten <=-1pt, shorten >=-1pt] (irB) -- (kB);
\draw [quivarrow, shorten <=-1pt, shorten >=-1pt] (ir1B) -- (irB);
\draw [quivarrow, shorten <=-1pt, shorten >=-1pt] (ir2B) -- (ir1B);
\draw [quivarrow, shorten <=-1pt, shorten >=-1pt] (i2B) -- (i3B);
\draw [quivarrow, shorten <=-1pt, shorten >=-1pt] (i1B) -- (i2B);
\draw [dashed, shorten <=-1pt, shorten >=-1pt] (i3B) -- (ir2B);

\end{scope}

\end{scope}

\begin{scope}[shift={(0,-3.6)},xscale=2]

\node at (-0.5,1.4) {(c)};

\node (i1D) at (0,0) {$\circ$};
\node [left] at (i1D) {\small {$i_1$}}; 

\node (irD) at (1,0) {$\circ$};
\node [right] at (irD) {\small {$i_r$}};

\node (kD) at (0.5,1) {$\bullet$};
\node[above] at (kD) {\small {$k$}};

\node (i2D) at (0,-1) {$\circ$};
\node [left] at (i2D) {\small {$i_2$}};

\node (i3D) at (0.33,-1.6) {$\circ$};

\node (ir1D) at (1,-1) {$\circ$};
\node [right] at (ir1D) {\small {$i_{r-1}$}};

\node (ir2D) at (0.67,-1.6) {$\circ$};

\node (C) at (0.5,-0.3) {\small $C$};

\draw [quivarrow, shorten <=-1pt, shorten >=-1pt] (kD) -- (i1D);
\draw [quivarrow, shorten <=-1pt, shorten >=-1pt] (irD) -- (kD);
\draw [quivarrow, shorten <=-1pt, shorten >=-1pt] (ir1D) -- (irD);
\draw [quivarrow, shorten <=-1pt, shorten >=-1pt] (ir2D) -- (ir1D);
\draw [quivarrow, shorten <=-1pt, shorten >=-1pt] (i2D) -- (i3D);
\draw [quivarrow, shorten <=-1pt, shorten >=-1pt] (i1D) -- (i2D);
\draw [dashed, shorten <=-1pt, shorten >=-1pt] (i3D) -- (ir2D);

\draw [->] (1.5,0.5) to node[above] {$\mu_k$} (1.8,0.5);

\begin{scope}[shift={(2.3,0)}]

\node (i1C) at (0,0) {$\circ$};
\node [left] at (i1C) {\small {$i_1$}}; 

\node (irC) at (1,0) {$\circ$};
\node [right] at (irC) {\small {$i_r$}};

\node (kC) at (0.5,1) {$\bullet$};
\node[above] at (kC) {\small {$k$}};

\node (i2C) at (0,-1) {$\circ$};
\node [left] at (i2C) {\small {$i_2$}};

\node (i3C) at (0.33,-1.6) {$\circ$};

\node (ir1C) at (1,-1) {$\circ$};
\node [right] at (ir1C) {\small {$i_{r-1}$}};

\node (ir2C) at (0.67,-1.6) {$\circ$};


\draw [quivarrow, shorten <=-1pt, shorten >=-1pt] (i1C) -- (kC);
\draw [quivarrow, shorten <=-1pt, shorten >=-1pt] (irC) -- (i1C);
\draw [quivarrow, shorten <=-1pt, shorten >=-1pt] (kC) -- (irC);
\draw [quivarrow, shorten <=-1pt, shorten >=-1pt] (i1C) -- (i2C);
\draw [quivarrow, shorten <=-1pt, shorten >=-1pt] (i2C) -- (i3C);
\draw [quivarrow, shorten <=-1pt, shorten >=-1pt] (ir2C) -- (ir1C);
\draw [quivarrow, shorten <=-1pt, shorten >=-1pt] (ir1C) -- (irC);
\draw [dashed, shorten <=-1pt, shorten >=-1pt] (i3C) -- (ir2C);

\end{scope}
\end{scope}

\end{tikzpicture}
$$
\begin{itemize}
\item[(d)] An oriented cycle containing exactly
one neighbour of $k$. Mutation at $k$
reverses the arrow between $k$ and its neighbour
in $C$.
\item[(e)] An oriented cycle containing no
neighbours of $k$. Mutation at $k$ does not
affect $C$.
\end{itemize}
\end{lem}

Recall that $B_Q$ is defined using
generators $s_i$ for $i\in I$. We
denote the corresponding generating
set for $B_{Q'}$ by $t_i$, $i\in I$.
Let $F_Q$ be the free group on the generators $s_i$ for $i\in I$.

\begin{defn}\label{def:mut-rule}
Let $\varphi_k:F_Q\rightarrow B_{Q'}$ be
the group homomorphism defined by
$$
\varphi_k(s_i)=
\begin{cases}
t_k t_i t_k^{-1} & \text{if } i\rightarrow k \text{ in } Q; \\
t_i & \text{otherwise.}
\end{cases}
$$
\end{defn}

\begin{prop} \label{p:groupinvariance}
The group homomorphism $\varphi_k$ induces
a group homomorphism (which we also denote
by $\varphi_k$) from $B_Q$ to $B_{Q'}$.

\begin{pf}
Let us write $\uu{i}=\varphi_k(s_i)$.
We must show that the elements $\uu{i}$ in
$B_{Q'}$ satisfy the defining relations of
$B_Q$. Note that the $t_i$ satisfy the defining relations for $B_{Q'}$.

Firstly, we check the relations (ii)
for an arrow incident with $k$.
Suppose that there is an arrow $i\rightarrow k$.
We have the following, using the fact that
$t_it_kt_i=t_kt_it_k$:
\begin{equation*}
\begin{split}
\uu{i} \uu{k} \uu{i} &= t_kt_it_kt_it_k^{-1} \\
&= t_k^2 t_i t_kt_k^{-1} = t_k^2 t_i.
\end{split}
\end{equation*}
Also,
$$
\uu{k} \uu{i} \uu{k} = t_k^2 t_k t_i t_k^{-1} t_k = t_k^2 t_i.$$
So
$$\uu{i} \uu{k} \uu{i}=\uu{k} \uu{i} \uu{k},$$
as required.

If there is an arrow $i\leftarrow k$, then we
have
$$\uu{i}\uu{k}\uu{i}=t_it_kt_i=t_kt_it_k=
\uu{k}\uu{i}\uu{k}.$$

Next, we consider relations (i) and (ii) for all
other arrows in $Q$. Relations of this kind
involving pairs of vertices which are not neighbours of $k$ follow immediately from the corresponding relations in $B_Q$.
If only one of the vertices in the relation is
a neighbour of $k$, the relation again follows
immediately since $t_k$ commutes with any
generator corresponding to a vertex not incident
with $k$ in $Q'$ (equivalently, in $Q$).
So we only need to consider the case
where both of the vertices in the pair
is incident with $k$ and we can use Lemma~\ref{l:localmutation}.

Going in either direction in part (a) of Lemma~\ref{l:localmutation}, the relation
$\uu{i}\uu{j}=\uu{j}\uu{i}$ follows from the
relation $t_it_j=t_jt_i$ in $B_{Q'}$, so we
consider part (b), firstly from left to right.
The cycle in $Q'$ gives the relation
$t_kt_i=t_jt_kt_it_jt_k^{-1}t_j^{-1}$.
Also applying the relation $t_k^{-1}t_j^{-1}t_k^{-1}=t_j^{-1}t_k^{-1}t_j^{-1}$, we obtain
\begin{equation*}
\begin{split}
\uu{i}\uu{j} &=t_kt_it_k^{-1}t_j \\
&= t_jt_kt_it_jt_k^{-1}t_j^{-1}t_k^{-1}t_j \\
&= t_jt_kt_it_k^{-1}=\uu{j}\uu{i}.
\end{split}
\end{equation*}
Going from right to left in part (b), we have,
using $t_jt_kt_j=t_kt_jt_k$, $t_it_j=t_jt_i$
and $t_it_kt_i=t_kt_it_k$,
\begin{equation*}
\begin{split}
\uu{j}\uu{i}\uu{j} &= t_kt_jt_k^{-1}t_it_kt_jt_k^{-1} \\
&= 
t_j^{-1}t_kt_jt_it_j^{-1}t_kt_j \\
&= t_j^{-1}t_kt_it_kt_j \\
&= t_j^{-1}t_it_kt_it_j \\
&= t_it_j^{-1}t_kt_jt_i \\
&= t_it_kt_jt_k^{-1}t_i \\
&= \uu{i} \uu{j} \uu{i}.
\end{split}
\end{equation*}

Next, we have to check that the $\uu{i}$ satisfy the relations of type (iii) for $Q$, so we need
to consider each the types of cycle described
in Lemma~\ref{l:cycles}. By Lemma \ref{lem:one-implies-all}, it is enough to check that,
for any given cycle in $Q$, one of the relations
in (iii) holds.

For part (a), we have
\begin{equation*}
\begin{split}
\uu{k}\uu{i}\uu{j}\uu{k} &=
t_kt_it_kt_jt_k^{-1}t_k \\
&= t_kt_it_kt_j,
\end{split}
\end{equation*}
while
\begin{equation*}
\begin{split}
\uu{i}\uu{j}\uu{k}\uu{i} &=
t_it_kt_jt_k^{-1}t_kt_i \\
&= t_it_kt_jt_i \\
&= t_it_kt_it_j,
\end{split}
\end{equation*}
which is equal to $\uu{k}\uu{i}\uu{j}\uu{k}$ as required.

For part (b) we have, applying a relation for the cycle in $Q'$ in the fourth step:
\begin{equation*}
\begin{split}
\uu{i_1}\uu{i_2}\cdots \uu{i_r}\uu{i_1}\uu{i_2}\cdots \uu{i_{r-2}}
&= t_kt_{i_1}t_k^{-1} t_{i_2}\cdots t_{i_r} t_kt_{i_1}t_k^{-1}t_{i_2}\cdots t_{i_{r-2}} \\
&= t_{i_1}^{-1}t_{k}t_{i_1} t_{i_2}\cdots t_{i_r} t_kt_{i_1}t_k^{-1}t_{i_2}\cdots t_{i_{r-2}} \\
&= t_{i_1}^{-1}t_{k}t_{i_1} t_{i_2}\cdots t_{i_r} t_kt_{i_1}t_{i_2}\cdots t_{i_{r-2}}t_k^{-1} \\
&= t_{i_1}^{-1}t_{i_1} t_{i_2}\cdots t_{i_r} t_kt_{i_1}t_{i_2}\cdots t_{i_{r-2}}t_{i_{r-1}}t_k^{-1} \\
&= t_{i_2}\cdots t_{i_r} t_kt_{i_1}t_k^{-1}t_{i_2}\cdots t_{i_{r-2}}t_{i_{r-1}} \\
&= \uu{i_2}\cdots \uu{i_r}\uu{i_1}\uu{i_2}\cdots \uu{i_{r-2}}\uu{i_{r-1}}.
\end{split}
\end{equation*}

For part (c), we have, applying a relation for
the cycle in $Q'$ in the fourth step:
\begin{equation*}
\begin{split}
\uu{i_1}\uu{i_2}\cdots \uu{i_{r-1}}\uu{i_r}\uu{k}\uu{i_1}\cdots \uu{i_{r-1}}
&=
t_{i_1}t_{i_2}\cdots t_{i_{r-1}}t_kt_{i_r}t_k^{-1}t_kt_{i_1}\cdots t_{i_{r-1}} \\
&=
t_{i_1}t_{i_2}\cdots t_{i_{r-1}}t_kt_{i_r}t_{i_1}\cdots t_{i_{r-1}} \\
&=
t_{i_1}t_kt_{i_2}\cdots t_{i_{r-1}}t_{i_r}t_{i_1}\cdots t_{i_{r-1}} \\
&=
t_{i_1}t_kt_{i_1}t_{i_2}\cdots t_{i_{r-1}}t_{i_r}t_{i_1}\cdots t_{i_{r-2}} \\
&=
t_kt_{i_1}t_kt_{i_2}\cdots t_{i_{r-1}}t_{i_r}t_{i_1}\cdots t_{i_{r-2}} \\
&=
t_kt_{i_1}t_{i_2}\cdots t_{i_{r-1}}t_kt_{i_r}t_k^{-1}t_kt_{i_1}\cdots t_{i_{r-2}} \\
&= \uu{k}\uu{i_1}\uu{i_2}\cdots \uu{i_{r-1}}\uu{i_r}\uu{k}\uu{i_1}\cdots \uu{i_{r-2}},
\end{split}
\end{equation*}
and we are done.
\end{pf}
\end{prop}

\begin{thm}\label{thm:groupiso}
 $\varphi_k:B_Q\to B_{Q'}$ is a group isomorphism.
 
\begin{pf}
As mutation is an involution, we can consider the composition
$$\varphi_k:B_Q\arr{\varphi_k} B_{Q'}\arr{\varphi_k} B_Q.$$
Fix some $i\in I$.  Note that mutation at $k$ does not change whether $i$ and $k$ are connected in the quiver: it just swaps the direction of any arrow that may exist between $i$ and $k$.  So if we have $i\to k$ then $s_i\mapsto t_kt_it_k\inv\mapsto s_ks_is_k\inv$.  If we have $i\leftarrow k$ then $s_i\mapsto t_i\mapsto s_ks_is_k\inv$.  And if there is no arrow between $i$ and $k$ then $s_i\mapsto t_i\mapsto s_i$.  But in this case $s_i$ and $s_k$ commute, so $s_i=s_ks_is_k\inv$.  Hence in every case $\varphi_k(s_i)=s_ks_is_k\inv$, so $\varphi_k$ is just a conjugation map and therefore $\varphi_k:B_Q\to B_{Q'}$ is an isomorphism.
\end{pf}
\end{thm}

\begin{rmk}\label{rmk:phi-inv}
The inverse of $\varphi_k$ is the group isomorphism $\varphi_k^{-1}:B_{Q'}\arr\sim B_Q$ defined by
$$
\varphi_k^{-1}(t_i)=
\begin{cases}
s_k^{-1} s_i s_k & \text{if } i\rightarrow k \text{ in } Q; \\
s_i & \text{otherwise.}
\end{cases}
$$
\end{rmk}

Noting Remark~\ref{r:artin}, we have the following:


\begin{thm} \label{t:groupinvariance}
If $Q$ is a mutation-Dynkin quiver of type $\Delta$ then $B_Q\cong B_\Delta$.
\end{thm}

\section{Topological interpretation of the generators}\label{sec:riemann}

\subsection{Braid groups}

In this section we consider quivers $Q$ which are mutation-equivalent to an orientation of the Dynkin diagram of type $\Delta$, where $\Delta=A_n$ or $D_n$. By Theorem~\ref{t:groupinvariance}, $B_Q$ is isomorphic to the Artin braid group $B_{\Delta}$ of
the same Dynkin type. In other words, $B_Q$
gives a presentation of $B_{\Delta}$. In this section we give a geometric interpretation of this
presentation.

We associate an oriented Riemann surface $S$ (with boundary) together with marked points $M$ to $\Delta$, as follows.
If $\Delta=A_n$, we take $S$ to be a disk with $n-3$
marked points on its boundary, as in~\cite{fz1,fz2}.
If $\Delta=D_n$, we take $S$ to be a disk with one marked point in its interior and $n$ marked points on its boundary, as in~\cite{fst08,schiffler08}.
In each case, it was shown that every quiver of the corresponding mutation type arises from a triangulation of $(S,M)$ (tagged, in the type $D_n$ case) in the following way.
We follow~\cite{fst08},
in a generality great enough to cover both cases (noting that there is at most one interior marked point).

A (simple) \emph{arc} in $(S,M)$ is a curve in $S$ (considered up to isotopy) whose endpoints are marked points in $M$ and which does not have any self-crossing points, except possibly at its endpoints.
Apart from these endpoints, it must be disjoint from $M$ and the boundary of $S$, and it must not cut 
out an unpunctured one- or two-sided polygon.

Two arcs are said to be \emph{compatible}
if they are non-crossing in the interior
of $S$.
A maximal set of compatible arcs is a \emph{triangulation}.

A \emph{tagged arc}
in $(S,M)$ is an arc which does not
cut out a once-punctured monogon; each of its ends is tagged, either plain or notched. Plain tags are omitted, while notched tags are displayed using the
bow-tie symbol $\bowtie$.
An end incident with a boundary marked point is always tagged plain.
Two tagged arcs $\alpha,\beta$ are \emph{compatible} if
\begin{enumerate}
\item[(i)] the untagged arcs underlying $\alpha$ and $\beta$ are compatible, and
\item[(ii)] if the untagged versions of $\alpha$ and $\beta$ are different but share an endpoint, then the corresponding ends
of $\alpha$ and $\beta$ are tagged in the same way.
\end{enumerate}
A \emph{tagged triangulation}
$\T$ of $(S,M)$ is a maximal collection
of tagged arcs in $(S,M)$. Note that
if none of the marked points in $M$
lies in the interior of $S$, every end of
an arc in a tagged triangulation must
be tagged plain, and tagged triangulations
of $S$ can be identified with triangulations
of $S$.

The set $M$ of marked points divides the boundary 
components of $(S,M)$ into connected components,
which we call \emph{boundary arcs}.
Note that the boundary arcs do not
lie in a triangulation or tagged
triangulation of $(S,M)$, by definition.

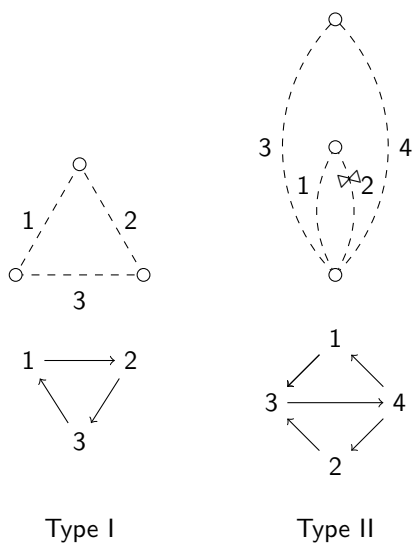
\begin{figure}
$$
\begin{tikzpicture}
[scale=1.7,baseline=(bb.base), quivarrow/.style={black, -latex}, arc/.style={black,dashed}, braid/.style={black}] 

\path (0,0) node (bb) {}; 

\begin{pgfonlayer}{background}
\draw[arc] (0,0) -- (1,0) -- (0.5,0.866) -- (0,0);
\end{pgfonlayer}

\filldraw [draw=black, fill=white] (0,0) circle (0.05);
\filldraw [draw=black, fill=white] (0.5,0.866) circle (0.05);
\filldraw [draw=black, fill=white] (1,0) circle (0.05);

\node (E1) at (0.1,0.43) {$1$};
\node (E2) at (0.9,0.43) {$2$};
\node (E3) at (0.5,-0.2) {$3$};

\begin{scope}[shift={(0,-1.1)}]

\node (Q1) at (0.1,0.43) {$1$};
\node (Q2) at (0.9,0.43) {$2$};
\node (Q3) at (0.5,-0.2) {$3$};

\draw[->] (Q1) -- (Q2);
\draw[->] (Q2) -- (Q3);
\draw[->] (Q3) -- (Q1);

\node at (0.5,-0.9) {Type I};

\end{scope}

\begin{scope}[shift={(2.5,0)}]

\filldraw [draw=black, fill=white] (0,0) circle (0.05);
\filldraw [draw=black, fill=white] (0,1) circle (0.05);
\filldraw [draw=black, fill=white] (0,2) circle (0.05);

\node (F1) at (-0.25,0.7) {$1$};
\node (F2) at (0.25,0.7) {$2$};
\node (F3) at (-0.55,1) {$3$};
\node (F4) at (0.55,1) {$4$};

\begin{pgfonlayer}{background}
\draw[arc] (0,0) to[out=45,in=-45] (0,2);
\draw[arc] (0,0) to[out=135,in=-135] (0,2);
\draw[arc] (0,0) to[out=60,in=-60] node[near end,sloped] {$\vertbowtie$} (0,1);
\draw[arc] (0,0) to[out=120,in=-120]  (0,1);
\end{pgfonlayer}

\begin{scope}[shift={(-0.5,-1)}]
\node (R3) at (0,0) {$3$};
\node (R4) at (1,0) {$4$};
\node (R1) at (0.5,0.5) {$1$};
\node (R2) at (0.5,-0.5) {$2$};

\draw[->] (R1) -- (R3);
\draw[->] (R3) -- (R4);
\draw[->] (R4) -- (R2);
\draw[->] (R4) -- (R1);
\draw[->] (R1) -- (R3);
\draw[->] (R2) -- (R3);
\node at (0.5,-1) {Type II};

\end{scope}
\end{scope}
\end{tikzpicture}
$$
\caption{Puzzle pieces for tagged triangulations in types $A_n$ and $D_n$ and the corresponding quivers}
\label{f:puzzlepieces}
\end{figure}

The tagged triangulation $\T$ can be built up by
gluing together puzzle pieces of the two types
shown in Figure~\ref{f:puzzlepieces} (see~\cite[Rk.~4.2]{fst08}) by gluing together along
boundary arcs. Note that the puzzle piece of type
II can only occur in the type $D_n$ case, and then
it occurs exactly once.

If $\alpha$ is an arc in a tagged triangulation
$\T$, then the \emph{flip} of $\T$ at $\alpha$
is the unique tagged triangulation containing
$\T\setminus \{\alpha\}$ but not containing $\alpha$.
By~\cite{fst08}, the set of tagged triangulations of $(S,M)$ is connected under flips.

The \emph{quiver} $Q_{\T}$ of a tagged triangulation $\T$ has vertices corresponding to the arcs in $\T$.
The quiver is built up by associating a
quiver to each puzzle piece; see Figure~\ref{f:puzzlepieces}. If a boundary arc
in the puzzle piece is also a boundary arc of
$(S,M)$, then the corresponding vertex in the
quiver is omitted, together with all incident arrows. The quivers are then glued together by identifying vertices whenever the corresponding edges are glued together in the puzzle pieces.

In order to discuss braid groups, we need to
consider more general curves in $(S,M)$.
We define a \emph{path} in $(S,M)$ to be a
(possibly non-simple) curve whose endpoints lie
in $S$ (not necessarily in $M$).

\begin{defn}
Let $\T$ be a tagged triangulation of
$(S,M)$. We associate a graph to $\T$,
which we call the \emph{braid graph}
$G_{\T}$ of $\T$, as follows. The vertices $V_{\T}$
of $G_{\T}$ are in bijection with the connected components of
the complement of $\T$ in $(S,M)$
and, whenever two such connected components
have a common tagged arc on their boundaries, there is an edge
in $G_{\T}$ between the corresponding
vertices. Thus the edges in $G_{\T}$
are in bijection with the arcs in $\T$.

We choose an embedding $\iota$ of $G_{\T}$
into $(S,M)$, mapping each vertex to an
interior point of the corresponding
connected component of the complement
of $\T$ in $(S,M)$ and each edge to a path
between the images of its endpoints transverse to the corresponding 
edge in $\T$. We identify $G_{\T}$ with
its image under $\iota$.
\end{defn}

Note that in the type $A$ case the braid graph is the tree from Section 3.1 of \cite{ccs}.

We associate an orbifold $X$
to $S$ as follows. In the type $A_n$
case, we just take $X=S$, and in the type $D_n$ case we take
$X$ to be $S$ with the interior
marked point of $S$ interpreted as a cone
point of order two. In each case, the
set $M$ of marked points induces a
corresponding set of
marked points in $X$, which we also denote
by $M$. Each arc or tagged arc $\alpha$ in $(S,M)$ induces a corresponding arc or
tagged arc in $(X,M)$ which we
also denote by $\alpha$. Thus each
(tagged) triangulation $\T$ of $(S,M)$
induces a corresponding set $\T$ of
(tagged) arcs in $(X,M)$.

Note also that orbifolds have been used to
model cluster algebras in~\cite{fst12}.
In this approach, the model for $B_n$ is an orbifold with a cone point of order $2$, regarded as a folding of $D_n$,
where $D_n$ is modelled by a disk with a single interior marked point (see also Lecture 15 of~\cite{thurston}, which was given by A. Felikson).

We denote by $X^{\circ}$ the orbifold $X$
with the cone point (if there is one) removed (so $X^{\circ}=X$ in type $A_n$).
Given any set $V$ of vertices in $X^{\circ}$, we may consider the corresponding \emph{braid group}, $\Gamma(X,V)$ as defined in~\cite{allcock02}. Each element of $\Gamma(X,V)$ (or \emph{braid}) can be
regarded as a permutation $g$ of $V$ together with a tuple $\gamma=(\gamma_v)_{v\in V}$ of paths $\gamma_v:[0,1]\rightarrow X^{\circ}$ such
that $\gamma_v(0)=v$ and $\gamma_v(1)=g(v)$ for each $v\in V$. In addition, for each $t\in [0,1]$, the points
$\gamma_v(t)$ for $v\in V$ must all be distinct
for all $v\in V$. Braids are considered
up to isotopy, and two braids can
be multiplied by composing the paths in
a natural way; we compose braids from right
to left, as for functions.

\begin{rmk} \label{r:isotopy}
Suppose $V$ and $V'$ are two sets of points in $X^{\circ}$ and there
is a bijection $\rho:V\rightarrow V'$.
Suppose also that there is a set of
paths $\delta_v:[0,1]\rightarrow X^{\circ}$, for $v\in V$, with $\delta_v(0)=v$
and $\delta_v(1)=\rho(v)$ for all $v\in V$.
Suppose furthermore that that the points
$\gamma_v(t)$ for $v\in V$ and $t\in [0,1]$ are all distinct. Then the maps $\delta_v$ induce a natural
isomorphism between $\Gamma(X,V)$ and $\Gamma(X,V')$.
\end{rmk}

\begin{defn} \label{d:pathbraid}
Each path $\pi$ in $X^{\circ}$ with endpoints $v_1,v_2$
in $V$ determines a braid $\sigma_{\pi}$ in $\Gamma(X,V)$ as follows (see~\cite[\S7]{foxneuwirth62}).
We thicken the path $\pi$ along
its length (avoiding the other vertices),
closing it off at the end points to form a (topological) disk. We give the boundary of the disk the clockwise orientation. The vertices $v_1$ and $v_2$ divide the boundary of the disk into two paths, one from $v_1$
to $v_2$ and the other from $v_2$ to $v_1$.
We define $\gamma_{v_1}$ to be the former
and $\gamma_{v_2}$ to be the latter.
See Figure~\ref{f:thickening}.
For $v\in V$ such that $v\not=v_1,v_2$, we define $\gamma_v(t)$ to be $v$ for all $t\in [0,1]$. Then $\sigma_{\pi}$ is the braid
$(\gamma_v)_{v\in V}$.
Note that $\sigma_{\pi}$ only depends on the isotopy class of the image of $\pi$ in $(X,V)$.
In particular, it is unchanged if $\pi$
is reversed.
\end{defn}
 
\begin{figure}
$$
\begin{tikzpicture}
[scale=1.7,baseline=(bb.base), quivarrow/.style={black, -latex}, arc/.style={black,dashed}, braid/.style={black}] 

\path (0,0) node (bb) {}; 

\draw (0,0) -- (0,1) -- (1,1.4) -- (2.4,1.4) -- (2.4,0.4) -- (1.6,0) -- (0,0);
\draw (0,1) -- (1.6,1) -- (2.4,1.4);
\draw (1.6,0) -- (1.6,1);
\draw[dotted] (0,0) -- (1,0.4) -- (1,1.4);
\draw[dotted] (1,0.4) -- (2.4,0.4);

\filldraw [draw=black, fill=black] (0.8,1.15) circle (0.03);

\filldraw [draw=black, fill=black] (1.4,1.15) circle (0.03);

\filldraw [draw=black, fill=black] (0.8,0.15) circle (0.03);

\filldraw [draw=black, fill=black] (1.4,0.15) circle (0.03);

\draw [braid] plot[smooth]  coordinates {(0.8,0.15) (0.9,0.55) (1.3,0.75) (1.4,1.15)};

\draw [braid] plot[smooth]  coordinates {(0.8,0.15) (0.9,0.55) (1.3,0.75) (1.4,1.15)};

\draw [braid] plot[smooth]  coordinates {(0.8,1.15) (0.9,0.75) (1.05,0.67)};
\draw[braid] plot[smooth] coordinates {(1.15,0.63) (1.3,0.55) (1.4,0.15)};

\draw[arc] plot[smooth] coordinates {(0.8,1.15) (1.4,1.15)};
\draw[arc] plot[smooth] coordinates {(0.8,0.15) (1.4,0.15)};

\node at (1.2,1.25) {$\pi$};
\node at (1.2,0.25) {$\pi$};

\end{tikzpicture}
$$
\caption{The braid $\sigma_{\pi}$}
\label{f:sigmaalpha}
\end{figure}
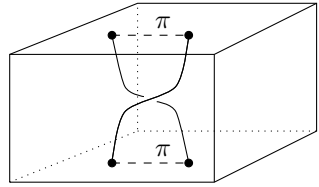

An example of a braid $\sigma_{\pi}$ is displayed as a picture
(in the same way as in~\cite{allcock02}) in Figure~\ref{f:sigmaalpha}. In this figure only, we display
$\pi$ as a dashed line to distinguish it from the braid $\sigma_{\pi}$. 

\begin{figure}
$$
\begin{tikzpicture}
[scale=3,baseline=(bb.base), quivarrow/.style={black, -latex}, arc/.style={black,dashed}, braid/.style={black}] 
\newcommand{\strarrow}{\arrow{angle 60}}

\path (0,0) node (bb) {}; 

\filldraw [draw=black, fill=black] (0,0) circle (0.015);

\filldraw [draw=black, fill=black] (1,0) circle (0.015);

\draw[braid] plot[smooth] coordinates {(0,0) (0.2,0.2) (0.5,0) (0.7,-0.3) (1,0)};

\node at (-0.15,0) {$v_1$};
\node at (1.15,0) {$v_2$};

\draw[braid] plot[smooth] coordinates {(0,0) (0.2,0.3) (0.5,0.1) (0.7,-0.2) (1,0)}[ postaction=decorate, decoration={markings, 
  mark= at position 0.5 with \strarrow}];

\node at (0.5,0.25) {$\gamma_{v_1}$};

\draw[braid] plot[smooth] coordinates {(1,0) (0.7,-0.4) (0.5,-0.1) (0.2,0.1) (0,0)}[ postaction=decorate, decoration={markings, 
  mark= at position 0.5 with \strarrow}];

\node at (0.45,-0.25) {$\gamma_{v_2}$};

\end{tikzpicture}
$$
\caption{Thickening of the path $\pi$ ($\pi$ is
the middle path)}
\label{f:thickening}
\end{figure}
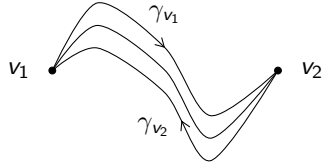

\subsection{Interpretation of the generators}

Let $\T$ be a triangulation of $(S,M)$.
Let $Q_{\T}$ be the quiver of $\T$.
Then $Q_{\T}$ has vertices $I$ corresponding to
the arcs in $\T$. We denote the arc in $\T$ associated to $i\in I$ by $\alpha_i$.
The corresponding edge in $G_{\T}$ is denoted
$\pi_i$.
Let $\sigma_i=\sigma_{\pi_i}\in \Gamma(X,V_{\T})$ be the corresponding braid.
We define $H_{\T}$ to be the subgroup of  
$\Gamma(X,V_{\T_0})$
generated by the braids $\sigma_i$ for
$i\in I$.

\begin{figure}
$$
\begin{tikzpicture}
[scale=1.7,baseline=(bb.base), quivarrow/.style={black, -latex}, arc/.style={black,dashed}, braid/.style={black}] 

\path (0,0) node (bb) {}; 

\begin{pgfonlayer}{background}
\draw[arc] (0,0) circle (1.6);

\draw[arc] (225:1.6) -- (110:1.6);
\draw[arc] (225:1.6) -- (75:1.6);
\draw[arc] (225:1.6) -- (38:1.6);
\draw[arc] (225:1.6) -- (0:1.6);
\draw[arc] (225:1.6) -- (-35:1.6);
\end{pgfonlayer}

\node at (225:1.8) {$P$};

\filldraw [draw=black, fill=white] (225:1.6) circle (0.05);
\filldraw [draw=black, fill=white] (110:1.6) circle (0.05);
\filldraw [draw=black, fill=white] (75:1.6) circle (0.05);
\filldraw [draw=black, fill=white] (38:1.6) circle (0.05);
\filldraw [draw=black, fill=white] (0:1.6) circle (0.05);
\filldraw [draw=black, fill=white] (-35:1.6) circle (0.05);

\node (A1) at (-1.15,0.6) {};
\filldraw [draw=black, fill=black] (A1) circle (0.05);
\node (B1) at (-0.32,1) {};
\filldraw [draw=black, fill=black] (B1) circle (0.05);
\node (C1) at (0.38,0.8) {};
\filldraw [draw=black, fill=black] (C1) circle (0.05);
\node (D1) at (1.02,0.25) {};
\filldraw [draw=black, fill=black] (D1) circle (0.05);
\node (E1) at (0.8,-0.65) {};
\filldraw [draw=black, fill=black] (E1) circle (0.05);
\node (F1) at (0,-1.3) {};
\filldraw [draw=black, fill=black] (F1) circle (0.05);

\draw [braid] plot  coordinates {(A1) (B1)};
\draw [braid] plot  coordinates {(B1) (C1)};
\draw [braid] plot  coordinates {(C1) (D1)};
\draw [braid] plot  coordinates {(D1) (E1)};
\draw [braid] plot  coordinates {(E1) (F1)};

\node at (-1.1,0) {$\alpha_1$};
\node at (-0.6,0.15) {$\alpha_2$};
\node at (0.2,0.23) {$\alpha_3$};
\node at (0.4,-0.37) {$\alpha_4$};
\node at (0,-0.92) {$\alpha_5$};

\begin{scope}[shift={(-1,-2.7)}]

\node (V1) at (0,0) {$1$};
\node (V2) at (0.5,0) {$2$};
\node (V3) at (1,0) {$3$};
\node (V4) at (1.5,0) {$4$};
\node (V5) at (2,0) {$5$};

\draw[->] (V2) -- (V1);
\draw[->] (V3) -- (V2);
\draw[->] (V4) -- (V3);
\draw[->] (V5) -- (V4);

\end{scope}

\begin{scope}[shift={(4,0)}]

\begin{pgfonlayer}{background}
\draw[arc] (0,0) circle (1.6);

\draw [arc] plot[smooth]  coordinates {(0,-1.6) (-0.25,-1) (-0.25,-0.5) (0,0)};
\draw [arc] plot[smooth]  coordinates {(0,-1.6) (0.25,-1) (0.25,-0.5) (0,0)}[postaction=decorate,
decoration={markings, mark= at position 0.85 with
\node[transform shape] {$\vertbowtie$};}];

\draw[arc] (-1.6,0) to[out=20,in=180] (0,0.35) to[out=0,in=50] (0,-1.6);

\draw[arc] (-1.6,0) to[out=30,in=180] (0.3,0.5) to[out=0,in=90] (-60:1.6);

\draw[arc] (-1.6,0) to[out=40,in=180] (0.4,0.7) to[out=0,in=110] (-20:1.6);

\draw[arc] (-1.6,0) to[out=50,in=180] (0.5,0.9) to[out=0,in=110] (10:1.6);

\draw[arc] (-1.6,0) to[out=60,in=-150] (70:1.6);


\end{pgfonlayer}

\node at (180:1.8) {$P$};
\node at (270:1.8) {$Q$};

\filldraw [draw=black, fill=white] (0,-1.6) circle (0.05);
\filldraw [draw=black, fill=white] (0,0) circle (0.05);
\filldraw [draw=black, fill=white] (-60:1.6) circle (0.05);
\filldraw [draw=black, fill=white] (-20:1.6) circle (0.05);
\filldraw [draw=black, fill=white] (10:1.6) circle (0.05);
\filldraw [draw=black, fill=white] (70:1.6) circle (0.05);
\filldraw [draw=black, fill=white] (-1.6,0) circle (0.05);

\node (A2) at (-0.7,1.2) {};
\filldraw [draw=black, fill=black] (A2) circle (0.05);

\node (B2) at (0.6,1.15) {};
\filldraw [draw=black, fill=black] (B2) circle (0.05);

\node (C2) at (1.2,0.5) {};
\filldraw [draw=black, fill=black] (C2) circle (0.05);

\node (D2) at (0.9,0.2) {};
\filldraw [draw=black, fill=black] (D2) circle (0.05);

\node (E2) at (0.55,0.15) {};
\filldraw [draw=black, fill=black] (E2) circle (0.05);

\node (F2) at (-0.8,0) {};
\filldraw [draw=black, fill=black] (F2) circle (0.05);

\node (G2) at (0,-0.8) {};
\filldraw [draw=black, fill=black] (G2) circle (0.05);

\draw [braid] plot  coordinates {(A2) (B2)};
\draw [braid] plot  coordinates {(B2) (C2)};
\draw [braid] plot  coordinates {(C2) (D2)};
\draw [braid] plot  coordinates {(D2) (E2)};
\draw [braid] plot  coordinates {(F2) (G2)};

\draw[braid] (E2) to[out=170,in=20] (F2);
\draw[braid] (F2) to[out=10,in=170] (0.1,0.1) to[out =-10,in=70] (G2);

\node at (0.1,1.4) {$\alpha_1$};
\node at (0.2,1) {$\alpha_2$};
\node at (1.17,-0.1) {$\alpha_3$};
\node at (1,-0.7) {$\alpha_4$};
\node at (0.55,-0.9) {$\alpha_5$};
\node at (-0.35,-1.2) {$\alpha_6$};
\node at (0.05,-1.2) {$\alpha_7$};

\begin{scope}[shift={(-1.2,-2.7)}]

\node (V1) at (0,0) {$1$};
\node (V2) at (0.5,0) {$2$};
\node (V3) at (1,0) {$3$};
\node (V4) at (1.5,0) {$4$};
\node (V5) at (2.1,0) {$5$};
\node (V6) at (2.5,0.4) {$6$};
\node (V7) at (2.5,-0.4) {$7$};

\draw[->] (V2) -- (V1);
\draw[->] (V3) -- (V2);
\draw[->] (V4) -- (V3);
\draw[->] (V5) -- (V4);
\draw[->] (V6) -- (V5);
\draw[->] (V7) -- (V5);

\end{scope}

\end{scope}

\end{tikzpicture}
$$
\caption{Initial triangulations and the corresponding
braid graphs and quivers}
\label{f:initialtriangulations}
\end{figure}
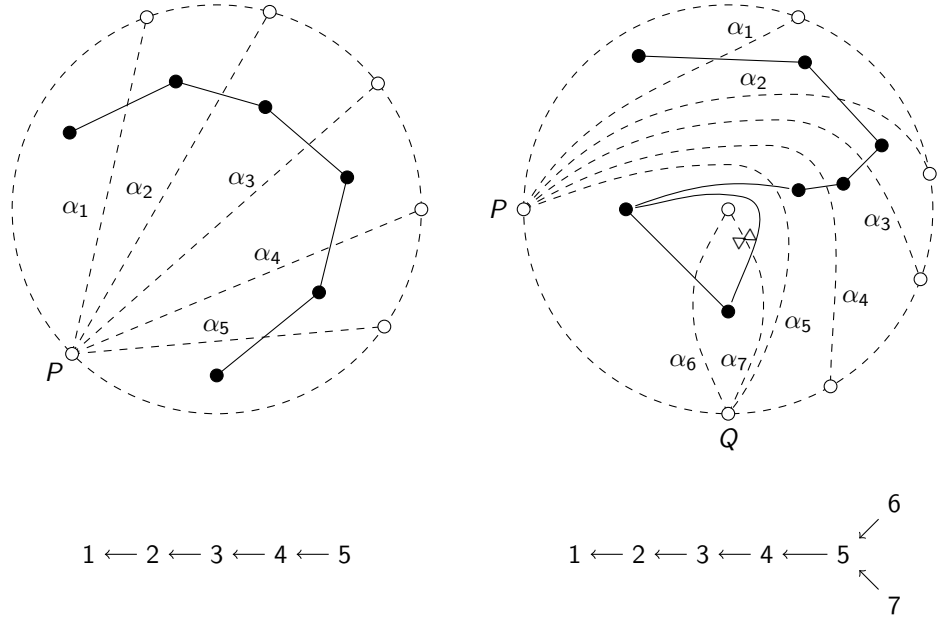

Let $\T_0$ be an initial triangulation of
$(S,M)$ defined as follows. In the type $A_n$
case, we choose a marked point $P$ in $M$ and
take noncrossing arcs between $P$ and each of the other marked points in $M$ not incident with a
boundary arc incident with $P$. In the type $D_n$
case, we choose two marked points $P,Q$ on the boundary of $S$. We take two arcs between the
interior marked point and $Q$, one tagged plain
at the interior marked point and the other one
tagged notched, and an arc between $P$ and
$Q$ (not homotopic to a boundary arc).
We then take (noncrossing)
arcs between $P$ and every other marked point
in $M$ on the boundary of $S$ not incident
with a boundary arc incident with $P$.
See Figure~\ref{f:initialtriangulations}.
Then the quiver $Q_{\T_0}$
associated to $Q_{\T_0}$ is a Dynkin quiver
of type $\Delta$. By Remark~\ref{r:artin},
$B_{Q_{\T_0}}$ is isomorphic to the Artin braid
group of type $\Delta$.

\begin{prop}\label{p:initialcase}
Let $\T_0$ be the triangulation of $(X,M)$
defined as above.
Then there is an isomorphism from $H_{\T_0}$
to $B_{Q_{\T_0}}$ taking the braid $\sigma_i$ to
the generator $s_i$ of $B_{Q_{\T_0}}$.
Furthermore, in type $A_n$, $H_{\T_0}$ coincides
with $\Gamma(X,V_{\T_0})$, while in type $D_n$,
$H_{\T_0}$ is a subgroup of $\Gamma(X,V_{\T_0})$ of index two.

\begin{pf}
For type $A_n$, see~\cite{foxneuwirth62}
and the explanation
in~\cite[\S4]{allcock02}.
For type $D_n$, note that the elements $\sigma_i$ for
$i\in I$ coincide with the generators $h_i$
defined in~\cite[\S1]{allcock02} (via an
isomorphism of the kind in Remark~\ref{r:isotopy}). The result then follows from~\cite[Thm.~1]{allcock02}.
\end{pf}
\end{prop}

The following lemma appears in~\cite[Th\'{e}or\`{e}me, part (iv)]{sergiescu93}.

\begin{lem}\label{l:threecycle}
Let $A,B,C$ be three distinct points in $X^{\circ}$ and suppose there is a topological disk in $X^{\circ}$, with $A$, $B$ and $C$ lying in order clockwise around its boundary. Let $AB$ denote the arc on this boundary between $A$ and $B$. We define $BC$ and $CA$ similarly.  Then $\sigma_{AB}\sigma_{BC}=\sigma_{BC}\sigma_{CA}$.
\end{lem}

\begin{thm} \label{t:braidisomorphism}
Let $\T$ be an arbitrary tagged triangulation of $(X,M)$. Then there is an isomorphism from $H_{\T}$
to $B_{Q_{\T}}$ taking the braid $\sigma_i$ to
the generator $s_i$ of $B_{Q_{\T}}$.
Furthermore, in type $A_n$, $H_{\T}$ coincides
with $\Gamma(X,V)$, while in type $D_n$,
$H_{\T}$ is a subgroup of $\Gamma(X,V)$ of index two.

\begin{pf}
The result holds for $\T=\T_0$ by Proposition~\ref{p:initialcase}. Note that any triangulation
can be obtained from $\T_0$ by applying a finite
number of flips of tagged triangulations.
We show that the theorem is true for
an arbitrary tagged triangulation $\T$ by induction on the number of flips required to obtain $\T$
from $\T_0$. To do this, it is enough to show that
if the theorem holds for a tagged triangulation
$\T$ and $\alpha_i$ is a tagged arc in $\T$ then
the theorem also holds for the flip of $\T$
at $\alpha_i$.

So we assume the result holds for a tagged
triangulation $\T$. Thus there is an isomorphism
$\psi_{\T}:H_{\T}\rightarrow B_{Q_{\T}}$ sending
$\sigma_i$ to $s_i$. We denote the corresponding
elements of $H_{\T'}$ by $\tau_i$ and $t_i$.
The tagged arcs in $\T$ are denoted by $\alpha_i$,
for $i\in I$, and we denote the corresponding tagged arcs in $\T'$ by $\beta_i$, for $i\in I$.
The edges of $G_{\T}$ are denote $\pi_i$, and
we denote the edges of $G_{\T'}$ by $\rho_i$.

We define:
$$\widetilde{\tau_i}=\begin{cases}
\sigma_k^{-1}\sigma_i\sigma_k, & \text{if $i\rightarrow k$ in $Q$;} \\
\sigma_i, & \text{otherwise}.
\end{cases}
$$
Then it is easy to see that
$H_{\T}$ is generated by the $\widetilde{\tau_i}$
for $i\in I$.

We consider the possible types of flip that can
occur, which are determined by the fact that $\T$
can be constructed out of puzzle pieces.
Suppose first that
$\T'$ is the flip of $\T$ at an arc $\alpha$
where two puzzle pieces of type (I) are glued
together. We label the corresponding vertices
in $I$ by $1,2,3,4,5$, for simplicity, and
suppose we are flipping at the edge in $\T$
dual to $\alpha_1$.
The braid graph local to the flip is shown
in the left hand diagram in Figure~\ref{f:flip1}. 
Applying Lemma~\ref{l:threecycle}, we see that the middle figure shows paths $\widetilde{\pi_i}$
with the property that $\widetilde{\tau_i}=\sigma_{\widetilde{\pi_i}}$ for $i=1,2,3,4,5$.

\begin{figure}
$$
\begin{tikzpicture}[scale=2,baseline=(bb.base), quivarrow/.style={black, -latex}, arc/.style={black,dashed}, braid/.style={black}] 

\path (0,0) node (bb) {}; 

\newcommand{\qedge}{0.5cm} 

\draw[->] (1.65,0.2) -- (1.85,0.2);
\draw[->] (4.15,0.2) -- (4.35,0.2);

\begin{pgfonlayer}{background}
\draw[arc] (0,0) -- (1,0) -- (1,1) -- (0,1) -- (0,0);

\draw[arc] (0,1) -- (1,0);
\end{pgfonlayer}

\node (C1) at (-0.3,0.5) {$\bullet$};
\node [left] at (C1) {\small {$C$}}; 

\node (D1) at (0.5,-0.3) {$\bullet$};
\node [below] at (D1) {\small {$D$}}; 

\node (F1) at (0.5,1.3) {$\bullet$};
\node [above] at (F1) {\small {$F$}}; 

\node (E1) at (1.3,0.5) {$\bullet$};
\node [right] at (E1) {\small {$E$}}; 

\node (A1) at (0.3,0.3) {$\bullet$};
\node at (0.37,0.23) {\small {$A$}}; 

\node (B1) at (0.7,0.7) {$\bullet$};
\node at (0.77,0.56) {\small {$B$}}; 

\draw [braid] plot  coordinates {(C1) (A1)};
\draw [braid] plot  coordinates {(A1) (D1)};
\draw [braid] plot  coordinates {(F1) (B1)};
\draw [braid] plot  coordinates {(B1) (E1)};
\draw [braid] plot  coordinates {(A1) (B1)};

\node at (0.35,0.47) {\small $\pi_1$};
\node at (-0.13,0.55) {\small $\pi_3$};
\node at (0.6,-0.15) {\small $\pi_2$};
\node at (1.13,0.65) {\small $\pi_5$};
\node at (0.66,1.13) {\small $\pi_4$};

\begin{scope}[shift={(0.5,-1.4)}]

\node (V1) at (0,0) {$1$};
\node (V5) at (\qedge,0) {$5$};
\node (V3) at (-\qedge,0) {$3$};
\node (V4) at (0,\qedge) {$4$};
\node (V2) at (0,-\qedge) {$2$};

\draw[->] (V1) -- (V2);
\draw[->] (V1) -- (V4);
\draw[->] (V5) -- (V1);
\draw[->] (V3) -- (V1);
\draw[->] (V4) -- (V5);
\draw[->] (V2) -- (V3);
\end{scope}

\begin{scope}[shift={(2.5,0)}]
\begin{pgfonlayer}{background}
\draw[arc] (0,0) -- (1,0) -- (1,1) -- (0,1) -- (0,0);

\draw[arc] (0,1) -- (1,0);
\end{pgfonlayer}

\node (C2) at (-0.3,0.5) {$\bullet$};
\node [left] at (C2) {\small {$C$}}; 

\node (D2) at (0.5,-0.3) {$\bullet$};
\node [below] at (D2) {\small {$D$}}; 

\node (F2) at (0.5,1.3) {$\bullet$};
\node [above] at (F2) {\small {$F$}}; 

\node (E2) at (1.3,0.5) {$\bullet$};
\node [right] at (E2) {\small {$E$}}; 

\node (A2) at (0.3,0.3) {$\bullet$};
\node at (0.37,0.23) {\small {$A$}}; 

\node (B2) at (0.7,0.7) {$\bullet$};
\node at (0.77,0.56) {\small {$B$}}; 

\draw [braid] plot  coordinates {(C2) (A2)};
\draw [braid] plot  coordinates {(B2) (D2)};
\draw [braid] plot  coordinates {(F2) (A2)};
\draw [braid] plot  coordinates {(B2) (E2)};
\draw [braid] plot  coordinates {(A2) (B2)};

\node at (0.53,0.64) {\small $\widetilde{\pi_1}$};
\node at (-0.13,0.54) {\small $\widetilde{\pi_3}$};
\node at (0.65,-0.15) {\small $\widetilde{\pi_2}$};
\node at (1.13,0.63) {\small $\widetilde{\pi_5}$};
\node at (0.57,1.13) {\small $\widetilde{\pi_4}$};
\end{scope}

\begin{scope}[shift={(5,0)}]
\begin{pgfonlayer}{background}
\draw[arc] (0,0) -- (1,0) -- (1,1) -- (0,1) -- (0,0);

\draw[arc] (0,0) -- (1,1);
\end{pgfonlayer}

\node (C3) at (-0.3,0.5) {$\bullet$};
\node [left] at (C3) {\small {$C$}}; 

\node (D3) at (0.5,-0.3) {$\bullet$};
\node [below] at (D3) {\small {$D$}}; 

\node (F3) at (0.5,1.3) {$\bullet$};
\node [above] at (F3) {\small {$F$}}; 

\node (E3) at (1.3,0.5) {$\bullet$};
\node [right] at (E3) {\small {$E$}}; 

\node (A3) at (0.3,0.7) {$\bullet$};
\node at (0.39,0.73) {\small {$A$}}; 

\node (B3) at (0.7,0.3) {$\bullet$};
\node at (0.77,0.22) {\small {$B$}}; 

\draw [braid] plot  coordinates {(C3) (A3)};
\draw [braid] plot  coordinates {(B3) (D3)};
\draw [braid] plot  coordinates {(F3) (A3)};
\draw [braid] plot  coordinates {(B3) (E3)};
\draw [braid] plot  coordinates {(A3) (B3)};

\node at (0.66,0.5) {\small $\rho_1$};
\node at (-0.13,0.64) {\small $\rho_3$};
\node at (0.67,-0.15) {\small $\rho_2$};
\node at (1.13,0.54) {\small $\rho_5$};
\node at (0.55,1.13) {\small $\rho_4$};

\begin{scope}[shift={(0.5,-1.4)}]

\node (V1) at (0,0) {$1$};
\node (V5) at (\qedge,0) {$5$};
\node (V3) at (-\qedge,0) {$3$};
\node (V4) at (0,\qedge) {$4$};
\node (V2) at (0,-\qedge) {$2$};

\draw[->] (V2) -- (V1);
\draw[->] (V4) -- (V1);
\draw[->] (V1) -- (V5);
\draw[->] (V1) -- (V3);
\draw[->] (V3) -- (V4);
\draw[->] (V5) -- (V2);
\end{scope}

\end{scope}

\end{tikzpicture}
$$
\caption{Flip involving an arc ($\alpha_1$) where two puzzle pieces of type I are glued}
\label{f:flip1}
\end{figure}
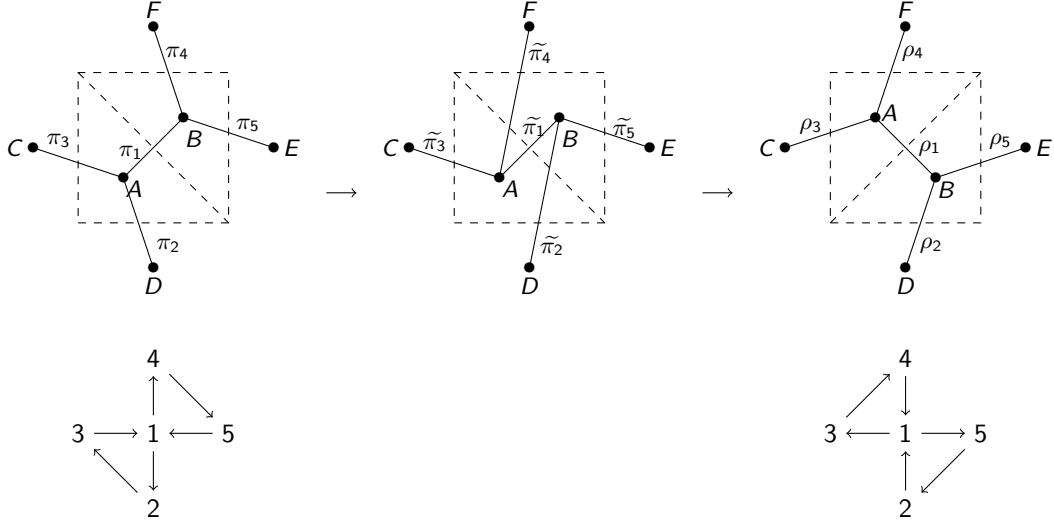


Rotating vertices $A$ and $B$ clockwise, to
get the right hand diagram in Figure~\ref{f:flip1},
we obtain, via Remark~\ref{r:isotopy},
an isomorphism from $H_{\T}$ to
$H_{\T'}$ taking $\widetilde{\tau_i}$ to
$\tau_i$ for all $i\in I$. The inverse is
an isomorphism from $H_{\T'}$ to $H_{\T}$
taking $\tau_i$ to $\sigma^{-1}_k\sigma_i\sigma_k$
if there is an arrow $i\rightarrow k$ in $Q$
and to $\sigma_i$ otherwise. Composing with
the isomorphism
$\varphi_k\circ \psi_{\T}$, where $\varphi_k$
is the isomorphism in Proposition~\ref{p:groupinvariance}, we obtain an isomorphism from $H_{\T'}$ to $B_{Q_{\T'}}$ taking $\tau_i$ to $t_i$ as required. This proves the required result in type A, so we are left with the
type D case, where puzzle pieces of type II may occur.

We next consider a flip inside a puzzle piece
of type II. We can apply essentially the same argument: see Figures~\ref{f:flip2} and~\ref{f:flip3}. Here we draw
the puzzle piece together with the two adjacent
triangles, necessarily of type I (since there
is only one cone point). We use the fact that
in the right hand diagram of Figure~\ref{f:flip3}, the resulting path $\widetilde{\pi_1}$ is isotopic to the path $\rho_1$ in $G_{\T'}$, using the fact that the cone point has order two.


\begin{figure}
$$
\begin{tikzpicture}[scale=1.7,baseline=(bb.base), quivarrow/.style={black, -latex}, arc/.style={black,dashed}, braid/.style={black}] 

\path (0,0) node (bb) {}; 

\newcommand{\qedge}{0.6cm} 

\draw[->] (1.38,0) -- (1.62,0);
\draw[->] (4.38,0) -- (4.62,0);

\filldraw [draw=black, fill=white] (1,0) circle (0.03);
\filldraw [draw=black, fill=white] (0,1) circle (0.03);
\filldraw [draw=black, fill=white] (-1,0) circle (0.03);
\filldraw [draw=black, fill=white] (0,-1) circle (0.03);
\filldraw [draw=black, fill=white] (0,0) circle (0.03);


\node (A1) at (-0.8,0) {$\bullet$};
\node [below] at (A1) {\small {$A$}}; 

\node (B1) at (0,0.5) {$\bullet$};
\node [above] at (B1) {\small {$B$}}; 

\node (C1) at (0.8,0) {$\bullet$};
\node [below] at (C1) {\small {$C$}}; 

\node (D1) at (0,-0.5) {$\bullet$};
\node [below] at (D1) {\small {$D$}}; 

\draw [braid] plot  coordinates {(A1) (B1)}[postaction=decorate, decoration={markings, mark= at position 0.15 with
\node [above] {\small $\pi_3$};}];

\draw [braid] plot  coordinates {(B1) (C1)}[postaction=decorate, decoration={markings, mark= at position 0.85 with
\node [above] {\small $\pi_4$};}];

\draw [braid] plot[smooth]  coordinates {(B1)  (-0.2,0) (D1)};
\node at (-0.33,0) {\small $\pi_1$};

\draw [braid] plot[smooth]  coordinates {(B1)  (0.2,0) (D1)};
\node at (0.36,0) {\small $\pi_2$};

\begin{pgfonlayer}{background}
\draw[arc] (0,0) circle (1);

\draw[arc] (0,-1) .. controls (0.8,-0.2) and (0.8,0.2) .. (0,1);
\draw[arc] (0,-1) .. controls (-0.8,-0.2) and (-0.8,0.2) .. (0,1);

\draw[arc] (0,-1) .. controls (0.3,-0.5) .. node[very near end,sloped] {$\vertbowtie$}(0,0);
\draw[arc] (0,-1) .. controls (-0.3,-0.5) .. 
(0,0);
\end{pgfonlayer}

\begin{scope}[shift={(0,-1.5)}]

\node (V1) at (0,0) {$1$};
\node (V3) at (-\qedge*0.8,-\qedge*0.8) {$3$};
\node (V4) at (\qedge*0.8,-\qedge*0.8) {$4$};
\node (V2) at (0,-\qedge*1.6) {$2$};

\draw[->] (V1) -- (V3);
\draw[->] (V3) -- (V4);
\draw[->] (V4) -- (V1);
\draw[->] (V4) -- (V2);
\draw[->] (V2) -- (V3);
\end{scope}

\begin{scope}[shift={(3,0)}]
\filldraw [draw=black, fill=white] (1,0) circle (0.03);
\filldraw [draw=black, fill=white] (0,1) circle (0.03);
\filldraw [draw=black, fill=white] (-1,0) circle (0.03);
\filldraw [draw=black, fill=white] (0,-1) circle (0.03);
\filldraw [draw=black, fill=white] (0,0) circle (0.03);

\node (A2) at (-0.8,0) {$\bullet$};
\node [below] at (A2) {\small {$A$}}; 

\node (B2) at (0,0.5) {$\bullet$};
\node [above] at (B2) {\small {$B$}}; 

\node (C2) at (0.8,0) {$\bullet$};
\node [below] at (C2) {\small {$C$}}; 

\node (D2) at (0,-0.5) {$\bullet$};
\node [below] at (D2) {\small {$D$}};

\draw [braid] plot  coordinates {(A2) (D2)}[postaction=decorate, decoration={markings, mark= at position 0.55 with
\node [below] {\small $\widetilde{\pi_3}$};}];

\draw [braid] plot  coordinates {(B2) (C2)}[postaction=decorate, decoration={markings, mark= at position 0.85 with \node [above] {\small $\widetilde{\pi_4}$};}];

\draw [braid] plot[smooth]  coordinates {(B2)  (-0.2,0) (D2)};
\node at (-0.33,0) {\small $\widetilde{\pi_1}$};

\draw [braid] plot[smooth]  coordinates {(B2)  (0.2,0) (D2)};
\node at (0.36,0) {\small $\widetilde{\pi_2}$};

\begin{pgfonlayer}{background}
\draw[arc] (0,0) circle (1);

\draw[arc] (0,-1) .. controls (0.8,-0.2) and (0.8,0.2) .. (0,1);
\draw[arc] (0,-1) .. controls (-0.8,-0.2) and (-0.8,0.2) .. (0,1);

\draw[arc] (0,-1) .. controls (0.3,-0.5) .. node[very near end,sloped] {$\vertbowtie$}(0,0);
\draw[arc] (0,-1) .. controls (-0.3,-0.5) .. 
(0,0);
\end{pgfonlayer}

\end{scope}

\begin{scope}[shift={(6,0)}]
\filldraw [draw=black, fill=white] (1,0) circle (0.03);
\filldraw [draw=black, fill=white] (0,1) circle (0.03);
\filldraw [draw=black, fill=white] (-1,0) circle (0.03);
\filldraw [draw=black, fill=white] (0,-1) circle (0.03);
\filldraw [draw=black, fill=white] (0,0) circle (0.03);

\node (A3) at (-0.8,0) {$\bullet$};
\node [below] at (A3) {\small {$A$}}; 

\node (B3) at (0.3,0) {$\bullet$};
\node [above] at (B3) {\small {$B$}}; 

\node (C3) at (0.8,0) {$\bullet$};
\node [below] at (C3) {\small {$C$}}; 

\node (D3) at (-0.3,0) {$\bullet$};
\node [above] at (D3) {\small {$D$}}; 

\draw [braid] plot  coordinates {(A3) (D3)}[postaction=decorate, decoration={markings, mark= at position 0.6 with
\node [below] {\small $\rho_3$};}];

\draw [braid] plot  coordinates {(B3) (C3)}[postaction=decorate, decoration={markings, mark= at position 0.6 with \node [above] {\small $\rho_4$};}];

\draw [braid] plot[smooth]  coordinates {(B3)  (0,0.4) (D3)};
\node at (-0.13,0.47) {\small $\rho_1$};

\draw [braid] plot[smooth]  coordinates {(B3)  (0,-0.4) (D3)};
\node at (0.15,-0.48) {\small $\rho_2$};

\begin{pgfonlayer}{background}
\draw[arc] (0,0) circle (1);

\draw[arc] (0,-1) .. controls (0.8,-0.2) and (0.8,0.2) .. (0,1);
\draw[arc] (0,-1) .. controls (-0.8,-0.2) and (-0.8,0.2) .. (0,1);

\draw[arc] (0,-1) -- node[near end,sloped] {$\vertbowtie$} (0,0);

\draw[arc] (0,1) -- node[near end,sloped] {$\vertbowtie$} (0,0);
\end{pgfonlayer}

\begin{scope}[shift={(0,-1.5)}]

\node (V1) at (0,0) {$1$};
\node (V3) at (-\qedge*0.8,-\qedge*0.8) {$3$};
\node (V4) at (\qedge*0.8,-\qedge*0.8) {$4$};
\node (V2) at (0,-\qedge*1.6) {$2$};

\draw[->] (V3) -- (V1);
\draw[->] (V1) -- (V4);
\draw[->] (V4) -- (V2);
\draw[->] (V2) -- (V3);
\end{scope}

\end{scope}

\end{tikzpicture}
$$
\caption{Flip (at $\alpha_1$) inside a puzzle piece of type II, first case}
\label{f:flip2}
\end{figure}
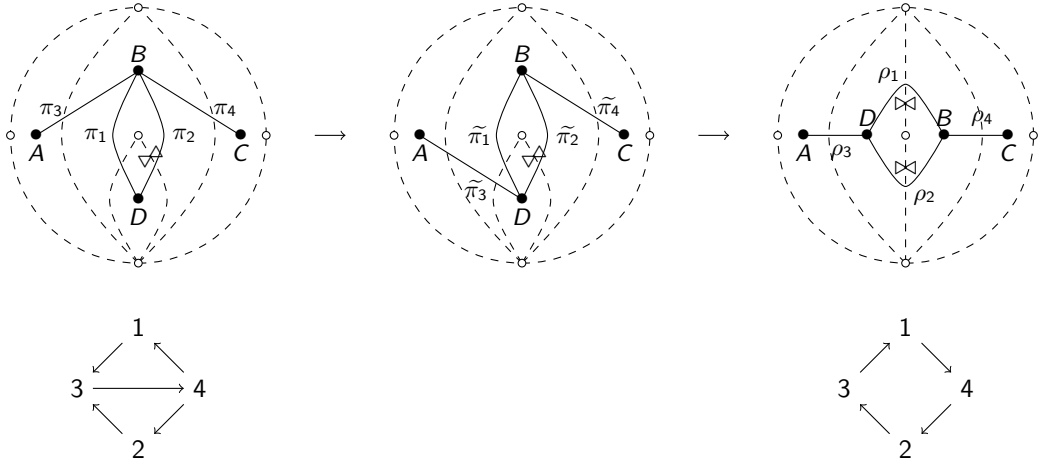

\begin{figure}
$$
\begin{tikzpicture}[scale=1.7,baseline=(bb.base),
  quivarrow/.style={black, -latex},
  arc/.style={black,dashed}, braid/.style={black}] 
\path (0,0) node (bb) {}; 

\newcommand{\qedge}{0.6cm} 

\draw[->] (1.38,0) -- (1.62,0);
\draw[->] (4.38,0) -- (4.62,0);

\filldraw [draw=black, fill=white] (1,0) circle (0.03);
\filldraw [draw=black, fill=white] (0,1) circle (0.03);
\filldraw [draw=black, fill=white] (-1,0) circle (0.03);
\filldraw [draw=black, fill=white] (0,-1) circle (0.03);
\filldraw [draw=black, fill=white] (0,0) circle (0.03);


\node (A1) at (-0.8,0) {$\bullet$};
\node [below] at (A1) {\small {$A$}}; 

\node (B1) at (0,0.5) {$\bullet$};
\node [above] at (B1) {\small {$B$}}; 

\node (C1) at (0.8,0) {$\bullet$};
\node [below] at (C1) {\small {$C$}}; 

\node (D1) at (0,-0.5) {$\bullet$};
\node [below] at (D1) {\small {$D$}}; 

\draw [braid] plot  coordinates {(A1) (B1)}[postaction=decorate, decoration={markings, mark= at position 0.15 with
\node [above] {\small $\pi_3$};}];

\draw [braid] plot  coordinates {(B1) (C1)}[postaction=decorate, decoration={markings, mark= at position 0.85 with
\node [above] {\small $\pi_4$};}];

\draw [braid] plot[smooth]  coordinates {(B1)  (-0.2,0) (D1)};
\node at (-0.33,0) {\small $\pi_1$};

\draw [braid] plot[smooth]  coordinates {(B1)  (0.2,0) (D1)};
\node at (0.36,0) {\small $\pi_2$};

\begin{pgfonlayer}{background}
\draw[arc] (0,0) circle (1);

\draw[arc] (0,-1) .. controls (0.8,-0.2) and (0.8,0.2) .. (0,1);
\draw[arc] (0,-1) .. controls (-0.8,-0.2) and (-0.8,0.2) .. (0,1);

\draw[arc] (0,-1) .. controls (0.3,-0.5) .. node[very near end,sloped] {$\vertbowtie$}(0,0);
\draw[arc] (0,-1) .. controls (-0.3,-0.5) .. 
(0,0);
\end{pgfonlayer}

\begin{scope}[shift={(0,-1.5)}]

\node (V1) at (0,0) {$1$};
\node (V3) at (-\qedge*0.8,-\qedge*0.8) {$3$};
\node (V4) at (\qedge*0.8,-\qedge*0.8) {$4$};
\node (V2) at (0,-\qedge*1.6) {$2$};

\draw[->] (V1) -- (V3);
\draw[->] (V3) -- (V4);
\draw[->] (V4) -- (V1);
\draw[->] (V4) -- (V2);
\draw[->] (V2) -- (V3);
\end{scope}

\begin{scope}[shift={(3,0)}]
\filldraw [draw=black, fill=white] (1,0) circle (0.03);
\filldraw [draw=black, fill=white] (0,1) circle (0.03);
\filldraw [draw=black, fill=white] (-1,0) circle (0.03);
\filldraw [draw=black, fill=white] (0,-1) circle (0.03);
\filldraw [draw=black, fill=white] (0,0) circle (0.03);

\node (A2) at (-0.8,0) {$\bullet$};
\node [below] at (A2) {\small {$A$}}; 

\node (B2) at (0,0.5) {$\bullet$};
\node [above] at (B2) {\small {$B$}}; 

\node (C2) at (0.8,0) {$\bullet$};
\node [below] at (C2) {\small {$C$}}; 

\node (D2) at (0,-0.5) {$\bullet$};
\node [below] at (D2) {\small {$D$}};

\draw [braid] plot[smooth]  coordinates {(A2) (-0.4,0.23) (0,0.25) (0.1,0) (D2)};

\node at (-0.67,0.26) {\small $\widetilde{\pi_3}$};

\draw [braid] plot  coordinates {(B2) (C2)}[postaction=decorate, decoration={markings, mark= at position 0.85 with \node [above] {\small $\widetilde{\pi_4}$};}];

\draw [braid] plot[smooth]  coordinates {(B2)  (-0.2,0) (D2)};
\node at (-0.33,0) {\small $\widetilde{\pi_1}$};

\draw [braid] plot[smooth]  coordinates {(B2)  (0.2,0) (D2)};
\node at (0.36,0) {\small $\widetilde{\pi_2}$};

\begin{pgfonlayer}{background}
\draw[arc] (0,0) circle (1);

\draw[arc] (0,-1) .. controls (0.8,-0.2) and (0.8,0.2) .. (0,1);
\draw[arc] (0,-1) .. controls (-0.8,-0.2) and (-0.8,0.2) .. (0,1);

\draw[arc] (0,-1) .. controls (0.3,-0.5) .. node[very near end,sloped] {$\vertbowtie$}(0,0);
\draw[arc] (0,-1) .. controls (-0.3,-0.5) .. 
(0,0);
\end{pgfonlayer}
\end{scope}

\begin{scope}[shift={(6,0)}]
\filldraw [draw=black, fill=white] (1,0) circle (0.03);
\filldraw [draw=black, fill=white] (0,1) circle (0.03);
\filldraw [draw=black, fill=white] (-1,0) circle (0.03);
\filldraw [draw=black, fill=white] (0,-1) circle (0.03);
\filldraw [draw=black, fill=white] (0,0) circle (0.03);

\node (A3) at (-0.8,0) {$\bullet$};
\node [below] at (A3) {\small {$A$}}; 

\node (B3) at (0.3,0) {$\bullet$};
\node [above] at (B3) {\small {$B$}}; 

\node (C3) at (0.8,0) {$\bullet$};
\node [below] at (C3) {\small {$C$}}; 

\node (D3) at (-0.3,0) {$\bullet$};
\node [above] at (D3) {\small {$D$}}; 

\draw [braid] plot  coordinates {(A3) (D3)}[postaction=decorate, decoration={markings, mark= at position 0.6 with
\node [below] {\small $\rho_3$};}];

\draw [braid] plot  coordinates {(B3) (C3)}[postaction=decorate, decoration={markings, mark= at position 0.6 with \node [above] {\small $\rho_4$};}];

\draw [braid] plot[smooth]  coordinates {(B3)  (0,0.4) (D3)};
\node at (-0.13,0.47) {\small $\rho_2$};

\draw [braid] plot[smooth]  coordinates {(D3)  (-0.1,0.22) (0,0.2) (0.15,0) (0,-0.2) (-0.15,0) (0,0.2) (0.1,0.22) (B3)};
\node at (0.15,-0.25) {\small $\rho_1$};

\begin{pgfonlayer}{background}
\draw[arc] (0,0) circle (1);

\draw[arc] (0,-1) .. controls (0.8,-0.2) and (0.8,0.2) .. (0,1);
\draw[arc] (0,-1) .. controls (-0.8,-0.2) and (-0.8,0.2) .. (0,1);

\draw[arc] (0,-1) -- (0,0);

\draw[arc] (0,1) --  (0,0);
\end{pgfonlayer}

\begin{scope}[shift={(0,-1.5)}]

\node (V1) at (0,0) {$1$};
\node (V3) at (-\qedge*0.8,-\qedge*0.8) {$3$};
\node (V4) at (\qedge*0.8,-\qedge*0.8) {$4$};
\node (V2) at (0,-\qedge*1.6) {$2$};

\draw[->] (V3) -- (V2);
\draw[->] (V2) -- (V4);
\draw[->] (V4) -- (V1);
\draw[->] (V1) -- (V3);
\end{scope}

\end{scope}
\end{tikzpicture}
$$
\caption{Flip (at $\alpha_2$) inside a puzzle piece of type II, second case}
\label{f:flip3}
\end{figure}
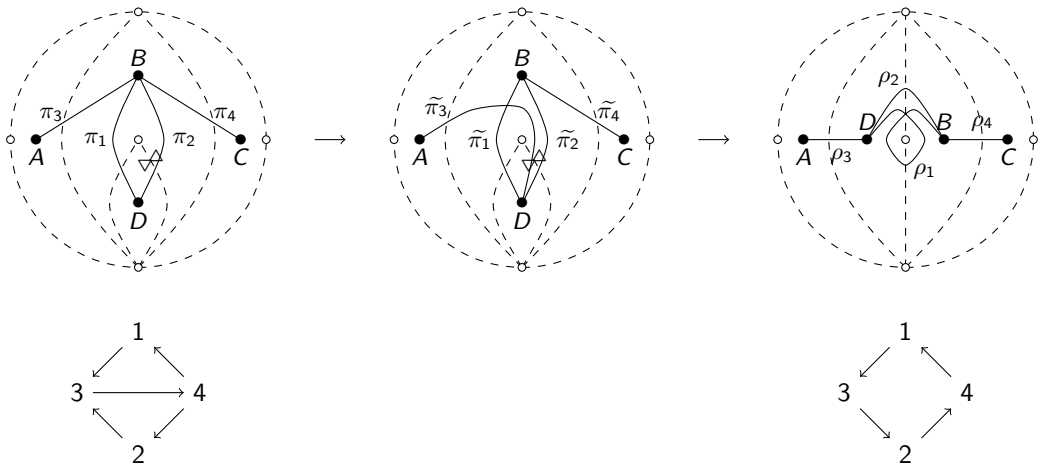


Note that the adjoining type I puzzle pieces (in
Figures~\ref{f:flip2} and~\ref{f:flip3}) may not occur, but the
argument is easily modified to cover these cases.
We also need to consider the flips from the
right hand diagram in each case to the
corresponding left hand one. We omit the
details: a similar argument can be applied in these cases.

Finally, we need to consider a flip involving an
arc where a puzzle piece of type I and a puzzle
piece of type II have been glued together.
These cases are shown in Figures~\ref{f:flip4}
and~\ref{f:flip5}: Figure~\ref{f:flip4} illustrates
the case where the puzzle piece of type I is
on the left of the puzzle piece of type II (when
it is drawn as shown), while Figure~\ref{f:flip5} illustrates the case where it is on the right. Again, a similar argument applies in the case of flips from the right hand diagram to the left hand one in these cases.
\end{pf}
\end{thm}

\begin{figure}
$$
\begin{tikzpicture}[scale=1.7,baseline=(bb.base),
  quivarrow/.style={black, -latex},
  arc/.style={black,dashed}, braid/.style={black}] 
\path (0,0) node (bb) {}; 

\newcommand{\qedge}{0.6cm} 

\draw[->] (1.38,0) -- (1.62,0);
\draw[->] (4.38,0) -- (4.62,0);

\filldraw [draw=black, fill=white] (1,0) circle (0.03);
\filldraw [draw=black, fill=white] (0,1) circle (0.03);
\filldraw [draw=black, fill=white] (-1,0) circle (0.03);
\filldraw [draw=black, fill=white] (0,-1) circle (0.03);
\filldraw [draw=black, fill=white] (0,0) circle (0.03);

\node (A1) at (-0.8,0) {$\bullet$};
\node [below] at (A1) {\small {$A$}}; 

\node (B1) at (0,0.5) {$\bullet$};
\node [above] at (B1) {\small {$B$}}; 

\node (C1) at (0.8,0) {$\bullet$};
\node [below] at (C1) {\small {$C$}}; 

\node (D1) at (0,-0.5) {$\bullet$};
\node [below] at (D1) {\small {$D$}}; 

\node (E1) at (-1.2,0.6) {$\bullet$};
\node [above] at (E1) {\small {$E$}};

\node (F1) at (-1.2,-0.6) {$\bullet$};
\node [below] at (F1) {\small {$F$}};

\draw [braid] plot  coordinates {(A1) (B1)}[postaction=decorate, decoration={markings, mark= at position 0.15 with
\node [above] {\small $\pi_3$};}];

\draw [braid] plot  coordinates {(B1) (C1)}[postaction=decorate, decoration={markings, mark= at position 0.85 with
\node [above] {\small $\pi_4$};}];

\draw [braid] plot[smooth]  coordinates {(B1)  (-0.2,0) (D1)};
\node at (-0.33,0) {\small $\pi_1$};

\draw [braid] plot[smooth]  coordinates {(B1)  (0.2,0) (D1)};
\node at (0.36,0) {\small $\pi_2$};

\draw [braid] plot coordinates
{(E1) (A1)};
\node at (-1.1,0.26) {\small $\pi_5$};

\draw [braid] plot coordinates
{(F1) (A1)};
\node at (-1.1,-0.26) {\small $\pi_6$};

\begin{pgfonlayer}{background}
\draw[arc] (0,0) circle (1);

\draw[arc] (0,-1) .. controls (0.8,-0.2) and (0.8,0.2) .. (0,1);
\draw[arc] (0,-1) .. controls (-0.8,-0.2) and (-0.8,0.2) .. (0,1);

\draw[arc] (0,-1) .. controls (0.3,-0.5) .. node[very near end,sloped] {$\vertbowtie$}(0,0);
\draw[arc] (0,-1) .. controls (-0.3,-0.5) .. 
(0,0);
\end{pgfonlayer}

\begin{scope}[shift={(\qedge*0.5,-1.7)}]

\node (V1) at (0,0) {$1$};
\node (V3) at (-\qedge*0.8,-\qedge*0.8) {$3$};
\node (V4) at (\qedge*0.8,-\qedge*0.8) {$4$};
\node (V2) at (0,-\qedge*1.6) {$2$};
\node (V5) at (-\qedge*1.6,0) {$5$};
\node (V6) at (-\qedge*1.6,-\qedge*1.6) {$6$};

\draw[->] (V1) -- (V3);
\draw[->] (V3) -- (V4);
\draw[->] (V4) -- (V1);
\draw[->] (V4) -- (V2);
\draw[->] (V2) -- (V3);
\draw[->] (V5) -- (V3);
\draw[->] (V3) -- (V6);
\draw[->] (V6) -- (V5);

\end{scope}

\begin{scope}[shift={(3,0)}]
\filldraw [draw=black, fill=white] (1,0) circle (0.03);
\filldraw [draw=black, fill=white] (0,1) circle (0.03);
\filldraw [draw=black, fill=white] (-1,0) circle (0.03);
\filldraw [draw=black, fill=white] (0,-1) circle (0.03);
\filldraw [draw=black, fill=white] (0,0) circle (0.03);

\node (A2) at (-0.8,0) {$\bullet$};
\node [below] at (A2) {\small {$A$}}; 

\node (B2) at (0,0.5) {$\bullet$};
\node [above] at (B2) {\small {$B$}}; 

\node (C2) at (0.8,0) {$\bullet$};
\node [below] at (C2) {\small {$C$}}; 

\node (D2) at (0,-0.5) {$\bullet$};
\node [below] at (D2) {\small {$D$}}; 

\node (E2) at (-1.2,0.6) {$\bullet$};
\node [above] at (E2) {\small {$E$}};

\node (F2) at (-1.2,-0.6) {$\bullet$};
\node [below] at (F2) {\small {$F$}};

\draw [braid] plot  coordinates {(A2) (B2)}[postaction=decorate, decoration={markings, mark= at position 0.6 with
\node [above] {\small $\widetilde{\pi_3}$};}];

\draw [braid] plot[smooth]  coordinates {(C2) (0.4,0.7) (0,0.85) (-0.4,0.7) (A2)};

\node at (0.8,0.27) {\small $\widetilde{\pi_4}$};
[postaction=decorate, decoration={markings, mark= at position 0.15 with \node [right] {\small $\widetilde{\pi_4}$};}];

\draw [braid] plot[smooth]  coordinates {(B2)  (-0.2,0) (D2)};
\node at (-0.33,0) {\small $\widetilde{\pi_1}$};

\draw [braid] plot[smooth]  coordinates {(B2)  (0.2,0) (D2)};
\node at (0.36,0) {\small $\widetilde{\pi_2}$};

\draw [braid] plot coordinates
{(E2) (A2)};
\node at (-1.1,0.26) {\small $\widetilde{\pi_5}$};

\draw [braid] plot coordinates
{(F2) (B2)};
\node at (-0.7,-0.38) {\small $\widetilde{\pi_6}$};

\begin{pgfonlayer}{background}
\draw[arc] (0,0) circle (1);

\draw[arc] (0,-1) .. controls (0.8,-0.2) and (0.8,0.2) .. (0,1);
\draw[arc] (0,-1) .. controls (-0.8,-0.2) and (-0.8,0.2) .. (0,1);

\draw[arc] (0,-1) .. controls (0.3,-0.5) .. node[very near end,sloped] {$\vertbowtie$}(0,0);
\draw[arc] (0,-1) .. controls (-0.3,-0.5) .. 
(0,0);
\end{pgfonlayer}
\end{scope}

\begin{scope}[shift={(6,0)}]
\filldraw [draw=black, fill=white] (1,0) circle (0.03);
\filldraw [draw=black, fill=white] (0,1) circle (0.03);
\filldraw [draw=black, fill=white] (-1,0) circle (0.03);
\filldraw [draw=black, fill=white] (0,-1) circle (0.03);
\filldraw [draw=black, fill=white] (0,0) circle (0.03);

\node (A3) at (-0.4,0.6) {$\bullet$};
\node [above] at (A3) {\small {$A$}}; 

\node (B3) at (-0.5,-0.3) {$\bullet$};
\node [above left] at (B3) {\small {$B$}}; 

\node (C3) at (0.8,0) {$\bullet$};
\node [below] at (C3) {\small {$C$}}; 

\node (D3) at (0,-0.5) {$\bullet$};
\node [below] at (D3) {\small {$D$}}; 

\node (E3) at (-1.2,0.6) {$\bullet$};
\node [above] at (E3) {\small {$E$}};

\node (F3) at (-1.2,-0.6) {$\bullet$};
\node [below] at (F3) {\small {$F$}};

\draw [braid] plot coordinates {(E3) (A3)};
\node at (-0.94,0.7) {\small $\rho_5$};

\draw [braid] plot  coordinates {(F3) (B3)};
\node at (-0.9,-0.66) {\small $\rho_6$};

\draw [braid] plot  coordinates {(B3)  (A3)};
\node at (-0.53,0.4) {\small $\rho_3$};

\draw [braid] plot[smooth]  coordinates {(A3) (0,0.8) (0.5,0.7) (C3)};
\node at (0.75,0.4) {\small $\rho_4$};

\draw [braid] plot[smooth]  coordinates {(B3) (0,0.15) (0.15,0)(D3)};
\node at (0.27,-0.2) {\small $\rho_2$};

\draw [braid] plot  coordinates {(B3)(D3)};
\node at (-0.3,-0.47) {\small $\rho_1$};

\begin{pgfonlayer}{background}
\draw[arc] (0,0) circle (1);

\draw[arc] (0,-1) .. controls (0.8,-0.2) and (0.8,0.2) .. (0,1);

\draw[arc] (0,-1) .. controls (0.3,-0.5) .. (0,0);
\draw[arc] (0,-1) .. controls (-0.3,-0.5) .. 
(0,0);

\draw [arc] plot[smooth]  coordinates
{(-1,0) (0,0.25) (0.4,0) (0.4,-0.3) (0.3,-0.6)
 (0,-1)};

\end{pgfonlayer}

\begin{scope}[shift={(\qedge*0.5,-1.7)}]

\node (V1) at (0,0) {$1$};
\node (V3) at (-\qedge*0.8,-\qedge*0.8) {$3$};
\node (V4) at (\qedge*0.8,-\qedge*0.8) {$4$};
\node (V2) at (0,-\qedge*1.6) {$2$};
\node (V5) at (-\qedge*1.6,0) {$5$};
\node (V6) at (-\qedge*1.6,-\qedge*1.6) {$6$};

\draw[->] (V3) -- (V5);
\draw[->] (V3) -- (V1);
\draw[->] (V3) -- (V2);
\draw[->] (V6) -- (V3);
\draw[->] (V4) -- (V3);
\draw[->] (V2) -- (V6);

\draw[->] (V5) to[out=30,in=90] (V4);
\draw[->] (V1) to[out=180,in=70] (V6);
\end{scope}

\end{scope}
\end{tikzpicture}
$$
\caption{Flip involving an arc ($\alpha_3$) where puzzle pieces of type I and II are glued, first case}
\label{f:flip4}
\end{figure}
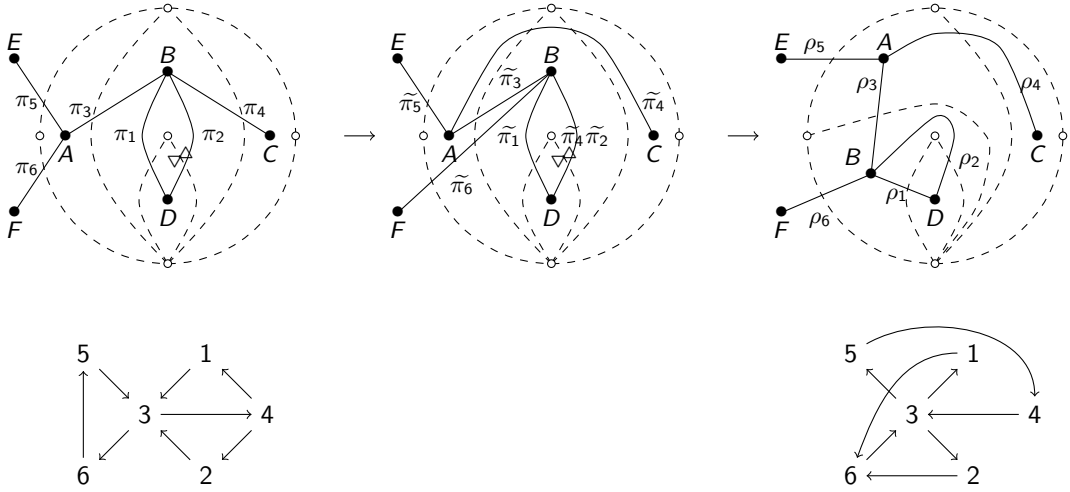

\begin{figure}
$$
\begin{tikzpicture}[scale=1.7,baseline=(bb.base),
  quivarrow/.style={black, -latex},
  arc/.style={black,dashed}, braid/.style={black}] 
\path (0,0) node (bb) {}; 

\newcommand{\qedge}{0.6cm} 

\draw[->] (1.38,0) -- (1.62,0);
\draw[->] (4.38,0) -- (4.62,0);

\filldraw [draw=black, fill=white] (1,0) circle (0.03);
\filldraw [draw=black, fill=white] (0,1) circle (0.03);
\filldraw [draw=black, fill=white] (-1,0) circle (0.03);
\filldraw [draw=black, fill=white] (0,-1) circle (0.03);
\filldraw [draw=black, fill=white] (0,0) circle (0.03);

\node (A1) at (-0.8,0) {$\bullet$};
\node [below] at (A1) {\small {$A$}}; 

\node (B1) at (0,0.5) {$\bullet$};
\node [above] at (B1) {\small {$B$}}; 

\node (C1) at (0.8,0) {$\bullet$};
\node [below] at (C1) {\small {$C$}}; 

\node (D1) at (0,-0.5) {$\bullet$};
\node [below] at (D1) {\small {$D$}}; 

\node (E1) at (1.2,0.6) {$\bullet$};
\node [above] at (E1) {\small {$E$}};

\node (F1) at (1.2,-0.6) {$\bullet$};
\node [below] at (F1) {\small {$F$}};

\draw [braid] plot  coordinates {(A1) (B1)}[postaction=decorate, decoration={markings, mark= at position 0.15 with
\node [above] {\small $\pi_3$};}];

\draw [braid] plot  coordinates {(B1) (C1)}[postaction=decorate, decoration={markings, mark= at position 0.85 with
\node [above] {\small $\pi_4$};}];

\draw [braid] plot[smooth]  coordinates {(B1)  (-0.2,0) (D1)};
\node at (-0.33,0) {\small $\pi_1$};

\draw [braid] plot[smooth]  coordinates {(B1)  (0.2,0) (D1)};
\node at (0.36,0) {\small $\pi_2$};

\draw [braid] plot coordinates
{(E1) (C1)};
\node at (1.1,0.26) {\small $\pi_5$};

\draw [braid] plot coordinates
{(F1) (C1)};
\node at (1.1,-0.26) {\small $\pi_6$};

\begin{pgfonlayer}{background}
\draw[arc] (0,0) circle (1);

\draw[arc] (0,-1) .. controls (0.8,-0.2) and (0.8,0.2) .. (0,1);
\draw[arc] (0,-1) .. controls (-0.8,-0.2) and (-0.8,0.2) .. (0,1);

\draw[arc] (0,-1) .. controls (0.3,-0.5) .. node[very near end,sloped] {$\vertbowtie$}(0,0);
\draw[arc] (0,-1) .. controls (-0.3,-0.5) .. 
(0,0);
\end{pgfonlayer}

\begin{scope}[shift={(-\qedge*0.5,-1.7)}]

\node (V1) at (0,0) {$1$};
\node (V3) at (-\qedge*0.8,-\qedge*0.8) {$3$};
\node (V4) at (\qedge*0.8,-\qedge*0.8) {$4$};
\node (V2) at (0,-\qedge*1.6) {$2$};
\node (V5) at (\qedge*1.6,0) {$5$};
\node (V6) at (\qedge*1.6,-\qedge*1.6) {$6$};

\draw[->] (V1) -- (V3);
\draw[->] (V3) -- (V4);
\draw[->] (V4) -- (V1);
\draw[->] (V4) -- (V2);
\draw[->] (V2) -- (V3);
\draw[->] (V4) -- (V5);
\draw[->] (V6) -- (V4);
\draw[->] (V5) -- (V6);

\end{scope}

\begin{scope}[shift={(3,0)}]
\filldraw [draw=black, fill=white] (1,0) circle (0.03);
\filldraw [draw=black, fill=white] (0,1) circle (0.03);
\filldraw [draw=black, fill=white] (-1,0) circle (0.03);
\filldraw [draw=black, fill=white] (0,-1) circle (0.03);
\filldraw [draw=black, fill=white] (0,0) circle (0.03);

\node (A2) at (-0.8,0) {$\bullet$};
\node [below] at (A2) {\small {$A$}}; 

\node (B2) at (0,0.5) {$\bullet$};
\node [above] at (B2) {\small {$B$}}; 

\node (C2) at (0.8,0) {$\bullet$};
\node [below] at (C2) {\small {$C$}}; 

\node (D2) at (0,-0.5) {$\bullet$};
\node [below] at (D2) {\small {$D$}}; 

\node (E2) at (1.2,0.6) {$\bullet$};
\node [above] at (E2) {\small {$E$}};

\node (F2) at (1.2,-0.6) {$\bullet$};
\node [below] at (F2) {\small {$F$}};

\draw [braid] plot  coordinates {(A2) (B2)}[postaction=decorate, decoration={markings, mark= at position 0.2 with
\node [above] {\small $\widetilde{\pi_3}$};}];

\draw [braid] plot  coordinates {(C2) (B2)};
\node at (0.4,0.37) {\small $\widetilde{\pi_4}$};

\draw [braid] plot  coordinates {(B2)  (E2)};
\node at (0.5,0.65) {\small $\widetilde{\pi_6}$};

\draw [braid] plot  coordinates {(C2) (F2)};
\node at (1.2,-0.4) {\small $\widetilde{\pi_5}$};

\draw [braid] plot[smooth] coordinates
{(D2) (-0.2,0) (0,0.2) (C2)};
\node at (-0.28,-0.15) {\small $\widetilde{\pi_1}$};

\draw [braid] plot coordinates
{(D2) (C2)};
\node at (0.38,-0.42) {\small $\widetilde{\pi_2}$};

\begin{pgfonlayer}{background}
\draw[arc] (0,0) circle (1);

\draw[arc] (0,-1) .. controls (0.8,-0.2) and (0.8,0.2) .. (0,1);
\draw[arc] (0,-1) .. controls (-0.8,-0.2) and (-0.8,0.2) .. (0,1);

\draw[arc] (0,-1) .. controls (0.3,-0.5) .. node[very near end,sloped] {$\vertbowtie$}(0,0);
\draw[arc] (0,-1) .. controls (-0.3,-0.5) .. 
(0,0);
\end{pgfonlayer}
\end{scope}

\begin{scope}[shift={(6,0)}]
\filldraw [draw=black, fill=white] (1,0) circle (0.03);
\filldraw [draw=black, fill=white] (0,1) circle (0.03);
\filldraw [draw=black, fill=white] (-1,0) circle (0.03);
\filldraw [draw=black, fill=white] (0,-1) circle (0.03);
\filldraw [draw=black, fill=white] (0,0) circle (0.03);

\node (B3) at (0.4,0.6) {$\bullet$};
\node [above] at (B3) {\small {$B$}}; 

\node (C3) at (0.5,-0.3) {$\bullet$};
\node [above left] at (C3) {\small {$C$}}; 

\node (A3) at (-0.8,0) {$\bullet$};
\node [below] at (A3) {\small {$A$}}; 

\node (D3) at (0,-0.5) {$\bullet$};
\node [below] at (D3) {\small {$D$}}; 

\node (E3) at (1.2,0.6) {$\bullet$};
\node [above] at (E3) {\small {$E$}};

\node (F3) at (1.2,-0.6) {$\bullet$};
\node [below] at (F3) {\small {$F$}};

\draw [braid] plot  coordinates {(E3) (B3)};
\node at (0.94,0.68) {\small $\rho_6$};

\draw [braid] plot  coordinates {(F3) (C3)};
\node at (0.92,-0.64) {\small $\rho_5$};

\draw [braid] plot  coordinates {(B3)  (C3)};
\node at (0.55,0.4) {\small $\rho_4$};

\draw [braid] plot[smooth]  coordinates {(B3) (0,0.8) (-0.5,0.7) (A3)};
\node at (-0.78,0.4) {\small $\rho_3$};

\draw [braid] plot[smooth]  coordinates {(C3) (0,0.15) (-0.15,0) (D3)};
\node at (-0.22,-0.2) {\small $\rho_1$};

\draw [braid] plot  coordinates {(C3) (D3)};
\node at (0.32,-0.49) {\small $\rho_2$};

\begin{pgfonlayer}{background}
\draw[arc] (0,0) circle (1);

\draw[arc] (0,-1) .. controls (-0.8,-0.2) and (-0.8,0.2) .. (0,1);

\draw[arc] (0,-1) .. controls (-0.3,-0.5) .. (0,0);
\draw[arc] (0,-1) .. controls (0.3,-0.5) .. 
(0,0);

\draw [arc] plot[smooth]  coordinates
{(1,0) (0,0.25) (-0.4,0) (-0.4,-0.3) (-0.3,-0.6) (0,-1)};

\end{pgfonlayer}

\begin{scope}[shift={(-\qedge*0.5,-1.7)}]

\node (V1) at (0,0) {$1$};
\node (V3) at (-\qedge*0.8,-\qedge*0.8) {$3$};
\node (V4) at (\qedge*0.8,-\qedge*0.8) {$4$};
\node (V2) at (0,-\qedge*1.6) {$2$};
\node (V5) at (\qedge*1.6,0) {$5$};
\node (V6) at (\qedge*1.6,-\qedge*1.6) {$6$};

\draw[->] (V1) -- (V4);
\draw[->] (V5) -- (V4);
\draw[->] (V2) -- (V4);
\draw[->] (V4) -- (V6);
\draw[->] (V4) -- (V3);
\draw[->] (V6) -- (V2);

\draw[->] (V6) to[out=110,in=-10] (V1);
\draw[->] (V5) to[out=135,in=80] (V3);
\end{scope}

\end{scope}
\end{tikzpicture}
$$
\caption{Flip involving an arc ($\alpha_4$) where puzzle pieces of type I and II are glued, second case}
\label{f:flip5}
\end{figure}
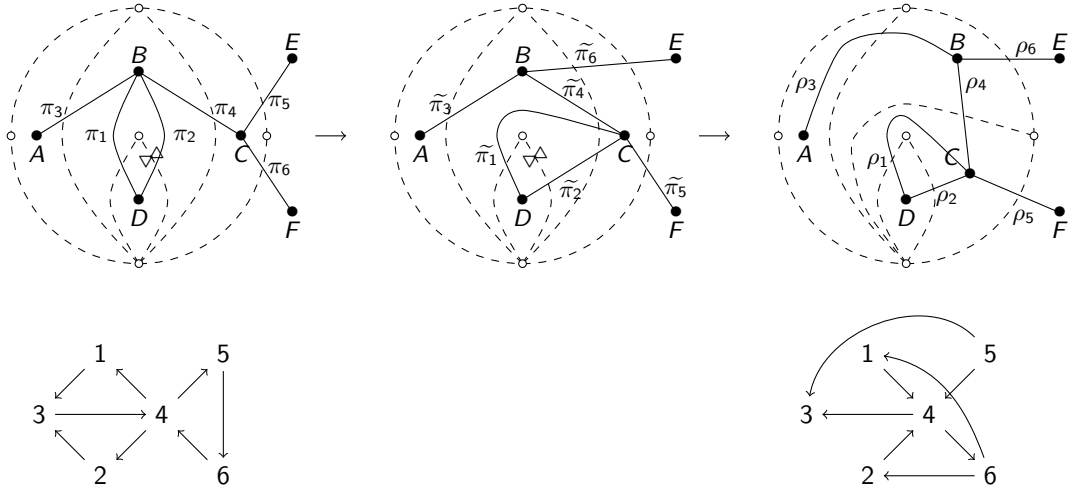






\section{Actions on categories}\label{sec:categ}
\subsection{Quivers with potential}\label{ss:qps}

Fix an algebraically closed field $\F$.  
To any quiver $Q$ we can associate the path algebra $\F Q$, which, as an $\F$-vector space, has basis given by all paths in $Q$ of length $\geq0$, and the multiplication of two paths $p_1$ and $p_2$ is their concatenation $p_1p_2$ if $p_1$ ends and $p_2$ starts at the same vertex, and is zero otherwise.

Let $\F Q_{\geq n}$ be the ideal of $\F Q$ generated by the paths in $Q$ of length at least $n$.  We can take the  completion $\widehat{\F Q}$ of $\F Q$ with respect to $\F Q_{\geq1}$, which is defined as follows:
$$\widehat{\F Q} = \varprojlim_n \frac{\F Q}{\F Q_{\geq n}}=\left\{(a_n+\F Q_{\geq n})_{n=1}^\infty\st a_n\in \F Q, \varphi_n(a_n+\F Q_{\geq n})=a_{n-1}+\F Q_{\geq n-1}\right\}$$
where the limit is taken along the chain of epimorphisms
$$ \frac{\F Q}{\F Q_{\geq 1}} \stackrel{\twoheadleftarrow}{\varphi_2} \frac{\F Q}{\F Q_{\geq 2}} \stackrel{\twoheadleftarrow}{\varphi_3} \frac{\F Q}{\F Q_{\geq 3}} \stackrel{\twoheadleftarrow}{}  \cdots$$

Let $\widehat{\F Q}_{\text{cyc}}$ denote the subspace of (possibly infinite) linear combinations of cycles in $Q$.
Recall that a \emph{potential} for a quiver $Q$ is an element
$W$ of $\widehat{\F Q}_{\text{cyc}}$,
regarded up to cyclic equivalence
(and for which no two cyclically equivalent paths in $Q$ occur in the decomposition of $W$). 
The pair
$(Q,W)$ is called a \emph{quiver with potential}~\cite{dwz1},
which we occasionally abbreviate to QP.
The following definition is {\cite[Definition 4.2]{dwz1}}.
\begin{defn}[Derksen-Weyman-Zelevinsky]
Let $Q_1$ and $Q_2$ be two quivers with the same vertex set $I$ and $(Q_1,W_1)$ and $(Q_2,W_2)$ be two QPs.  A \emph{right equivalence} between $(Q_1,W_1)$ and $(Q_2,W_2)$ is an algebra isomorphism $\varphi: \widehat{\F Q_1}\to \widehat{\F Q_2}$ such that $\varphi(W_1)$ is cyclically equivalent to $W_2$ and $\varphi$ is the identity when restricted to the semisimple subalgebra $\F I$ of $\widehat{\F Q_1}$.
\end{defn}

A quiver with potential $(Q,W)$ with $W$ containing paths
of length two or more is \emph{trivial} if $Q$ is a disjoint union of $2$-cycles
and there is an algebra automorphism
of $\widehat{kQ}$ preserving the span of the arrows of $Q$ (a \emph{change of arrows}) which takes $W$ to the sum of the $2$-cycles in $Q$. A quiver with potential $(Q,W)$
is said to be \emph{reduced} if $W$ is a
linear combination of cycles in $Q$ of length $3$ or more.

The \emph{splitting theorem}~\cite[Thm.\ 4.6]{dwz1} states that every quiver with potential can be
written as a direct sum of a reduced
quiver with potential and a trivial quiver with potential which are unique
up to right equivalence.

Let $(Q,W)$ be a quiver with potential, and let $k$ be a vertex of $Q$ not involved in any $2$-cycles. By replacing $W$ with a cyclically equivalent potential on $Q$
if necessary, we can assume that none
of the cycles in the decomposition of
$W$ start or end at $k$. We denote by
$\widetilde{\mu_k}(Q,W)$ the
non-reduced mutation of $(Q,W)$ at $k$ in $Q$, as defined in~\cite[\S5]{dwz1}.
Then the right equivalence class of $\widetilde{\mu_k}(Q,W)$ is determined by
the right equivalence class of $(Q,W)$
by~\cite[Thm.\ 5.2]{dwz1}. The \emph{mutation} $\mu_k(Q,W)$ of $(Q,W)$
at $k$ is then defined to be the reduced
component of $\widetilde{\mu_k}(Q,W)$,
and is uniquely determined up to right
equivalence, given the right equivalence class of $(Q,W)$.

As before, we will say that
a quiver with potential $(Q,W)$ is \emph{Dynkin} if the
underlying unoriented graph of $Q$ is
an orientation of a Dynkin quiver (and
hence $W=0$). We shall say that
a quiver with potential $(Q',W')$ is \emph{mutation-Dynkin}
if it can be obtained by repeatedly mutating a Dynkin quiver with potential in the above sense.
For the rest of Section~\ref{ss:qps}
\emph{we will restrict to Dynkin types
$A$ and $D$}.

Let $(S,M)$ be the Riemann surface with
marked points associated to $\Delta$ as
in Section \ref{sec:riemann}. So, if $\Delta=A_n$,
we take $S$ to be a disk with $n-3$ points on its boundary, and if $\Delta=D_n$, we take $S$ to be a disk with one marked point in its interior and $n$ marked points on its boundary.

\begin{figure}
$$
\begin{tikzpicture}
[scale=2,baseline=(bb.base), quivarrow/.style={black, -latex}, arc/.style={black,dashed}, braid/.style={black}] 

\path (0,0) node (bb) {}; 


\node at (-1,0.7) {(a)};

\begin{pgfonlayer}{background}
\draw[arc] (0,0) -- (1,0) -- (0.5,0.866) -- (0,0);
\end{pgfonlayer}

\filldraw [draw=black, fill=white] (0,0) circle (0.05);
\filldraw [draw=black, fill=white] (0.5,0.866) circle (0.05);
\filldraw [draw=black, fill=white] (1,0) circle (0.05);

\node (E1) at (0.1,0.43) {$1$};
\node (E2) at (0.9,0.43) {$2$};
\node (E3) at (0.5,-0.2) {$3$};

\begin{scope}[shift={(2,0.2)}]

\node (Q1) at (0.1,0.43) {$1$};
\node (Q2) at (0.9,0.43) {$2$};
\node (Q3) at (0.5,-0.2) {$3$};

\draw[->] (Q1) -- (Q2);
\draw[->] (Q2) -- (Q3);
\draw[->] (Q3) -- (Q1);

\node at (0.5,0.53) {$\lettera$};
\node at (0.8,0.05) {$\letterb$};
\node at (0.2,0.05) {$\letterc$};

\end{scope}

\begin{scope}[shift={(4,0.2)}]

\node at (0,0.6) {$W_{\T}:\,\lettera\letterb\letterc$};
\node at (0,0.2) {$W_{\T}':\,\lettera\letterb\letterc$};

\end{scope}


\begin{scope}[shift={(0.5,-2.7)}]

\node at (-1.5,1.7) {(b)};

\filldraw [draw=black, fill=white] (0,0) circle (0.05);
\filldraw [draw=black, fill=white] (0,1) circle (0.05);
\filldraw [draw=black, fill=white] (0,2) circle (0.05);

\node (F1) at (-0.25,0.7) {$1$};
\node (F2) at (0.25,0.7) {$2$};
\node (F3) at (-0.55,1) {$3$};
\node (F4) at (0.55,1) {$4$};

\begin{pgfonlayer}{background}
\draw[arc] (0,0) to[out=45,in=-45] (0,2);
\draw[arc] (0,0) to[out=135,in=-135] (0,2);
\draw[arc] (0,0) to[out=60,in=-60] node[near end,sloped] {$\vertbowtie$} (0,1);
\draw[arc] (0,0) to[out=120,in=-120]  (0,1);
\end{pgfonlayer}

\begin{scope}[shift={(1.6,1)}]
\node (R3) at (0,0) {$3$};
\node (R4) at (1,0) {$4$};
\node (R1) at (0.5,0.5) {$1$};
\node (R2) at (0.5,-0.5) {$2$};

\draw[->] (R1) -- (R3);
\draw[->] (R3) -- (R4);
\draw[->] (R4) -- (R2);
\draw[->] (R4) -- (R1);
\draw[->] (R2) -- (R3);

\node at (0.18,0.34) {$\lettera_2$};
\node at (0.83,0.33) {$\lettera_1$};
\node at (0.2,-0.4) {$\letterb_2$};
\node at (0.8,-0.4) {$\letterb_1$};
\node at (0.5,0.1) {$\letterc$};

\end{scope}

\begin{scope}[shift={(4,1)}]

\node at (0,0.6) {$W_{\T}:\,\lettera_1\lettera_2\letterc+\letterb_1\letterb_2\letterc$};
\node at (0,0.2) {$W_{\T}':\,\lettera_1\lettera_2\letterc+\letterb_1\letterb_2\letterc$};

\end{scope}
\end{scope}


\begin{scope}[shift={(0.35,-3.9)}]

\node at (-1.35,0.7) {(c)};

\begin{pgfonlayer}{background}
\draw[arc] (0:0.8) -- (72:0.8) -- (144:0.8) -- (216:0.8) -- (288:0.8) -- (0:0.8);
\node at (14:0.4) {$3$};
\node at (86:0.4) {$2$};
\node at (158:0.4) {$1$};
\node at (230:0.4) {$5$};
\node at (302:0.4) {$4$};

\draw[arc] (0,0) -- (0:0.8);
\draw[arc] (0,0) -- (72:0.8);
\draw[arc] (0,0) -- (144:0.8);
\draw[arc] (0,0) -- (216:0.8);
\draw[arc] (0,0) -- (288:0.8);
\end{pgfonlayer}

\filldraw [draw=black, fill=white] (0:0.8) circle (0.05);
\filldraw [draw=black, fill=white] (72:0.8) circle (0.05);
\filldraw [draw=black, fill=white] (144:0.8) circle (0.05);
\filldraw [draw=black, fill=white] (216:0.8) circle (0.05);
\filldraw [draw=black, fill=white] (288:0.8) circle (0.05);
\filldraw [draw=black, fill=white] (0,0) circle (0.05);

\begin{scope}[shift={(2.3,0.2)}]

\node (Q3) at (14:0.6) {$3$};
\node (Q2) at (86:0.6) {$2$};
\node (Q1) at (158:0.6) {$1$};
\node (Q5) at (230:0.6) {$5$};
\node (Q4) at (302:0.6) {$4$};

\draw[->] (Q1) -- (Q2);
\draw[->] (Q2) -- (Q3);
\draw[->] (Q3) -- (Q4);
\draw[->] (Q4) -- (Q5);
\draw[->] (Q5) -- (Q1);

\node at (50:0.66) {$\letterb$};
\node at (122:0.66) {$\lettera$};
\node at (194:0.66) {$\lettere$};
\node at (266:0.66) {$\letterd$};
\node at (338:0.66) {$\letterc$};

\end{scope}

\begin{scope}[shift={(4,0.2)}]

\node at (0,0.6) {$W_{\T}:\,\lettera\letterb\letterc\letterd\lettere$};
\node at (0,0.2) {$W_{\T}':\,-\lettera\letterb\letterc\letterd\lettere$};

\end{scope}

\end{scope}


\begin{scope}[shift={(0.35,-6)}]

\node at (-1.35,0.7) {(d)};

\begin{pgfonlayer}{background}
\draw[arc] (0:0.8) -- (72:0.8) -- (144:0.8) -- (216:0.8) -- (288:0.8) -- (0:0.8);
\node at (14:0.4) {$3$};
\node at (86:0.4) {$2$};
\node at (158:0.4) {$1$};
\node at (230:0.4) {$5$};
\node at (302:0.4) {$4$};

\draw[arc] (0,0) -- node[near start,sloped] {$\vertbowtie$} (0:0.8);
\draw[arc] (0,0) -- node[near start,sloped] {$\vertbowtie$} (72:0.8);
\draw[arc] (0,0) -- node[near start,sloped] {$\vertbowtie$} (144:0.8);
\draw[arc] (0,0) -- node[near start,sloped] {$\vertbowtie$} (216:0.8);
\draw[arc] (0,0) -- node[near start,sloped] {$\vertbowtie$} (288:0.8);
\end{pgfonlayer}

\filldraw [draw=black, fill=white] (0:0.8) circle (0.05);
\filldraw [draw=black, fill=white] (72:0.8) circle (0.05);
\filldraw [draw=black, fill=white] (144:0.8) circle (0.05);
\filldraw [draw=black, fill=white] (216:0.8) circle (0.05);
\filldraw [draw=black, fill=white] (288:0.8) circle (0.05);
\filldraw [draw=black, fill=white] (0,0) circle (0.05);

\begin{scope}[shift={(2.3,0.2)}]

\node (Q3) at (14:0.6) {$3$};
\node (Q2) at (86:0.6) {$2$};
\node (Q1) at (158:0.6) {$1$};
\node (Q5) at (230:0.6) {$5$};
\node (Q4) at (302:0.6) {$4$};

\draw[->] (Q1) -- (Q2);
\draw[->] (Q2) -- (Q3);
\draw[->] (Q3) -- (Q4);
\draw[->] (Q4) -- (Q5);
\draw[->] (Q5) -- (Q1);

\node at (50:0.66) {$\letterb$};
\node at (122:0.66) {$\lettera$};
\node at (194:0.66) {$\lettere$};
\node at (266:0.66) {$\letterd$};
\node at (338:0.66) {$\letterc$};

\end{scope}

\begin{scope}[shift={(4,0.2)}]

\node at (0,0.6) {$W_{\T}:\,\lettera\letterb\letterc\letterd\lettere$};
\node at (0,0.2) {$W_{\T}':\,\lettera\letterb\letterc\letterd\lettere$};
\end{scope}
\end{scope}

\end{tikzpicture}
$$
\caption{Terms in the potential $W_{\T}$}
\label{f:potentialterms}
\end{figure}
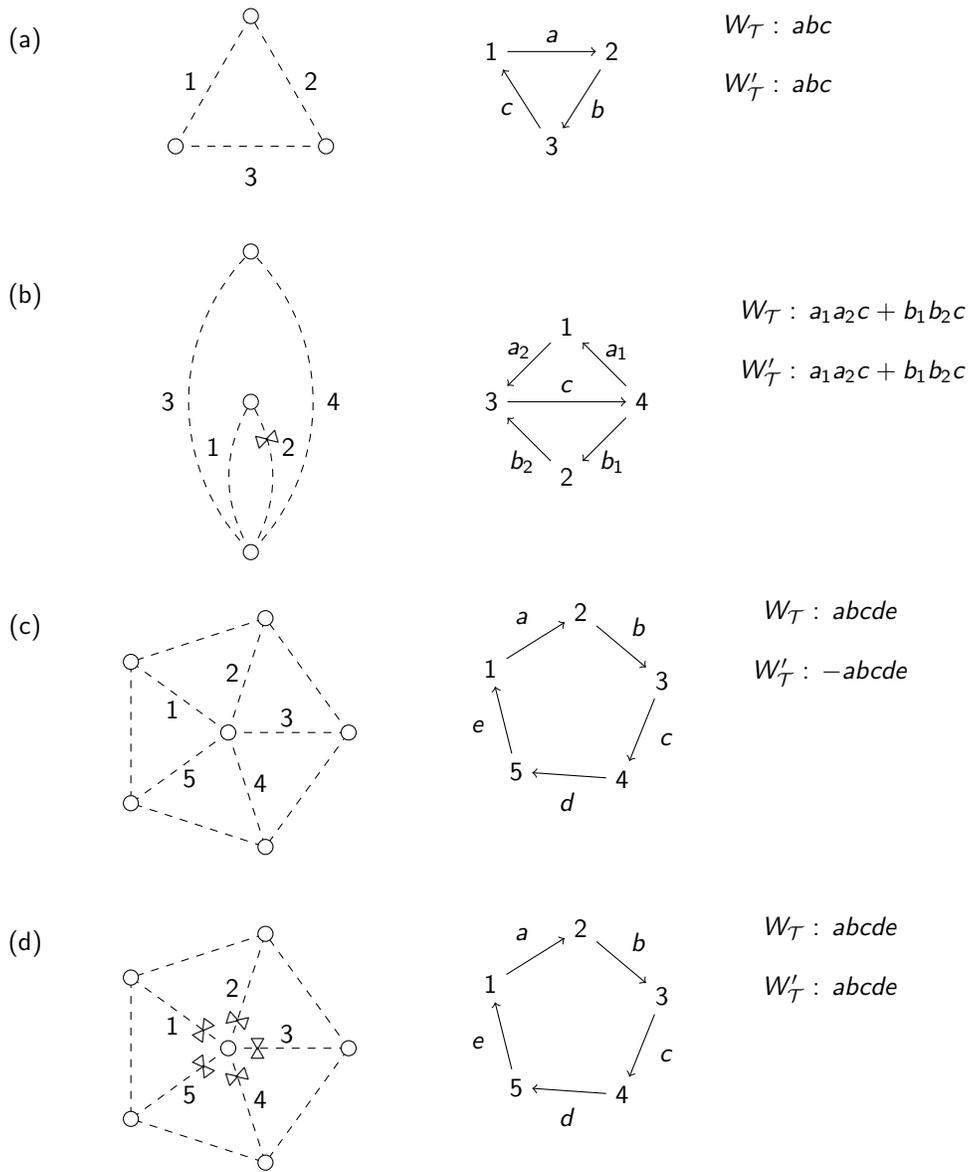

Let $Q$ be a mutation-Dynkin quiver.
By~\cite{fst08}, $Q=Q_{\T}$ for some tagged
triangulation $\T$ of $(S,M)$.
Let $W,W'$ be the sum of the terms coming from local configurations in $\T$
as shown in Figure~\ref{f:potentialterms}
(where in (c) and (d) there are at
least three arcs incident with the interior marked point).

Then $W_{\T}$ is the potential
given by taking the sum of the induced cycles in $Q_{\T}$ (i.e.\ induced subgraphs of
$Q_{\T}$ which are cycles),
and $W_{\T}'$ is the potential
associated to $\T$ in~\cite[\S3]{LF4},
taking the parameter associated to the
internal marked point (if there is one) to be equal to $-1$. Then we have the following:


\begin{lem} \label{l:rightequivalent}
The potentials $W_{\T}$ and $W_{\T}'$
are right equivalent.

\begin{pf}
We assume we are in case $D_n$, since
the two potentials coincide in case
$A_n$.
If the interior marked point is as in case (c) of Figure~\ref{f:potentialterms} (with
at least $3$ arcs incident with it), then
there is a unique triangle in $\T$
with sides $1$ and $2$.
We label the arrows in the corresponding
$3$-cycle in $W_{\T}$ or $W_{\T}'$
by $a,x,y$, in order around the cycle.
Then the automorphism $\varphi$ of
$\widehat{kQ_{\T}}$ negating $a$
and $x$ and taking each other arrow 
to itself gives a right equivalence between $W_{\T}$ and $W_{\T}'$, since $a$ and $x$ are not involved in any other terms in any of these potentials.

If the interior marked point is as in case (d), then $W_{\T}$ and $W_{\T}'$ coincide.
\end{pf}
\end{lem}

We recall the following special case of~\cite[Thm.\ 8.1]{LF4}.

\begin{thm} \cite{LF4} \label{t:flipmutate}
Let $\T,\T'$ be triangulations of $(S,M)$.
If $\T'$ is obtained from $\T$ by flipping at an arc $\alpha_k$ then $\mu_k(Q_{\T},W'_{\T})$ is right equivalent to $(Q_{\T'},W'_{\T'})$.
\end{thm}

By~\cite[Thm.\ 7.1]{dwz1}, it
follows from this that the quiver of $\mu_k(Q_{\T},W_{\T})$ coincides with the
quiver obtained from $Q_{\T}$ by Fomin-Zelevinsky quiver mutation at $k$.

Hence we can effectively ignore potentials:

\begin{prop}\label{prop:pot-is-sum-cyc}
Any mutation-Dynkin quiver with potential $(\widetilde{Q},\widetilde{W})$
of type $A$ or $D$ is right equivalent to
$(\widetilde{Q},W_{\widetilde{Q}})$,
where $W_{\widetilde{Q}}$ is the sum of all chordless cycles in
$\widetilde{Q}$.

\begin{pf}
Note that a Dynkin quiver with zero potential is
of the form $(Q_{\T},W_{\T})$ for some triangulation $\T$ (see~\cite{fst08}).
Suppose that $(\widetilde{Q},\widetilde{W})$ is obtained from a Dynkin quiver with zero potential  by iterated mutation in the sense
of~\cite{dwz1}.
Then, by Theorem~\ref{t:flipmutate} and Lemma~\ref{l:rightequivalent}, $(\widetilde{Q},\widetilde{W})$ is right equivalent to $(Q_{\T},W_{\T})$ for some triangulation $\T$ of $(S,M)$.
\end{pf}
\end{prop}

Note that an alternative proof of Proposition~\ref{prop:pot-is-sum-cyc}
would be to compute the mutation of a quiver with potential
$(Q_{\T},W'_{\T})$ directly, and show that it is right equivalent to
$(Q_{\T'},W'_{\T'})$. This is not too difficult to do, but requires consideration
of several cases and still requires arguments dealing with changes
of sign as in Lemma~\ref{l:rightequivalent}, so we instead refer to~\cite{LF4} above.

\subsection{Differential graded algebras and modules}
Let $\F$ be an algebraically closed field.  We think of $\F$ as a graded $\F$-algebra concentrated in degree $0$.  
If $V=\bigoplus V_i$ is a graded $\F$-module then let $V[j]$ be the graded $\F$-module with $(V[j])_i=V_{i+j}$. 
If $f:V\to W$ is a map of graded vector spaces with homogeneous components $f_i:V_i\to W_i$ then let $f[j]:V[j]\to W[j]$ be the map of graded vector spaces with homogeneous components $f[j]_i:V[j]_i\to W[j]_i$ defined by $f[j]_i(v)=(-1)^jf_{i+j}(v)$ for $v\in V[j]_i=V_{i+j}$.
Thus $[1]$ is an endofunctor of the category of graded $\F$-modules, called the \emph{shift functor}.

We say that a map $f:V\to W$ of graded vector spaces has degree $i$ to mean that $f$ is a map $V\to W[i]$.
We use the Koszul sign rule for graded $\F$-algebras, so if $f:V\to V'$ and $g:W\to W'$ are maps of graded $\F$-modules of degree $m$ and $n$ then
$$(f\otimes g)(v\otimes w)=(-1)^{in}f(v)\otimes g(w)$$
for $v\in V_i$ and $w\in W$.

A unital differential graded algebra (or dg-algebra, or dga) over $\F$ is a graded $\F$-algebra $A=\bigoplus_{i\in Z}A_i$ with multiplication $m:A\otimes_\F A\to A$ of degree $0$ together with a unit $\iota:\F\into A$ and an $\F$-linear differential $d:A\to A$ of degree $+1$.  These should satisfy the following relations:
\begin{itemize}
 \item the associativity relation $m\circ(1\otimes m)=m\circ(m\otimes 1)$;
 \item the boundary relation $d^2=0$;
 \item the Leibniz relation $d\circ m=m\circ (1\otimes d+d\otimes 1)$;
 \item the unital relation $m\circ (\id_A\otimes \iota)=m\circ (\iota\otimes \id_A)$, which should agree with the $\F$-algebra structure of $A$.
\end{itemize}
We often denote our dga by $(A,d)$, or simply by $A$.  Each dga $(A,d)$ has an underlying unital graded algebra, obtained by simply forgetting the differential, which we denote $u(A)$.

A left module $M$ for $A$ is a graded left $\F$-module $M$ with a left action $m_M:A\otimes M\to M$ of 
$u(A)$ 
together with a 
map $d_M:M\to M$ of degree $+1$, called a \emph{differential},
such that
$$d_M\circ m_M=m_M\circ (1\otimes d_M+d\otimes 1).$$
We always have the \emph{regular} module $M=A$ with $d_M=d$ and $m_M=m$.
Similarly, a right module $M$ for $A$ is a graded right $\F$-module $M$ with a right action $m_M:M\otimes A\to M$ of 
$u(A)$ 
together with a differential $d_M$ such that $d_M\circ m_M=m_M\circ (1\otimes d+d_M\otimes 1)$.  If $(M,d_M)$ is an $A$-module, then $(M[1],d_M[1])$ is also an $A$-module, which we sometimes just write as $M[1]$.  Modules for $A$ are modules for $u(A)$, simply by forgetting the differential.

A map $f:M\to N$ of left $A$-modules is a degree $0$ map of 
$u(A)$-modules
such that $f$ commutes with the differentials: $d_N\circ f=f\circ d_M$.  We thus obtain a category $A\Mod$ of left $A$-modules, and we write the morphism spaces in this category as $\Hom_{A\Mod}(M,N)$. 
$A\Mod$ is an $\F$-category: each morphism space is an $\F$-module. 

Given two differential algebras $(A,d_A)$ and $(B,d_B)$, an $A\da B$-bimodule $(M,d_M)$ is a graded $\F$-module which is a left $(A,d_A)$-module with left action $m^\ell$ and a right $(B,d_B)$-module with right action $m^r$ where the two actions commute: $m^r\circ(m^\ell\otimes \id_B)=m^\ell\circ(\id_A\otimes m^r)$.  We will always assume that $\F$ acts centrally. 
Under this assumption we can, and will, identify left $A$-modules with $A\da k$-bimodules and $A\da B$-bimodules with left $A\otimes_\F B^\op$-modules, where $B^\op$ denotes the algebra $B$ with the order of multiplication reversed.
A map of bimodules should commute with the differential on both the left and the right, and we obtain an $\F$-category $A\biMod B$ of $A\da B$-bimodules.

Given a map $f:M\to N$ of left $A$-modules, we can construct a new left $A$-module called the \emph{cone} of $f$, denoted $\cone(f)$.  As a left module for $u(A)$, we have $\cone(f)=N\oplus M[1]$.  The differential is given by:
$$\begin{pmatrix}
d_N &0\\
f[1] & d_{M[1]}
\end{pmatrix}.$$
If $L$ is isomorphic to $\cone(f)$ for some map $f:M\to N$, we say that $L$ is an \emph{extension} of $M$ by $L[-1]$.

We will use the following lemma, whose proof follows immediately from the definitions, repeatedly.
\begin{lem}\label{lem:coneFf-comm}
Let $f:M\to N$ be a map in $A\Mod$.
\begin{enumerate}
 \item 
Let $F:A\Mod\to B\Mod$ be an additive functor which commutes with the shift functor.  Then we have an isomorphism $\cone(Ff)\cong F\cone(f)$ in $B\Mod$.
 \item 
For any commutative diagram
$$\xymatrix{
M\ar[r]^f\ar[d]^{\varphi_M}_{\sim} & N\ar[d]^{\varphi_N}_{\sim}\\
M'\ar[r]^{f'} & N'
}$$
in $A\Mod$ where both $\varphi_M$ and $\varphi_N$ are isomorphisms, we have an isomorphism $\varphi_N\oplus\varphi_{M}[1]:\cone(f)\to\cone(f')$ of $A$-modules.
\end{enumerate}
\end{lem}

Let $(A,d_A)$, $(B,d_B)$, and $(C,d_C)$ be dgas.  If $(M,d_M)$ is an $A\da B$-bimodule and $(N,d_N)$ is an $A\da C$-bimodule then let $\Hom^i_{A}(M,N)$ be the space of all graded left $u(A)$-module maps $f:M\to N$ of degree $i$.  Note that we do not require that these maps commute with the differential.  We define $\Hom_{A}(M,N)=\bigoplus_{i\in \Z}\Hom_{A}^i(M,N)$, and this is a graded $u(B)\da u(C)$-bimodule.  We also have a version for right modules, which we write as $\Hom_{A^\op}(M,N)$.

Note the distinction between $\Hom_{A}(M,N)$ and the hom spaces in the category $A\Mod$. With the differential $d(f)=d_N\circ f-(-1)^i f\circ d_M$ for $f\in \Hom^i_{ A}(M,N)$, $\Hom_{A}(M,N)$ becomes a $B\da C$-bimodule.  Similarly, if $(M,d_M)$ is an $B\da A$-bimodule and $(N,d_N)$ is a $C\da A$-bimodule, $\Hom_{A^\op}(M,N)$ is a $C\da B$-bimodule. 
$\Hom_{ A}(-,-)$ is the internal hom in the bimodule category, and we can recover the hom spaces in $A\Mod$ as the $0$-cycles of $\Hom_{ A}(M,N)$.

If $(M,d_M)$ is an $A\da B$-bimodule and $(N,d_N)$ is a $B\da C$-bimodule then let $M\otimes_B N$ denote the space $M\otimes_{u(B)}N$.  It is a graded $u(A)\da u(C)$-bimodule: if $m\in M_i$ and $n\in N_j$ then $m\otimes n$ has degree $i+j$.  With the differential $d_M\otimes \id_N +\id_M\otimes d_N$, it becomes an $A\da C$-bimodule.

For an $A\da B$-bimodule $(M,d_M)$, we thus have functors $M\otimes_B-:B\Mod\to A\Mod$ and $\Hom_{ A}(M,-):A\Mod\to B\Mod$.
The functor $M\otimes_B-$ is left adjoint to $\Hom_{ A}(M,-)$.  For $(N,d_N)$ a left $A$-module, the counit $\ev_N:M\otimes_B\Hom_{ A}(M,N)\to N$ of the adjunction is the evaluation map, which acts as $x\otimes f\mapsto (-1)^{ij}f(m)$ for $x\in M_i$ and $f\in \Hom_{A}^j(M,N)$.

\subsection{Derived categories}

Our references are \cite{kel-derdg,kel-ondg}.

If $A$ is a graded vector space and $d$ is a differential, i.e., a degree $+1$ endomorphism of $A$ which satisfies $d^2=0$, then the $i$th homology of $A$, denoted $H_i(A)$, is the subquotient $\ker d_i/\im d_{i-1}$, where $d_i:A_i\to A_{i+1}$ denotes the restriction of $d$ to $A_i$.  If $(A,d)$ is a dga then the homology $H(A)=\bigoplus H_i(A)$ is a graded algebra, and if $M$ is a left $A$-module then $H(M)=\bigoplus H_i(M)$ is a left $H(A)$-module.  In fact, taking homology is a functor from the category of $A$-modules to the category of graded $H(A)$-modules.  We say that a left $A$-module $M$ is \emph{acyclic} if $H(M)=0$, and that a map $f:M\to N$ of $A$-modules is a \emph{quasi-isomorphism} if $H(f)$ is an isomorphism.

The \emph{category up to homotopy} of $A\Mod$, denoted $\K( A)$, is the $\F$-category whose objects are all left $A$-modules and whose morphism spaces, for $M,N\in A\Mod$, are $\Hom_{K(A)}(M,N)=H_0\Hom_{A}(M,N)$.  The \emph{derived category} of $A$, denoted $\D(A)$, is the $\F$-category obtained by localizing $\K( A)$ at the full subcategory of acyclic $A$-modules.  As a map of modules is a quasi-isomorphism if and only if its cone is acyclic, this is equivalent to localizing $\K( A)$ at the class of all quasi-isomorphisms.  So we have a canonical functor $\K( A)\to \D( A)$, which we call the projection functor.
The finite-dimensional derived category, denoted $\Dfd( A)$, is the full subcategory of $\D( A)$ on objects with finite-dimensional total homology, i.e., on $A$-modules $M$ such that $H(M)$ is a finite-dimensional $\F$-vector space.

Let $(A,d_A)$ be a dga.  We say that:
\begin{itemize}
\item $P\in A\Mod$ is \emph{indecomposable projective} if it is an indecomposable direct summand of the regular module,
\item $P\in A\Mod$ is \emph{relatively projective} if it is a direct sum of shifts of indecomposable projective modules, and
\item
$P\in A\Mod$ is \emph{cofibrant} if, for each surjective quasi-isomorphism $f:M\to N$, the map $\Hom_{A\Mod}(P,f):\Hom_{A\Mod}(P,M)\to \Hom_{A\Mod}(P,N)$ is surjective.  
\end{itemize}
The following result (see \cite[Section 3]{kel-derdg} and \cite[Proposition 2.13]{ky}) characterizes cofibrant modules.
\begin{prop}[Keller]\label{prop:char-cofib}
An $A$-module $P$ is cofibrant if and only if it is an iterated extension of a relatively projective module by other relatively projective modules, possibly infinitely many times.
\end{prop}

Let $A\cofib$ denote the full subcategory of $\K(A)$ on the cofibrant objects.
The projection functor $\K(A)\to \D(A)$ induces an equivalence $A\cofib\arr{\sim}\D(A)$.  Each $A$-module $M$ has a \emph{cofibrant replacement}, defined up to quasi-isomorphism and denoted $\p M$, which can be realized as the image of $M$ under the left adjoint $\D(A)\to \K(A)$ to the canonical projection functor \cite[Proposition 3.1]{kel-ondg}.

Let $(B,d_B)$ be another dga and let $F:A\Mod\to B\Mod$ be an additive functor.  Then $F$ preserves chain homotopies, so induces a functor $\K(F):\K(A)\to \K(B)$.  If $\K(F)$ preserves quasi-isomorphisms then, by the universal property of localization, it induces a functor $\D(F):\D(A)\to \D(B)$.  If $P\in A\biMod B$ is cofibrant as a left $A$-module then, by \cite[Theorem 3.1(a)]{kel-derdg} and \cite[Proposition 2.13]{ky}, $\Hom_A(P,-)$ preserves acyclic modules, and so preserves quasi-isomorphisms.  By imitating the proof of \cite[Theorem 3.1(a)]{kel-derdg} we see that if $P\in A\biMod B$ is cofibrant as a right $B$-module then $P\otimes _B-$ also preserves acyclic modules.  We often write $P\otimes _B-$ and $\Hom_A(P,-)$, instead of $\D(P\otimes _B-)$ and $\D(\Hom_A(P,-))$, for the induced functors $\D(B)\to \D(A)$ and $\D(A)\to \D(B)$.

For an arbitrary $M\in A\biMod B$, we can obtain a functor $M\dert_B-:\D(B)\to \D(A)$, known as the \emph{left derived functor} of $M\otimes_B-$, by composing the cofibrant replacement functor $\D(B)\to \K(B)$, the tensor functor $\K(M\otimes_B-):\K(B)\to \K(A)$, and the projection functor $\K(A)\to \D(B)$.  By \cite[Lemma 6.3(a)]{kel-derdg}, we have an isomorphism $M\dert_BN\cong \p M\otimes_B N$ for all $N\in \D(B)$.
The following basic, but useful, lemma says that this isomorphism is natural.
\begin{lem}\label{lem:derived-f}
Let $M\in\rMod B$.
\begin{enumerate}
 \item\label{lem:pt:derlr} We have a natural isomorphism of functors $\p M\otimes_B -\cong M\dert_B-$.
 \item If $M$ is cofibrant then we have a natural isomorphism of functors $M\otimes_B-\cong M\dert_B-$.
\end{enumerate}

\begin{pf}
\begin{enumerate}
 \item We need to show that for each $N\in B\Mod$ there is a quasi-isomorphism $\varphi_N:\p M\otimes_BN\to M\otimes_B\p N$ such that, for all maps $f:N\to N'$, the diagram
$$\xymatrix{
\p M\otimes_BN\ar[r]^{\varphi_N}\ar[d]^{\p M\otimes f} &M\otimes_B\p N\ar[d]^{M\otimes\p f}\\
\p M\otimes_BN'\ar[r]^{\varphi_{N'}} &M\otimes_B\p N'\\
}$$
commutes.  Consider the following diagram:
$$\xymatrix @C=10pt {
 &\p M\otimes_B \p N\ar@/_/[ddl]_{\p M\otimes \pi_N}\ar@/^/[ddr]^{\pi_M\otimes \p N}  &&&&&\\
 &&&&& \p M \otimes_B \p N'\ar@/_/[ddl]_{\p M\otimes \pi_{N'}}\ar@/^/[ddr]^{\pi_M\otimes \p N'} &\\
\p M\otimes_BN\ar@{-->}[rr]^{\varphi_N}\ar[drrrr]^{\p M\otimes f} &&M\otimes_B\p N\ar[drrrr]^{M\otimes\p f} &&&&\\
&&&&\p M\otimes_BN'\ar@{-->}[rr]^{\varphi_{N'}} &&M\otimes_B\p N'\\
}$$
As $\p M$ and $\p N$ are cofibrant and $\pi_M$ and $\pi_N$ are quasi-isomorphisms, both $\p M\otimes \pi_N$ and $\pi_M\otimes \p N$ are quasi-isomorphisms, so we can define $\varphi_N=(\pi_M\otimes \p N)\circ (\p M\otimes \pi_N)^{-1}$ and it is a quasi-isomorphism.  Then to check naturality we need to show that
$$(M\otimes \p f)\circ (\pi_M\otimes \p N)\circ (\p M\otimes \pi_N)^{-1} = (\pi_M\otimes \p N')\circ (\p M\otimes \pi_{N'})^{-1}\circ (\p M\otimes f).$$
By the bifunctoriality of the tensor product, the left hand side is equal to
$$(\pi_M\otimes \p N')\circ (\p M\otimes \p f) \circ (\p M\otimes \pi_N)^{-1}$$
so we just need to show that 
$$\p f \circ \pi_N^{-1}= \pi_{N'}^{-1}\circ  f$$
but this follows from the functoriality $ f \circ \pi_N= \pi_{N'}\circ \p f$ of the cofibrant replacement functor $\p $.

\item We just need to show that, for $M$ cofibrant, there is a natural isomorphism $M\otimes_B-\cong \p M\otimes_B-$, and then the result will follow by part \ref{lem:pt:derlr} of the lemma.  This follows because $\pi_M:\p M\to M$ is a quasi-isomorphism and by the bifunctoriality of the tensor product.
\end{enumerate}
\end{pf}
\end{lem}

If the functor $M\dert_B-$ is an equivalence $\D(B)\arr{\sim} \D(A)$, we say that $M$ is a \emph{tilting module}.

We say that a module $M$ is of \emph{finite projective dimension} if its cofibrant replacement is an iterated extension of finitely many shifted indecomposable projective modules.

The following basic lemma will also be useful.  It can be found as \cite[Lemma 6.2(a)]{kel-derdg}.  We include a proof for the convenience of the reader.

\begin{lem}\label{lem:homastens}
Let $M$ be a left $A$-module of finite projective dimension and let $P$ be its cofibrant replacement. 
Then we have a natural isomorphism of functors
$$\Hom_A(P,A)\otimes_A-\arr\sim\Hom_A(P,-):A\Mod\to \F\Mod.$$

\begin{pf}
First note that, for any $P\in A\Mod$, we always have a natural transformation
$$\Hom_A(P,A)\otimes_A-\to\Hom_A(P,-)$$
obtained by staring with the map
$$\ev\otimes1:\left(P\otimes_\F\Hom_A(P,A)\right)\otimes_AM\to A\otimes_AM,$$
using the associativity isomorphism to obtain a map
$$P\otimes_\F\left(\Hom_A(P,A)\otimes_AM\right)\to A\otimes_AM,$$
using the adjunction
$$\Hom_\F(\Hom_A(P,A)\otimes_AM, \Hom_A(P,A\otimes_AM))
\cong\Hom_A(P\otimes_\F\left(\Hom_A(P,A)\otimes_AM\right), A\otimes_AM),$$
and finally using the natural isomorphism $A\otimes_AM\cong M$.

To show that our natural transformation is an isomorphism, we use induction on the number of times we need to extend a summand of the regular module to obtain $P$.  We handle the base case as follows: the natural transformation is certainly an isomorphism when $P$ is the regular module and so, as hom functors commute with finite direct sums, it is an isomorphism for all summands of the regular module.  For our inductive step, suppose the lemma holds for $P_1$ and $P_2$, and let $P=\cone(f)$ for some map $f:P_1\to P_2$.  Then, for $M\in A\Mod$,
one can check that the map $\Hom_A(P,A)\otimes_AM\to\Hom_A(P,M)$ 
comes from the commutative diagram
$$\xymatrix @C=50pt {
\Hom_A(P_1,A)\otimes_AM\ar[d] &\Hom_A(P_2,A)\otimes_AM\ar[l]_{\Hom(f,A)\otimes M}\ar[d]\\
\Hom_A(P_1,M) &\Hom_A(P_2,M)\ar[l]_{\Hom(f,M)}
}$$
as in the construction from the second half of Lemma \ref{lem:coneFf-comm}, where the vertical maps come from the natural transformation described above.  Therefore, as both vertical maps are isomorphisms by induction, $\Hom_A(P,A)\otimes_AM\to\Hom_A(P,M)$ is an isomorphism.
\end{pf}
\end{lem}

\subsection{Spherical twists}\label{ss:spher-tw}
Our references are \cite{st,rz,gra1}.

Let $(A,d)$ be a dga and $M$ be a left $A$-module with finite dimensional total homology.  Let $d\in\Z$.  Following \cite{st}, we say that $M$ is $d$-spherical if:
\begin{itemize}
 \item $M$ is a $d$-Calabi-Yau object, i.e., we have an isomorphism $$\Hom_{\Dfd(A)}(M,N)\cong \Hom_{\Dfd(A)}(N,M[d])^*$$ which is functorial in $N$, and
 \item $\bigoplus_{i\in\Z}\Hom_{\Dfd(A)}(M,M[i])$ is isomorphic as a graded algebra to $\F[x]/\gen{x^2}$, with $x$ in degree $d$.
\end{itemize}
Associated to any spherical object $M$, we have a spherical twist functor $F_M:\Dfd(A)\to\Dfd(A)$ which is defined as follows.  First, let $P=\p M$ be a cofibrant replacement of $M$.  Then let $X_M$ be the cone of the following map of $A\da A$-bimodules:
$$P\otimes_\F\Hom_A(P,A)\arr{\ev}A$$
where the nonzero map is the obvious evaluation map.  As both $\Hom_A(P,A)$ and $A$ are cofibrant, $X_M$ is cofibrant as a right $A$-module.  Then we define the spherical twist at $M$ by
$$F_M=X_M\otimes_A-:\Dfd(A)\to\Dfd(A).$$
The spherical twist is an autoequivalence of $\Dfd(A)$ (so $X_M$ is a tilting module).

Note that, by Lemmas \ref{lem:coneFf-comm} and \ref{lem:homastens}, if $M$ has finite projective dimension then 
$$F_M(N)\cong P\otimes_\F\Hom_A(P,N)\arr{\ev}N.$$

We will need a simple result on the commutation relation of spherical twists with derived equivalences.  It is a generalization of \cite[Lemma 2.11]{st}.
\begin{prop}\label{prop:conj-sph}
Let $A,B$ be dgas.  Let $T\in B\biMod A$ be a tilting module and $\Phi=T\dert_A-:\Dfd(A)\to\Dfd(B)$ be the associated derived equivalence.  Let $M\in A\Mod$ have finite dimensional total homology and suppose it is $d$-spherical, for some $d\in\Z$.  Suppose that $\Phi(M)\in B\Mod$ has finite dimensional total homology.
Then $\Phi(M)$ is also $d$-spherical and we have an isomorphism of functors
$$\Phi\circ F_M\cong F_{\Phi(M)}\circ \Phi:\Dfd(A)\arr{\sim} \Dfd(B).$$
In particular, we have an isomorphism
$$F_{\Phi(M)}\cong \Phi\circ F_M\circ \Phi^{-1}:\Dfd(B)\arr{\sim} \Dfd(B)$$
where $\Phi^{-1}$ 
 is the quasi-inverse functor of $\Phi$.

\begin{pf}
As $\Phi:\Dfd(A)\to\Dfd(B)$ is a derived equivalence it has quasi-inverse $\Phi^{-1}:\Dfd(B)\to\Dfd(A)$, and so we have isomorphisms
$$\Hom_{\Dfd(B)}(\Phi(M),\Phi(M)[i])\cong \Hom_{\Dfd(A)}(M,M[i])$$
and
\begin{align*}
 \Hom_{\Dfd(B)}(\Phi(M),N) &\cong \Hom_{\Dfd(A)}(M,\Phi^{-1}(N))\\
   &\cong \Hom_{\Dfd(A)}(\Phi^{-1}(N),M[d])^*\\
   &\cong \Hom_{\Dfd(B)}(N,\Phi(M)[d])^*,
\end{align*}
the second natural in $N\in \Dfd(B)$, using the facts that $M$ is a $d$-Calabi-Yau object and the shift functor $[d]$ commutes with all triangulated functors.  Thus $\Phi(M)$ is $d$-spherical.

By Lemma \ref{lem:derived-f} we may assume that $T$ is cofibrant as a right $B$-module and that $\Phi=T\otimes_A-$.
We want to show that $T\otimes_A X_M\otimes_A-\cong X_{\Phi(M)}\otimes_B T\otimes_A-$, so it is enough to check that we have an isomorphism $T\otimes_A X_M\cong X_{\Phi(M)}\otimes_B T$ in $\Dfd(B\otimes_\F A^\op)$.  
To construct this isomorphism, we use the following 
extension of Lemma \ref{lem:coneFf-comm}, which follows from the triangulated 5-lemma:
for any commutative diagram
$$\xymatrix{
M\ar[r]^f\ar[d]^{\varphi_M}_{\sim} & N\ar[d]^{\varphi_N}_{\sim}\\
M'\ar[r]^{f'} & N'
}$$
in $B\biMod A$ with $\varphi_M$ and $\varphi_N$ both quasi-isomorphisms, 
we have a quasi-isomorphism 
$\varphi_N\oplus\varphi_{M[1]}:\cone(f)\to\cone(f')$.

As above, write $P=\p M$.  Then, by the first part of Lemma \ref{lem:coneFf-comm}, we just need to find two vertical maps which are quasi-isomorphisms and make the following diagram commute:
$$\xymatrix @C=56pt {
T\otimes_AP\otimes_\F \Hom_B(T\otimes_AP,B)\otimes_B T\ar[r]^(.68){\ev\otimes 1_T}\ar@{-->}[d]_{\sim} &B\otimes_B T\ar@{-->}[d]_{\sim}\\
T\otimes_AP\otimes_\F \Hom_A(P,A)\ar[r]^(.68){1_T\otimes\ev} &T\otimes_AA
}$$
Our plan is to do this in stages: we will show that the vertical maps in the following diagram exist, and are quasi-isomorphisms, and that the diagram commutes:
$$\xymatrix @C=56pt {
T\otimes_AP\otimes_\F \Hom_B(T\otimes_AP,B)\otimes_B T\ar[r]^(.68){\ev\otimes 1_T}\ar[d]_{\sim} &B\otimes_B T\ar[d]_{=}\\
T\otimes_AP\otimes_\F \Hom_B(T\otimes_AP,B\otimes_B T)\ar[r]^(.68){\ev}\ar[d]_{\sim} &B\otimes_B T\ar[d]_{\sim}\\
T\otimes_AP\otimes_\F \Hom_B(T\otimes_AP,T\otimes_AA)\ar[r]^(.68){\ev}\ar[d]_{\sim} & T\otimes_AA\ar[d]_{=}\\
T\otimes_AP\otimes_\F \Hom_A(P,A)\ar[r]^(.68){1_T\otimes\ev} &T\otimes_AA
}$$

Let us show that the first square commutes.  We'll introduce some temporary notation for the rest of this proof.  Let $F$ and $G$ denote the functors $F=T\otimes_A P\otimes_\F-$ and $G=\Hom_B(T\otimes_A P,-)$, so $F$ is left adjoint to $G$, and let $H$ denote the functor $-\otimes_BT$.  Then we have unit and counit natural transformations $\varepsilon:FG\to 1$ and $\eta:1\to GF$, and a natural isomorphism $\zeta:HF\arr\sim FH$ coming from the associativity isomorphism for tensor products.  We first need to define a map 
$$HFGB=T\otimes_AP\otimes_F\Hom_B(T\otimes_AP,B)\otimes_B T\to T\otimes_AP\otimes_F\Hom_B(T\otimes_AP,B\otimes_BT)=FGHB$$
We define this as the composite
$$HFGB\arr{\zeta GB}  FHGB\arr{F\eta_{HGB}} FGFHGB\arr{FG\zeta^{-1}GB} FGHFGB\arr{FGH\varepsilon_B} FGHB.$$
One checks that this is an isomorphism using the same argument as in Lemma \ref{lem:homastens}.
To see that the diagram commutes, we break it up into smaller diagrams as follows:
$$\xymatrix{
HFGB\ar[d]^\sim_{\zeta GB}\ar[rrr]^{H\varepsilon_B} &&& HB\\
FHGB\ar[r]^1\ar[d]_{F\eta_{HGB}} &FHGB\ar[dr]^\sim_{\zeta^{-1} GB} &&\\
FGFHGB\ar[d]^\sim_{FG\zeta^{-1} GB}\ar[ru]_{\varepsilon_{FHGB}} && HFGB\ar[uur]^{H\varepsilon_B} &\\
FGHFGB\ar[rrr]^{FGH\varepsilon_B}\ar[urr]^{\varepsilon_{HFGB}} &&& FGHB\ar[uuu]^{\varepsilon_{HB}}
}$$
Now we see that both squares commute by the naturality of $\varepsilon$, the triangle commutes by the triangle identity $\varepsilon_F\circ F\eta=1_F$, and the pentagon commutes because the isomorphisms are defined by $\zeta$ and its inverse.

To define the second square we use the obvious composite isomorphism $B\otimes_BT\arr\sim T\arr\sim T\otimes_AA$.  This commutes because the evaluation map is a counit, and therefore a natural transformation.

To show that the third square commutes, we introduce some more notation.  Let $F'$ and $G'$ denote the functors $F'=P\otimes_\F-$ and $G'=\Hom_A(P,-)$, so $F'$ is left adjoint to $G'$, and let $H'$ and $I'$ denote the functors $H'=T\otimes_A-$ and $I'=\Hom_B(T,-)$, so $H'$ is left adjoint (in fact, quasi-inverse) to $I'$.  We denote the counit and unit maps of the first adjunction by $\varepsilon':F'G'\to 1$ and $\eta':1\to G'F'$, and of the second adjunction by $\varepsilon'':H'I'\to1$ and $\eta'':1\to I'H'$.  Note that, because $H'$ induces an equivalence of derived categories, $\varepsilon''$ and $\eta''$ give quasi-isomorphisms when applied to any object.

Using the associativity isomorphism for tensor products we have a natural isomorphism of functors $F\cong H'F'$, and by the uniqueness of right adjoints (or by using the tensor-hom adjunctions directly) this gives us another natural isomorphism $G\cong G'I'$.

We now redraw our final square, breaking it up into smaller diagrams:
$$\xymatrix{
FGH'A\ar[rrr]^{\varepsilon'_{H'A}}\ar[d]^\sim &&&H'A\\
H'F'G'I'H'A\ar[rr]^{H'\varepsilon'_{I'H'A}} &&H'I'H'A\ar[ur]^{\varepsilon''H'A} &\\
H'F'G'A\ar[rrr]^{H'\varepsilon'_A}\ar[u]^{H'F'G'\eta''_A} &&&H'A\ar[ul]^{H'\eta''_A}\ar[uu]^1
}$$
Here, the top square commutes by definition of the isomorphisms $F\cong H'F'$ and $G\cong G'I'$, the triangle commutes by the triangle identity $\varepsilon''H'\circ H'\eta''=1_{H'}$, and the bottom square commutes by the naturality of $\varepsilon'$.
\end{pf}
\end{prop}

We now describe the braid relations for spherical twists, as in Propositions 2.12 and 2.13 of \cite{st} (see also \cite{rz,gra2}).
\begin{prop}\label{prop:braid-relns}
 Suppose that $M$ and $N$ are spherical objects of $\Dfd(A)$ and let
$$(M,N)=\dim_\F\bigoplus_{n\in\Z}\Hom_{\Dfd(A)}(M,N[n]).$$
 Let $F_M,F_N:\Dfd(A)\arr\sim\Dfd(A)$ be the associated spherical twists.
\begin{itemize}
 \item If $(M,N)=0$ then $F_M\circ F_N\cong F_N\circ F_M$;
 \item if $(M,N)=1$ then $F_M\circ F_N\circ F_M\cong F_N\circ F_M\circ F_N$.
\end{itemize}
\end{prop}

\subsection{Ginzburg dg-algebras}
There is a by now well-known method to associate a differential graded algebra to a quiver with potential (\cite[Section 5]{g-cy}, \cite[Section 2.6]{ky}).

Let $(Q,W)$ be a quiver with potential.  Construct a new quiver $\overline{Q}$ by adding arrows to $Q$: for each arrow $a:i\to j$ in $Q$ we add a new arrow $a^*:j\to i$, and for each vertex $i$ in $Q$ we add a new arrow $t_i:i\to i$.  We view $\overline{Q}$ as a graded quiver with the arrows of $Q$ in degree $0$, the arrows $a^*$ in degree $-1$, and the arrows $t_i$ in degree $-2$.  This induces a grading on the path algebra $\F\overline{Q}$ of $\overline{Q}$ such that the degree $0$ part $(\F\overline{Q})_0$ is just the path algebra of $\F Q$ of $Q$.  Let $\overline{J}$ denote the ideal of $\F\overline{Q}$ generated by the arrows of $\overline{Q}$, and let $\widehat{\F\overline{Q}}$ denote the completion of the graded algebra $\F\overline{Q}$ with respect to $\overline{J}$, as in Section \ref{ss:qps}.

We define a differential $d$ on $\widehat{\F\overline{Q}}$ by imposing that $d$ is zero on 
each idempotent $e_i$ associated to a vertex $i$ of $Q$,
specifying how $d$ acts on arrows of $\overline{Q}$,
and then extending to the rest of $\widehat{\F\overline{Q}}$ using the Leibniz rule and continuity.  For degree reasons, we must have $d(a)=0$ for each arrow $a$ of $Q$.  For arrows $a^*$, we set $d(a^*)=\partial_aW$, where $\partial_a$ denotes the cyclic derivative.  For arrows $t_i$, we set $d(t_i)=e_i\left(\sum aa^*-a^*a\right)e_i$, where we sum over all arrows $a$ of $Q$.  Then $\Gamma_{Q,W}=(\widehat{\F\overline{Q}},d)$ is called the \emph{Ginzburg dga} of $(Q,W)$.  

Note that \cite[Lemma 2.9]{ky} if $(Q_1,W_1)$ and $(Q_2,W_2)$ are right equivalent, then we have an isomorphism of dgas $\Gamma_{Q_1,W_1}\cong\Gamma_{Q_2,W_2}$.  Hence, if we are working with 
quivers with potential of mutation type $A$ or $D$, by Proposition \ref{prop:pot-is-sum-cyc} we only need to consider the Ginzburg dgas $\Gamma_{Q,W_Q}$, and so can denote them $\Gamma_Q$.

Keller and Yang showed \cite[Theorem 3.2]{ky} that QP-mutation lifts to equivalences of derived categories of Ginzburg dgas.  Suppose that $(Q,W)$ is a QP and that $(Q',W')=\mu_k(Q,W)$ for some $k\in I$.
\begin{thm}[Keller-Yang]
There is a tilting complex $T$ which gives an equivalence of triangulated categories
$$\mu_k=\Hom_{\Gamma_{Q',W'}}(T,-):\D(\Gamma_{Q',W'})\to\D(\Gamma_{Q,W})$$
and it restricts to an equivalence of triangulated categories
$$\mu_k=\Hom_{\Gamma_{Q',W'}}(T,-):\Dfd(\Gamma_{Q',W'})\to\Dfd(\Gamma_{Q,W}).$$
\end{thm}

Recall that, for a dga $A$, the finite-dimensional derived category $\Dfd(A)$ is $d$-Calabi-Yau if there exists a bifunctorial isomorphism
$$\Hom_{\Dfd(A)}(M,N)\cong \Hom_{\Dfd(A)}(N,M[d])^*$$
where $(-)^*$ denotes the $k$-linear dual.
We will need the following important result \cite[Theorem 6.3 and Theorem A.12]{ke-dcyc}:
\begin{thm}[Keller, Van den Bergh]\label{thm:ginz-cy}
The category $\Dfd(\Gamma_{Q,W})$ is $3$-Calabi-Yau.
\end{thm}

Let $(Q,W)$ be a QP and $\Gamma=\Gamma_{Q,W}$.
Associated to each vertex $i$ of $Q$, we have the one-dimensional simple left $\Gamma$-module, which we denote $S_i$.
In \cite[Section 2.14]{ky}, Keller and Yang explain how to construct the cofibrant replacement of $S_i$: as long as we remember the differential, we can proceed as if the Ginzburg dga were an ordinary hereditary algebra, and the underlying $u(\Gamma)$-module of $\p S_i$ is the direct sum of one copy of the projective $P_i$ and one copy of the shifted projective $P_j[1]$ for each arrow $j\to i$ in $Q$.  Using this, they show \cite[Lemma 2.15]{ky}:
\begin{lem}\label{lem:dimext}
Let $i,j\in I$ and $n\in Z$ and $\Gamma=\Gamma_{Q,W}$.  Then $\Hom_{\Dfd(\Gamma)}\left(S_i,S_j[n]\right)=0$ if $n\neq0,1,2,3$, and
$$\dim_\F\Hom_{\Dfd(\Gamma)}\left(S_i,S_j[n]\right)=\begin{cases}
\delta_{ij} &\text{ if }n=0;\\
\#\{\text{arrows }i\to j\text{ in }Q\} &\text{ if }n=1;\\
\#\{\text{arrows }j\to i\text{ in }Q\} &\text{ if }n=2;\\
\delta_{ij} &\text{ if }n=3.,
\end{cases}$$
where $\delta_{ij}$ is the Kronecker delta.
\end{lem}

\subsection{Relations between functors}\label{ss:reln-per}

By Theorem \ref{thm:ginz-cy}, every object of $\Dfd(\Gamma_{Q,W})$ is a $3$-Calabi-Yau object.  
By Lemma \ref{lem:dimext}, $$\bigoplus_{j\in\Z}\Hom_{\Dfd(\Gamma_{Q,W})}(S_i,S_i[j])\cong\F[x]/\gen{x^2}$$ with $x$ in degree $3$. 
Hence $S_i$ is $3$-spherical, and we have a spherical twist $F_{S_i}$ associated to $S_i$.  We will sometimes write $F_i$ instead of $F_{S_i}$.

Let $k$ be a vertex of $Q$, and write $(Q',W')=\mu_k(Q,W)$.  Then write $\Gamma=\Gamma_{Q,W}$ and $\Gamma'=\Gamma_{Q',W'}$ for the associated Ginzburg dgas.  Write $T_i$ for the left $\Gamma'$-module associated to the vertex $i$ of $Q'$ and $G_i$ for the associated autoequivalence $F_{T_i}$ of $\Dfd(\Gamma')$.  In this section we will study how the spherical twists $F_i:\Dfd(\Gamma)\arr\sim\Dfd(\Gamma)$ interact with the mutation functors $\mu_k:\Dfd(\Gamma')\arr\sim\Dfd(\Gamma)$.
Our key tools will be Proposition \ref{prop:conj-sph} and the results on the images of the simple modules under the mutation functors \cite[Lemma 3.12(a)]{ky}, which we will describe below.

If $A$ is a dga and $M,N\in\Dfd(A)$, we have a natural map
$$M\otimes_\F\Hom_{\Dfd(A)}(M,N)\to N$$
in $\Dfd(A)$ given by evaluation.  For any graded vector space $V$, we have biadjoint functors $-\otimes_FV$ and $-\otimes_F V^*$, and these respect the left $A$-module structure, so we also obtain a natural map
$$M\to N\otimes_\F\Hom_{\Dfd(A)}(M,N)^*$$
in $\Dfd(A)$.  Now let $L,N\in A\Mod$.  The \emph{universal extension of $N$ by $L$} is the cone of the natural map
$$N[-1]\to L\otimes_\F\Hom_{\Dfd(A)}(N[-1],L)^*$$
and the \emph{universal coextension of $L$ by $N$} is the cone of the natural map
$$N[-1]\otimes_\F\Hom_{\Dfd(A)}(N[-1],L)\to L.$$

The following result is contained in the proof of \cite[Lemma 3.12(a)]{ky}:
\begin{lem}[Keller-Yang]\label{lem:kysimples}
We have $\mu_k(T_k)\cong S_k[1]$ and, for $i\neq k$, $\mu_k(T_i)$ is the universal extension of $S_i$ by $S_k$.
\end{lem}

The following result should be compared to Definition \ref{def:mut-rule}.
\begin{prop}\label{prop:descr_fmuinv}
If $Q$ has no double arrows then we have a natural isomorphism of functors
\begin{equation*}
F_{\mu_k^{-1}(S_i)} \cong \left\{
\begin{array}{cl}
G_kG_iG_k^{-1} & \text{if } i\to k\text{ in }Q;\\
G_i & \text{otherwise}.
\end{array} \right.
\end{equation*}

\begin{pf}
We first use Lemma \ref{lem:kysimples} to calculate the images of the simple $\Gamma'$-modules under the inverse mutation functor $\mu^{-1}$, where $\mu=\mu_k$.  We know that $\mu(T_k)\cong S_k[1]$, so $\mu^{-1}(S_k)\cong T_k[-1]$.
By assumption, there is at most one arrow between any two vertices in $Q$.
If $i\neq k$ and there is no arrow $i\to k$ in $Q$
then, by Lemma \ref{lem:kysimples}, $\Hom_{\Dfd(\Gamma)}(S_i[-1],S_k)=0$ and so
$\mu(T_i)\cong\cone(S_i[-1]\to0)\cong S_i$, thus $\mu^{-1}(S_i)\cong T_i$.

If $i\neq k$ and there is an arrow
$i\to k$ in $Q$
then $\Hom_{\Dfd(\Gamma)}(S_i[-1],S_k)$ is $1$-dimensional and so 
$\mu(T_i)\cong\cone(S_i[-1]\to S_k)$, with the nonzero map determined up to a scalar.  We can then use Lemma \ref{lem:coneFf-comm} to calculate $\mu(\cone(T_k[-1]\to T_i)$: this is
$\cone(\mu(T_k)[-1]\to\mu(T_i))$ where, as $\mu$ is an equivalence, the map must again be nonzero and determined up to scalar.  We know that $\mu(T_k)[-1]\cong S_k$ and $\mu(T_i)$ is $S_i\oplus S_k$ with appropriate differential.  One can check that the injection $S_k\into S_i\oplus S_k$ respects the differentials, and so this must be our nonzero map.  This is quasi-isomorphic to the map $0\to S_i$, and so $\mu(\cone(T_k[-1]\to T_i))\cong S_i$ and hence $\mu^{-1}(S_i)\cong \cone(T_k[-1]\to T_i)$.  Note that this is the universal coextension of $T_i$ by $T_k$.

Now we check that the formula holds.  If $i=k$ then $F_{\mu^{-1}(S_i)}=F_{T_i[1]}$, and as the shift functor on $\Dfd(\Gamma')$ is naturally isomorphic to $\Gamma'[1]\otimes_{\Gamma'}-$ we see that, by Proposition \ref{prop:conj-sph}, $F_{T_i[1]}\cong[1]\circ G_i\circ [-1]\cong G_i$.
If $i\neq k$ and there is no arrow
$i\to k$ in $Q$
then $\mu^{-1}(S_i)\cong T_i$ so $F_{\mu^{-1}(T_i)}=G_i$.

Finally, suppose $i\neq k$ and there is an arrow
$i\to k$ in $Q$.  
As mutation at $k$ reverses all arrows incident with $k$, and can never change the number of arrows incident with $k$, there must be exactly one arrow $k\to i$ in $Q'$.
We first calculate $G_k(T_i)$: this is
$$\cone(\p T_k\otimes_\F\Hom_{\Gamma'}(\p T_k,T_i)\to T_i).$$
As $\Hom_{\Gamma'}(\p T_k,T_i)$ is a differential graded $\F$-module, it is quasi-isomorphic to its homology, which is the direct sum $\bigoplus\Hom_{K(\Gamma')}(T_k,T_i[n])$ with $\Hom_{K(\Gamma')}(T_k,T_i[n])\cong \Hom_{\Dfd(\Gamma')}(T_k,T_i[n])$ in degree $n$.  So by Lemma \ref{lem:kysimples} the homology is only nonzero in degree $1$, where it is $1$-dimensional, and thus
$$G_k(T_i)\cong
\cone( T_k\otimes_\F\F[-1]\to T_i).$$
So we see that $\mu^{-1}(S_i)\cong G_k(T_i)$ and thus, using Proposition \ref{prop:conj-sph} again, $F_{\mu^{-1}(S_i)} \cong G_kG_iG_k^{-1} $.
\end{pf}
\end{prop}

We are now able to show that our braid groups $B_Q$ act via spherical twists on the category $\Dfd(\Gamma)$.
\begin{thm}\label{thm:mutation-action}
Let $(Q,W)$ be a mutation-Dynkin quiver with potential of type $ADE$.  Then we have a group homomorphism
\begin{align*}
B_Q &\to \Aut\Dfd(\Gamma_{Q,W})\\
s_i &\mapsto F_i
\end{align*}
sending the group generator associated to the vertex $i\in I$ to the spherical twist at the simple $\Gamma_{Q,W}$-module $S_i$.

\begin{pf}
As $(Q,W)$ is mutation-Dynkin, it is obtained by mutating a quiver with potential $(Q'',0)$ finitely many times, where $Q''$ is a Dynkin quiver.  Then we have a group homomorphism $B_{Q''}\to \Aut\Dfd(\Gamma_{Q'',0})$ by Remark \ref{r:artin}, Proposition \ref{prop:braid-relns}, and Lemma \ref{lem:dimext}.  This gives the base case of an inductive argument, so we need to show that if the spherical twists $F_i$ on $\Gamma=\Gamma_{Q,W}$ satisfy the relations of $B_{Q}$ for a mutation-Dynkin quiver with potential $(Q,W)$ then the spherical twists $G_i$ on $\Gamma'=\Gamma_{Q',W'}$ satisfy the relations of $B_{Q'}$.

Assume the functors $F_i:\Dfd(\Gamma)\to \Dfd(\Gamma)$ satisfy the relations of $B_{Q}$ and let $\mu=\mu_k:\Dfd(\Gamma')\to\Dfd(\Gamma)$ be the Keller-Yang derived equivalence.  Then the functors $\mu^{-1}\circ F_i\circ \mu:\Dfd(\Gamma')\to\Dfd(\Gamma')$ also satisfy the relations of $B_{Q}$.  By Proposition \ref{prop:conj-sph} we have $\mu^{-1}\circ F_i\circ \mu\cong F_{\mu^{-1}(S_i)}$, i.e., the following diagram commutes:
$$\xymatrix{
\Dfd(\Gamma')\ar[r]^\mu\ar[d]^{F_{\mu^{-1}(S_i)}} &\Dfd(\Gamma)\ar[d]^{F_{S_i}} \\
\Dfd(\Gamma')\ar[r]^\mu &\Dfd(\Gamma)
}$$
So we have a group homomorphism $\rho:B_Q \arr\rho \Aut\Dfd(\Gamma')$ sending $s_i$ to $F_{\mu^{-1}(S_i)}$.
As, by Proposition \ref{prop:Dynkinquiv}, $Q$ has no double arrows, we can use Proposition \ref{prop:descr_fmuinv} to write $\rho$ as
\begin{align*}
 B_Q \arr\rho &\Aut\Dfd(\Gamma') \\
 s_i \mapsto  &
   \left\{
\begin{array}{cl}
G_kG_iG_k^{-1} & \text{if } i\to k\text{ in }Q;\\
G_i & \text{otherwise}.
\end{array} \right.
\end{align*}
Precomposing this with the group isomorphism $\varphi_k^{-1}:B_Q'\arr\sim B_Q$ of Remark \ref{rmk:phi-inv}, we obtain the group homomorphism
$$\xymatrix @R=10pt {
B_{Q'}\ar[r]^{\varphi_k^{-1}} &B_Q\ar[r]^\rho &\Aut\Dfd(\Gamma')\\
t_i\ar@{|->}[r] &
{\left\{
\begin{array}{cl}
s_k^{-1}s_is_k & \text{if } i\to k\text{ in }Q;\\
s_i & \text{otherwise};
\end{array} \right\}}
\ar@{|->}[r] &
{\left\{
\begin{array}{cl}
G_k^{-1}G_kG_iG_k^{-1}G_k & \text{if } i\to k\text{ in }Q;\\
G_i & \text{otherwise};
\end{array} \right\}}
\cong G_i
}$$
as required.
\end{pf}
\end{thm}

\begin{rmk}\label{rmk:faithful}
Known results on the faithfulness of braid group actions can be transferred to our setting.  Suppose $Q''$ is an orientation of an $ADE$ graph and the usual action $B_{Q''}\to\Aut\Dfd(\Gamma_{Q'',0})$ is faithful.  From the proof of Theorem \ref{thm:mutation-action} we see that our actions of $B_Q$ where $Q$ is of mutation type $ADE$ are just built by precomposing group isomorphisms with the action of $B_{Q''}$, and so these are also faithful under this assumption.

It was shown by Seidel and Thomas \cite[Theorem 2.18]{st}, building on work of Khovanov and Seidel \cite{ks}, that given a collection of $d$-spherical objects, with $d\geq 2$, in a type $A_n$-configuration the action of the braid group by spherical twists is faithful. 
Thus the actions of Theorem \ref{thm:mutation-action} are faithful in mutation type $A$.  The faithfulness result of Seidel and Thomas was
 extended to all collections of $2$-spherical objects in type $ADE$ configurations by Brav and Thomas \cite{bt}, using the Garside structure of the braid monoid, 
but it is not immediately clear how to generalize their argument to the $3$-Calabi-Yau situation.
\end{rmk}

\begin{rmk}
Although we have shown that our braid groups of mutation-Dynkin quivers can be realized categorically, this is not a categorification of our earlier results because we cannot decategorify (see, for example, \cite{bd}): we cannot recover Theorem \ref{thm:groupiso} from Theorem \ref{thm:mutation-action} because we use Theorem \ref{thm:groupiso} to prove Theorem \ref{thm:mutation-action}.  The problem is that, for an arbitrary mutation-Dynkin quiver with potential $(Q,W)$, we do not in advance know the relations satisfied by the spherical twist functors $F_i$.  This question will be addressed in a forthcoming paper \cite{gra9}.
\end{rmk}

\begin{rmk} \label{r:singletwist}
The arguments of \cite{gra2} generalize to show that, if vertices $i$ and $j$ of $Q$ are joined by an arrow, then $F_iF_jF_i\cong F_jF_iF_j$ can be realized as a single periodic twist.  Similarly, one can show that if $i\to j\to k\to i$ is a $3$-cycle in $Q$ then $F_1F_2F_3F_1\cong F_2F_3F_1F_2 \cong F_3F_1F_2F_3$ can be realized as a single periodic twist.  This exhausts the possibilities in type $A$; type $D$ will be studied further in \cite{gra9}.
\end{rmk}

\textbf{Funding:} \\
This work was supported by the Engineering and Physical Sciences Research Council 
[grant number EP/G007497/1] and the Mathematical Sciences Research Institute,
Berkeley.

\textbf{Acknowledgements:} \\
The second author would like to thank K. Baur and A. King for useful discussions concerning orbifold interpretations of cluster algebras. Both authors would like
to thank the referee, and I. Webster, for helpful comments on an earlier version of
the paper. Some of the work for this paper was done while the second author was
Guest Professor in 2014 at the Department of Mathematical Sciences, Norwegian University of Science and Technology, Trondheim, N-7491, Norway. \\


School of Mathematics, University of East Anglia, Norwich Research Park,
Norwich, NR4 7TJ, UK. \\
\href{mailto:j.grant@uea.ac.uk}{\nolinkurl{j.grant@uea.ac.uk}}

School of Mathematics, University of Leeds, Leeds LS2 9JT, U.K.\\
\href{mailto:marsh@maths.leeds.ac.uk}{\nolinkurl{marsh@maths.leeds.ac.uk}}

\end{document}